\documentclass[12pt]{amsart}
\usepackage{amssymb}
\usepackage{amsthm}
\usepackage{amscd}
\usepackage{mathrsfs}

\usepackage{pifont}

\usepackage{latexsym}
\usepackage{amsfonts}

\usepackage{amsmath}
\usepackage{epic,eepic}
\usepackage{graphicx}

\theoremstyle{plain}
\newtheorem{lemma}{Lemma}[subsection]
\newtheorem{prop}[lemma]{Proposition}
\newtheorem{thm}[lemma]{Theorem}
\newtheorem{cor}[lemma]{Corollary}
\newtheorem{aplemma}{Lemma~A.\hspace{-1.5mm}}
\newtheorem{approp}{Proposition~A.\hspace{-1.5mm}}
\newtheorem{apthm}{Theorem~A.\hspace{-1.5mm}}
\newtheorem{apcor}{Corollary~A.\hspace{-1.5mm}}
\newtheorem{intthm}{Theorem}

\usepackage[all]{xy}

\newtheorem{conj}[lemma]{Conjecture}

\theoremstyle{definition}

\newtheorem{rema}[lemma]{Remark}

\newtheorem{remb}{Remark}

\newtheorem{defi}[lemma]{Definition}
\newtheorem{exa}[lemma]{Example}
\newtheorem{aprem}{Remark~A.\hspace{-1.5mm}}
\newtheorem{apdefi}{Definition~A.\hspace{-1.5mm}}
\newcommand{\bde}{\begin{defi}}
\newcommand{\ede}{\end{defi}\vspace{1mm}}
\newcommand{\ble}{\begin{lemma}}
\newcommand{\ele}{\end{lemma}}
\newcommand{\bpr}{\begin{prop}}
\newcommand{\epr}{\end{prop}}
\newcommand{\bt}{\begin{thm}}
\newcommand{\et}{\end{thm}}
\newcommand{\bco}{\begin{cor}}
\newcommand{\eco}{\end{cor}}
\newcommand{\bre}{\begin{rema}}
\newcommand{\ere}{\end{rema}}
\newcommand{\brea}{\begin{rema}}
\newcommand{\erea}{\end{rema}\vspace{1mm}}
\newcommand{\breb}{\begin{remb}}
\newcommand{\ereb}{\end{remb}\vspace{1mm}}
\newcommand{\bex}{\begin{exa}}
\newcommand{\eex}{\end{exa}}
\newcommand{\bpf}{\begin{proof}}
\newcommand{\epf}{\end{proof}\vspace{1mm}}

\newcommand{\bade}{\begin{apdefi}}
\newcommand{\eade}{\end{apdefi}}
\newcommand{\bale}{\begin{aplemma}}
\newcommand{\eale}{\end{aplemma}}
\newcommand{\bapr}{\begin{approp}}
\newcommand{\eapr}{\end{approp}}
\newcommand{\bat}{\begin{apthm}}
\newcommand{\eat}{\end{apthm}}
\newcommand{\baco}{\begin{apcor}}
\newcommand{\eaco}{\end{apcor}}
\newcommand{\bare}{\begin{aprem}}
\newcommand{\eare}{\end{aprem}}

\newcommand{\be}{\begin{enumerate}}
\newcommand{\ee}{\end{enumerate}}
\newcommand{\bcd}{\[\begin{CD}}
\newcommand{\ecd}{\end{CD}\]}
\newcommand{\bit}{\begin{itemize}}
\newcommand{\eit}{\end{itemize}}
\newcommand{\bq}{\begin{quote}}
\newcommand{\eq}{\end{quote}}
\newcommand{\ba}{\begin{array}}
\newcommand{\ea}{\end{array}}

\newcommand{\mcB}{\mathcal{B}}
\newcommand{\mcC}{\mathcal{C}}

\newcommand{\mcE}{\mathcal{E}}
\newcommand{\mcF}{\mathcal{F}}
\newcommand{\mcG}{\mathcal{G}}
\newcommand{\mcH}{\mathcal{H}}
\newcommand{\mcI}{\mathcal{I}}

\newcommand{\mcK}{\mathcal{K}}
\newcommand{\mcL}{\mathcal{L}}
\newcommand{\mcM}{\mathcal{M}}
\newcommand{\mcN}{\mathcal{N}}
\newcommand{\mcO}{\mathcal{O}}
\newcommand{\mcP}{\mathcal{P}}

\newcommand{\mcR}{\mathcal{R}}
\newcommand{\mcS}{\mathcal{S}}
\newcommand{\mcT}{\mathcal{T}}
\newcommand{\mcU}{\mathcal{U}}
\newcommand{\mcV}{\mathcal{V}}

\newcommand{\mbB}{\mathbb{B}}
\newcommand{\mbC}{\mathbb{C}}

\newcommand{\mbE}{\mathbb{E}}
\newcommand{\mbF}{\mathbb{F}}
\newcommand{\mbG}{\mathbb{G}}

\newcommand{\mbP}{\mathbb{P}}
\newcommand{\mbQ}{\mathbb{Q}}
\newcommand{\mbR}{\mathbb{R}}

\newcommand{\mbU}{\mathbb{U}}

\newcommand{\mbX}{\mathbb{X}}

\newcommand{\mbZ}{\mathbb{Z}}

\newcommand{\msE}{\mathscr{E}}

\newcommand{\msP}{\mathscr{P}}

\newcommand{\msX}{\mathscr{X}}

\newcommand{\mfC}{\mathfrak{C}}

\newcommand{\mfM}{\mathfrak{M}}
\newcommand{\mfN}{\mathfrak{N}}
\newcommand{\mfO}{\mathfrak{O}}
\newcommand{\mfP}{\mathfrak{P}}

\newcommand{\mfR}{\mathfrak{R}}
\newcommand{\mfS}{\mathfrak{S}}
\newcommand{\mfT}{\mathfrak{T}}
\newcommand{\mfU}{\mathfrak{U}}
\newcommand{\mfV}{\mathfrak{V}}
\newcommand{\mfW}{\mathfrak{W}}
\newcommand{\mfX}{\mathfrak{X}}
\newcommand{\mfY}{\mathfrak{Y}}

\newcommand{\mfa}{\mathfrak{a}}
\newcommand{\mfb}{\mathfrak{b}}
\newcommand{\mfc}{\mathfrak{c}}
\newcommand{\mfd}{\mathfrak{d}}
\newcommand{\mfe}{\mathfrak{e}}

\newcommand{\mfg}{\mathfrak{g}}
\newcommand{\mfh}{\mathfrak{h}}

\newcommand{\mfl}{\mathfrak{l}}

\newcommand{\mfn}{\mathfrak{n}}
\newcommand{\mfo}{\mathfrak{o}}
\newcommand{\mfp}{\mathfrak{p}}

\newcommand{\mfs}{\mathfrak{s}}
\newcommand{\mft}{\mathfrak{t}}

\newcommand{\migi}{\rightarrow}
\newcommand{\longmigi}{\longrightarrow}

\newcommand{\isom}{\stackrel{\sim}{\migi}}

\newcommand{\migiincl}{\hookrightarrow}

\newcommand{\migisurj}{\twoheadrightarrow}

\newcommand{\mr}{\mathrm}
\newcommand{\hidden}[1]{\,}
\newcommand{\SSP}{\vspace{3mm}}
\newcommand{\LSP}{\vspace{5mm}}

\usepackage{color}

\def\pmoutlinefnt#1{\setbox0=\hbox{#1}%
   \setbox1=\hbox{\kern-.020em\copy0\kern-\wd0\kern.020em\copy0%
   \kern-\wd0\kern.020em\copy0}$
   \copy1\kern-\wd1\raise.020em\copy1\kern-\wd1\raise-.020em\copy1%
   \color[rgb]{1,1,1}\kern-\wd0\kern-.020em\box0
$}

\def\pmshadowfont#1{\setbox0=\hbox{#1}%
   \setbox1=\hbox{\kern-.020em\copy0\kern-\wd0\kern.020em\copy0%
   \kern-\wd0\kern.020em\copy0}$
   \copy1\kern-\wd1\raise.020em\copy1\kern-\wd1\raise-.020em\copy1%
   \color[rgb]{1,1,1}\kern-\wd0\kern-.065em\raise.015em\copy0
$}

\newcommand{\vin}{\rotatebox{90}{$\in$}}

\begin{document}

\title[Topological quantum field theory for dormant opers]{
Topological quantum field theory for dormant opers}
\author{Yasuhiro Wakabayashi}
\date{}
\markboth{Yasuhiro Wakabayashi}{}
\maketitle	
\footnotetext{Y. Wakabayashi: Department of Mathematics, Tokyo Institute of Technology, 2-12-1 Ookayama, Meguro-ku, Tokyo 152-8551, JAPAN;}
\footnotetext{e-mail: {\tt wkbysh@math.titech.ac.jp};}
\begin{abstract}
The purpose of the present paper is to develop the enumerative geometry of dormant $G$-opers for a semisimple algebraic group $G$. In the present paper, we construct a compact moduli stack admitting a perfect obstruction theory by introducing the notion of a dormant faithful twisted $G$-oper (or a ``{\it $G$-do'per}\," for short). The resulting virtual fundamental class induces a semisimple $2$d TQFT (= $2$-dimensional topological quantum field theory) counting the number of $G$-do'pers. This $2$d TQFT gives an analogue of the Witten-Kontsevich theorem describing the intersection numbers of psi classes on the moduli stack of $G$-do'pers.
\end{abstract}
\tableofcontents 

\section*{Introduction} \SSP

The purpose of the present paper is to develop    the {\it enumerative geometry of the moduli
       of dormant $G$-opers} (i.e., $G$-opers with vanishing $p$-curvature) 
  for a semisimple algebraic group $G$ in characteristic $p>0$.
The formulations and background knowledge of dormant $G$-opers used    in the present  paper  under the assumption that $G$ is of adjoint type  were discussed in the author's paper  ~\cite{Wak5}.
In the  present paper, we generalize the previous work and  construct a compact moduli stack admitting a perfect obstruction theory  by introducing  the notion of a {\it dormant faithful twisted $G$-oper} (or  a ``{\it $G$-do'per}\," for short) defined on a stacky log curve. 
The resulting virtual  fundamental class (in the cases of some classical types $G$) induces 
 a semisimple $2$d TQFT (= $2$-dimensional topological quantum field theory), which 
  arises 
from the nature of algebraic geometry in 
 positive characteristic unlike many of the other examples of  TQFTs constructed  in 
    geometry.
   This result may be thought of as an improvement  of ~\cite[Proposition 7.33]{Wak5}.
In particular, 
an explicit description of (the characters of) the corresponding Frobenius algebra allows us 
to  perform a computation  for the counting problem of $G$-do'pers, i.e., the Verlinde formula for $G$-do'pers.
This  $2$d TQFT   gives also   an analogue of 
the  Witten-Kontsevich theorem,   describing    the  intersection numbers of psi classes  on the moduli stack of $G$-do'pers.

Thus, the results proved in the present paper show that
the moduli stack of $G$-do'pers 
 has many features similar to certain spaces  dealt with in enumerative geometry related to high energy physics, say,   the moduli stack of $r$-spin curves (denoted usually by $\mcB_{g, r} (\mbG_m, \omega_{\mr{log}}^{1/r})$  or $\overline{\mfM}_{g,r}^{1/r}$) and the moduli stack of stable maps into a suitable variety $V$ (denoted by $\overline{\mfM}_{g,r} (V)$).  
In the rest of this Introduction, we shall provide   more detailed discussions, including the content of the present paper.

\LSP
\subsection*{0.1} \label{S01}
Recall that a  {\it dormant $G$-oper} 
 is,  roughly speaking, a
$G$-bundle
   on an algebraic curve in characteristic $p >0$  equipped with a connection satisfying certain conditions, including the condition that its  $p$-curvature vanishes identically.
Various properties of dormant $\mr{PGL}_2$-opers and their moduli 
  were  discussed
  by S. Mochizuki (cf.  ~\cite{Mzk1}, ~\cite{Mzk2}) in the context of  {\it $p$-adic Teichm\"{u}ller theory}.
If 
$G = \mr{PGL}_n$ or $\mr{SL}_n$ for a general $n$
 (but the underlying curve  is assumed to be unpointed and  smooth  over an algebraically closed field), then
these objects
   has been studied  by K. Joshi, S. Ramanan, E. Z. Xia, J. K. Yu, C. Pauly, T. H. Chen, X. Zhu et al.
(cf.  ~\cite{CZ}, ~\cite{JRXY}, ~\cite{JP}, and  ~\cite{Jo14}).

As discussed in these  references,   
dormant $G$-opers and  their moduli
have  diverse aspects and 
occur naturally in mathematics.
A detailed understanding  of them  from  the viewpoint of enumerative geometry
 will be  applied to the counting  problems  for various objects (e.g., lattice points inside a  convex rational polytope, edge-colorings of a trivalent graph, and Frobenius-destabilized bundles, etc.) appearing in  some areas of mathematics linked to 
    the theory of dormant opers, see ~\cite{JP}, ~\cite{LO},  ~\cite{Os1}, and ~\cite{Wak2}.

\LSP
\subsection*{0.2} \label{S011}
As an example, let us explain 
 one aspect of 
 the enumerative geometry of dormant $G$-opers concerning
  the algebraic-solution problem of  linear differential equations in positive characteristic.

Let
 $X$ be 
a geometrically connected,  proper, and  smooth curve over  a perfect  field  $k$ of characteristic $p$ with function field $K$.
 Consider a monic  linear  ordinary differential operator  $D$  of order $n>1$  with regular singularities  defined on $X$; it  may be expressed locally  as follows:
\begin{align} \label{ee650}
D := \frac{d^n}{d x^n}  +q_1 \frac{d^{n-1}}{d x^{n-1}} + \cdots + q_{n-1} \frac{d}{dx} + q_n,  
\end{align}
where  $q_1, \cdots, q_n \in K$ and $x$ denotes a local coordinate in $X$.
Since the  $p$-powers $K^p$ of elements of $K$ coincides with  the constant field in $K$ (i.e., the kernel of the universal derivation $d : K \migi \Omega_{K/k}$), the set of solutions to the equation $Dy=0$ in $K$  forms a $K^p$-vector space.

We shall say that the differential equation $D y=0$ has {\bf a full  set of  solutions}
if it has $n$ solutions  in $K$ linearly independent over $K^p$.
The study of differential equations having many (algebraic)  solutions was originally considered
in the complex case,  tackled and developed  since  the 1870s by many mathematicians: H. A. Schwarz (for the hypergeometric equations), L. I. Fuchs, P. Gordan, and C. F.  Klein (for the second order equations), C. Jordan (for the $n$-th order) et al.
We are interested in the  positive characteristic analogue   of this
traditional study,
 and, in particular, want to know {\it how many  differential equations $D y = 0$ (associated with  $D$ as in (\ref{ee650})) have  a full set of solutions}.

Let us consider a special case, i.e., the case  
 of monic, linear, and   second order differential operators on   the projective line $\mbP_{k} := \mr{Proj} (k [s, t])$ having at most three regular singular points.
The classical theory of Riemann schemes shows that 
any such operator may be transformed (e.g., via pull-back by  an automorphism of $\mbP_{k}$) into
a {\it Gauss' hypergeometric differential operator}
\begin{align}
D_{a, b, c} := \frac{d^2}{d x^2} + \left(\frac{c}{x} + \frac{1 -c + a + b}{x-1} \right)\frac{d}{dx} +\frac{ab}{x (x-1)} 
\end{align}
determined by some  triple  $(a, b,c) \in k^{\times 3}$, where $x := s/t$ and $k^{\times 3}$ denotes the product of $3$ copies of $k$.

We shall  use the notation  $\widetilde{(-)}$  to denotes the inverse of the  bijective restriction  $\{1, \cdots, p \} \isom \mbF_p := \mbZ/p\mbZ$  of the natural quotient $\mbZ \migisurj \mbF_p$.
If the equation $D_{a, b, c} y=0$ has a full set of  solutions, then
the set $\left\{y_{a,b,c} (x), x^{1-\widetilde{c}} y_{a-c+1, b-c+1, 2-c} (x) \right\}$
$\subseteq k (x)$ forms a basis of the solutions.
Here,  $y_{a, b, c} (x)$ denotes  a polynomial of $x$ defined by the following  truncated hypergeometric series
\begin{align}
y_{a,b,c} (x) := & \  1 + \frac{a \cdot b}{1 \cdot c}x + \frac{a \cdot  (a+1) \cdot b \cdot  (b+1)}{1 \cdot 2 \cdot c \cdot (c +1)}x^2 \\
& +  \frac{a \cdot  (a+1) \cdot (a+2) \cdot b \cdot  (b+1) \cdot (b+2)}{1 \cdot 2 \cdot 3\cdot c \cdot (c +1) \cdot (c+2)}x^3 + \cdots,  \notag
\end{align}
where we stop the series as soon as the numerator vanishes.

According to ~\cite[\S\,1.6]{Ihara1} (or ~\cite[\S\,6.4]{Katz2}), 
the equation $D_{a, b, c} y=0$ has a full set of solutions if and only if $(a,b,c)$ lies in $ \mbF_p^{\times 3}$
   and  either $\widetilde{b} \geq  \widetilde{c} > \widetilde{a}$ or $\widetilde{a} \geq  \widetilde{c} > \widetilde{b}$ is satisfied.
   In particular, 
after a straightforward calculation, we see that {\it there exists precisely    $\frac{p^3-p}{3}$ hypergeometric equations having   a full set of solutions}.

\LSP
\subsection*{0.3} \label{S021}
Here, we shall recall the relationship between dormant opers and  differential operators
  in terms of connections on 
  vector bundles.
To each differential operator  $D$ as in (\ref{ee650}), one can  associate, in a well-known manner, 
a  connection on a vector bundle
 expressed locally as follows:
\begin{align} \label{ee300}
\nabla = \frac{d}{dx} -   
\begin{pmatrix} 
-q_1 & -q_2  & -q_3 & \cdots  & -q_{n-1} & -q_n \\
 1 & 0 & 0 & \cdots    & 0 & 0 \\
  0 & 1 & 0&  \cdots   & 0 & 0 \\
   0 & 0 & 1&  \cdots  & 0 & 0 \\ 
  \vdots    & \vdots  & \vdots & \ddots  &  \vdots & \vdots \\ 
       0 & 0  & 0 & \cdots    & 1 & 0
     \end{pmatrix}.
\end{align}
The assignment $y \mapsto {^t(} \frac{d^{n-1}y}{dx^{n-1}}, \cdots, \frac{dy}{dx},  y)$ gives a bijective  correspondence between  
the solutions of the equation $D y=0$ 
and the horizontal  sections of this vector bundle with respect to $\nabla$.

A vector bundle equipped with 
 a connection of the form (\ref{ee300})
 provides an oper.
Indeed,  any $\mr{PGL}_n$-oper can be represented, via a gauge transformation, 
by  such a connection 
  (cf.   ~\cite[Theorem D]{Wak5}).
Under a certain assumption on the subprincipal  symbol of $D$ (i.e., the assumption that the $1$-st order differential equation  determined from the subprincipal symbol using  the manner discussed in  ~\cite[Remark 4.30]{Wak5} has a nonzero solution), 
 the differential equation $D y =0$ associated with  $D$ as in (\ref{ee650}) has a full set of solutions if  and only if 
 the associated $\mr{PGL}_n$-oper is dormant  (cf. ~\cite[(6.0.5), Proposition]{Katz2}).

We  go back to the hypergeometric case.
Let   $\msP_{k} := (\mbP_{k}, \{ [0], [1], [\infty ] \})$ be  the projective line with three marked points determined by $0$, $1$, and $\infty$.
For a triple $(a, b, c) \in k^{\times 3}$,
denote by $\msE_{a,b,c}^\spadesuit$ the dormant $\mr{PGL}_2$-oper  on  $\msP_{k}$ induced by $D_{a,b,c}$.
If we set $\mr{Ex}_{a, b, c} := (1-c, c-a-b, b-a)$ (i.e., the exponent differences at $0$, $1$, and $\infty$),
then the isomorphism class of $\msE_{a,b,c}^\spadesuit$ is determined by  this data regarded as a triple of elements in $k/\{ \pm1\}$ (:= the quotient set of $k$ by the equivalence relation generated by $v \sim - v$ for every $v \in k$).
On the other hand, 
it is verified  that
\begin{align}
\mcE^\spadesuit_{a, b, c} \cong  \mcE^\spadesuit_{a', b', c'} \ \ \Longleftrightarrow \ \  \mr{Ex}_{a, b, c} = \mr{Ex}_{a', b', c'} \ \text{in} \ (k/\{ \pm1\})^{\times 3}.
\end{align}
 Hence, 
since $\mr{Ex}_{a, b, c} \neq (0, 0, 0)$ if $\msE_{a,b,c}^\spadesuit$ is dormant, 
the assignment $D_{a, b, c} \mapsto \msE^\spadesuit_{a, b, c}$ determines  a $2^3$-to-$1$ correspondence
\begin{align} \label{QG2010}
\begin{pmatrix}
\text{the set of hypergeometric} \\
\text{operators $D_{a,b,c}$ such that the equation}\\
\text{$D_{a,b,c} y =0$ has a full set of solutions}
\end{pmatrix}
\stackrel{2^3 : 1}{\longleftrightarrow}
\begin{pmatrix}
\text{the set of isomorphism classes} \\
\text{of dormant $\mr{PGL}_2$-opers on $\msP_{k}$}
\end{pmatrix}.
\end{align}
In particular,  the italicized assertion at the end of the previous subsection implies that 
the number of (isomorphism classes of) dormant $\mr{PGL}_2$-opers on $\msP_{k}$ is exactly equal to   $\frac{p^3-p}{24}$.

The same computations  were  verified 
by S. Mochizuki (cf. 
~\cite[Chap.\,V, Corollary 3.7]{Mzk2}),
H. Lange-C. Pauly (cf. ~\cite[Theorem 2]{LP}), and B. Osserman (cf. ~\cite[Theorem 1.2]{Os4}) by  applying different methods.
At any rate, the simplest case of the counting problem for  dormant opers 
was completely resolved. 
When one sets about trying to compute  the number of differential operators or the corresponding opers in more general cases (i.e., the case of general $n$ and $X$),
 the factorization property of their  moduli stacks  discussed in  
 the present paper (or ~\cite[Theorem F]{Wak5}) will provide an effective way to do this.

\LSP
\subsection*{0.4} \label{S101}
In the rest of this Introduction, we  describe  briefly the results obtained in the present paper.
 
 Let $G$ be a split semisimple algebraic group over a perfect field $k$ of characteristic $p$.
Then,  we define a {\it dormant faithful twisted $G$-oper} (or  a {\it $G$-do'per} for short) 
as  a   $G$-oper on a stacky log curve with vanishing $p$-curvature satisfying a certain  representability condition  (cf. Definitions \ref{Dy0351}, \ref{Dy035}, and \ref{DD035}).  
Under the assumption that   $G$ is of adjoint type,  the notion of a $G$-do'per 
coincides with   the classical notion of a  $\mfg$-oper, where $\mfg$ denotes the Lie algebra of $G$,  in the sense of   ~\cite[Definitions 2.1 and 3.15]{Wak5} (cf. Remark \ref{ppp037}).

Given a pair of nonnegative integers $(g, r)$ with $2g-2+r>0$.
we denote  the usual  moduli stack classifying $r$-pointed stable curves of genus $g$ over $k$ by $\overline{\mfM}_{g,r}$.
Then, 
one obtains the category  fibered in groupoids
\begin{align}
\mfO \mfp_{G, g,r}^{^\mr{Zzz...}}
\end{align}
 (cf. (\ref{ee215}))   classifying $G$-do'pers on stacky log curves whose  coarse moduli spaces are  classified by $\overline{\mfM}_{g,r}$.
Here, 
let us consider two additional conditions  $(*)_{G}$,  $(**)_{G}$ on $G$   described as follows:
\begin{itemize}
\item[]
\begin{itemize} 
 \item[$(*)_{G}$ :] 
 If  $h$ denotes  the integer $h_{\mbG}$ defined in ~\cite[Eq.\,(104)]{Wak5}, where ``$\mbG$"  is taken to be the adjoint group of $G$, then the inequality   $p > 2 h$ holds;
\vspace{1mm}
 \item[$(**)_{G}$ :] 
 $G$ is 
 of classical type $A_n$ (with $2n < p-2$ or $(n, p)=(1, 3)$), $B_l$ (with $4l < p-2$), or $C_m$ (with $4m < p$).
\vspace{1.5mm}
\end{itemize}
\end{itemize}
Notice that $(**)_{G}$ implies $(*)_{G}$.
The following assertion is the first main result of the present paper, which 
concerns 
 the structure of the moduli  space $\mfO \mfp_{G, g,r}^{^\mr{Zzz...}}$;
 this result   generalizes ~\cite[Theorems C and G]{Wak5}.

\SSP
\begin{intthm} [cf.  Theorems \ref{P05} and \ref{P01}] \label{Tqqq}
 \begin{itemize}
 \item[(i)]
 Assume that  $G$ satisfies the condition $(*)_{G}$.
 Then, $\mfO \mfp_{G, g,r}^{^\mr{Zzz...}}$ may be represented by a nonempty  proper  Deligne-Mumford stack over $k$ which is finite over  $\overline{\mfM}_{g,r}$ and has an irreducible component that dominates $\overline{\mfM}_{g,r}$.
 Moreover,  $\mfO \mfp_{G, g,r}^{^\mr{Zzz...}}$ admits  a perfect obstruction theory, and hence, 
  has a virtual fundamental class $[\mfO \mfp_{G, g,r}^{^\mr{Zzz...}}]^\mr{vir}$.
 \item[(ii)]
 Assume further that $G$ satisfies the condition  $(**)_{G}$.
 Then,  $\mfO \mfp_{G, g,r}^{^\mr{Zzz...}}$  is generically \'{e}tale  over $\overline{\mfM}_{g,r}$ (i.e., any irreducible component of  $\mfO \mfp_{G, g,r}^{^\mr{Zzz...}}$  dominating  $\overline{\mfM}_{g,r}$ has a dense open substack which is \'{e}tale  over $\overline{\mfM}_{g,r}$), and moreover, has generic stabilizer isomorphic to the center $Z$ of $G$.
 \end{itemize}
 \end{intthm}
\SSP

The above theorem asserts the existence  of   a canonical  compact moduli space of dormant $G$-opers  admitting a virtual fundamental class.
Recall that virtual fundamental classes of moduli spaces play a central role in enumerative geometry as they represent a major ingredient in the construction of deformation  invariants, e.g., Gromov-Witten and   Donaldson-Thomas  invariants.
Virtual fundamental classes of Quot schemes  were  constructed in  ~\cite[Theorem 1]{MO}.
Although there exists a connection  between certain  Quot schemes  and the moduli space of $\mr{PGL}_n$-do'pers on a sufficiently general  curve (cf. ~\cite[\S\,9]{Wak5} or an unpublished version\footnote{This version is available at:  {\tt https://arxiv.org/pdf/1311.4359v1.pdf}} of ~\cite[\S\,11]{Jo14}),
 the  virtual class $[\mfO \mfp_{G, g,r}^{^\mr{Zzz...}}]^\mr{vir}$ does not seem to come, {\it a priori}, from Quot schemes because of its construction.

\LSP
\subsection*{0.5} \label{S105}
Next, by means of the resulting virtual fundamental class, we construct   a {\it CohFT} (= a {\it cohomological field theory}) associated with  $G$-do'pers, which   in fact  forms  a {\it $2$d TQFT} (= a {\it $2$-dimensional topological quantum field theory}).

Cohomological field theories 
 were introduced by M. Kontsevich and Y. Manin in ~\cite{KM} to axiomatize the properties of  Gromov-Witten classes of a given target variety over the field of complex numbers $\mbC$.
For instance, the trivial CohFT arises from the Gromov-Witten theory of one point.
As it turns out,  this notion is more general, in the sense that not all CohFTs come from Gromov-Witten theory. 
Indeed, one may find some examples of CohFTs, including  Hodge CohFTs and CohFTs arising from Witten's $r$-spin classes or Fan-Jarvis-Ruan-Witten (FJRW) theory.

On the other hand, in the axiomatic formulation due to Atiyah (cf. ~\cite{At}),
an $N$-dimensional topological quantum field theory (or just an $N$d TQFT) is a rule $\Lambda$ which to each closed oriented manifold $\Sigma$ of dimension $N-1$ associates a vector space $\Lambda_\Sigma$,
and to each oriented manifold of dimension $N$ whose boundary is $\Sigma$ associates a vector in $\Lambda_\Sigma$.
This rule satisfies  a collection of axioms which express that topologically equivalent manifolds have isomorphic associated vector spaces, and that disjoint unions of manifolds go to tensor products of vector spaces, etc.
 One may consider each  $2$d TQFT as a special kind of  CohFT, i.e., a CohFT valued  in the $0$-th cohomology.

Hereinafter, suppose that $k$ is algebraically closed.
In \S\,\ref{SSS666}, we introduce (cf. Definition \ref{D07}) the definition of a CohFT by means of  the $l$-adic \'{e}tale cohomology of 
$\overline{\mfM}_{g,r}$
and recall the notion of a $2$d TQFT.
One of   key ingredients  in the  construction of  the desired CohFT (or $2$d TQFT)  is 
the  stack $\mfR \mfa \mfd$   (cf. (\ref{EE900}))  defined as follows.
Let $\mft$ be the Lie algebra  of a fixed maximal torus of $G$ (viewed as a $k$-scheme).
Denote by $\mft^F_{\mr{reg}}$ (cf. (\ref{ee500})) the subscheme of  $\mft$  consisting of the Frobenius-invariant regular elements.
The natural action of the Weyl group $W$ on $\mft^F_{\mr{reg}}$ yields 
 the quotient scheme $\mft_\mr{reg}^F/W$.
If $Z$ denotes the center of $G$, then
the trivial $Z$-action  gives rise to  the quotient stack $[(\mft_\mr{reg}^F/W)/Z]$. 
This stack induces  the  stack of cyclotomic gerbes 
 $\overline{\mcI}_\mu ([(\mft_\mr{reg}^F/W)/Z])$, which we 
     denote by $\mfR \mfa \mfd$.
The connected components of  $\overline{\mcI}_\mu ([(\mft_\mr{reg}^F/W)/Z])$ are  indexed by a certain finite set $\Delta$ (cf. (\ref{ee142})). 
Each connected component $
\mfR \mfa \mfd_{\pmb{\rho}}$ (indexed by  $\pmb{\rho} \in \Delta$)
is canonically isomorphic to  the classifying stack $\mcB \overline{Z}$ of a certain  quotient group $\overline{Z}$ of $Z$ (cf. Proposition \ref{Lw001}).

Another key ingredient is  called  {\it  radius} (cf. Definition \ref{y035}); this  is an invariant  associated with each $G$-do'per (or more generally, each faithful twisted $G$-oper) and each marked point of the underlying curve.
 The notion of   radius  introduced  in the present paper  generalizes the classical definition  given  in ~\cite[Definition 1.2]{Mzk2}  and ~\cite[Definition 2.29]{Wak5}.
The reason for the name of this notion is  based on 
 a certain  analogy between the theory of ``{\it pants}" in  the classical Teichm\"{u}ller theory  and the structure theory of $\mr{PGL}_2$-opers on a Riemann sphere minus three points.
In fact, we may make  this analogy  by identifying  the  radii in our sense with  the size of the respective  holes in the Riemann sphere.

The assignment from each $G$-do'per  to its  radius  at the $i$-th marked point ($i =1, \cdots, r$) may be realized geometrically as a morphism
$\mr{ev}_i  : \mfO \mfp^{^\mr{Zzz...}}_{G,g,r}   \migi \mfR \mfa \mfd  \left(= \coprod_{\pmb{\rho} \in \Delta} \mfR \mfa \mfd_{\pmb{\rho}} \right)$.
That is to say,   the point classifying a $G$-do'per  is mapped, via  $\mr{ev}_i$,   to  a point in the connected component $\mfR \mfa \mfd_{\pmb{\rho}}$ indexed by  its radius $\pmb{\rho}$ at the $i$-th marked point. 
Thus, we obtain    a diagram of  stacks 
\begin{align}
\begin{CD}
\mfO \mfp^{^\mr{Zzz...}}_{G, g,r}  @> \pi_{g,r} >> \overline{\mfM}_{g,r}
\\
@V \mr{ev}_i VV @.
\\
\mfR \mfa \mfd, @.
\end{CD}
\end{align}
where $\pi_{g,r}$ denotes the natural projection.

By means of
the various objects
   just discussed, the constituents in  the desired CohFT are  constructed.
Given a prime $l$ different from $p$, let us consider the $l$-adic \'{e}tale cohomology  
\begin{align}
\mcV := \widetilde{H}^*_{\text{\'{e}t}}(\mfR \mfa \mfd, \overline{\mbQ}_l) 
\end{align}
of the stack $\mfR \mfa \mfd$,  where $\overline{\mbQ}_l$ denotes
   the algebraic closure of $\mbQ_l$ 
  and we write $\widetilde{H}^*_{\text{\'{e}t}}(-, \overline{\mbQ}_l) := \bigoplus_{i=0}^\infty H^i_{\text{\'{e}t}}(-, \overline{\mbQ}_l (\lfloor{\textstyle \frac{i}{2}}\rfloor))$.
We obtain  a collection   $\{ \Lambda_{G, g,r} \}_{g,r \geq 0, 2g-2+r >0}$ consisting of the $\overline{\mbQ}_l$-linear morphisms
\begin{align} \label{ee201}
\Lambda_{G, g,r} : & \hspace{2mm} \mcV^{\otimes r} \hspace{2mm}  \migi \hspace{2mm}  \widetilde{H}^*_{\text{\'{e}t}} (\overline{\mfM}_{g,r}, \overline{\mbQ}_l) \\
&  \hspace{3mm}  \vin    \hspace{26mm} \vin  \notag  \\[-4pt]
& \bigotimes_{i=1}^r v_{i}  \mapsto   \left(\pi_{g,r*}^\mr{hom} \left(\left(\prod_{i=1}^r \mr{ev}^*_i (v_i) \right) \cap \mr{cl}^{3g-3+r} ([\mfO \mfp^{^\mr{Zzz...}}_{G,g,r}]^\mr{vir})\right) \right)^\blacklozenge \notag
\end{align}
(cf.
 (\ref{ee502}) for its precise definition).
 Also, a $\overline{\mbQ}_l$-bilinear pairing $\eta : \mcV \times \mcV \migi \overline{\mbQ}_l$ (cf. (\ref{ee80255})) and a distinguished element $e_{\pmb{\varepsilon}}$ of $\mcV$ (cf. (\ref{fff444})) are defined. 
Then, the second main result is as follows.

\SSP
\begin{intthm} [cf.  Theorem \ref{T10}, (ii)] \label{TBB}
 Assume that $G$ satisfies  the condition $(**)_{G}$.
 Then, the collection of data 
 \begin{align}
\Lambda_{G}  := (\mcV, \eta, e_{\pmb{\varepsilon}},  \{ \Lambda_{G, g,r} \}_{g,r \geq 0, 2g-2+r >0})
 \end{align}
(cf.  (\ref{ee148})) 
  forms a  CohFT (with flat identity) valued in $\widetilde{H}_{\text{\'{e}t}}^0 (\overline{\mfM}_{g,r}, \overline{\mbQ}_l)$, namely,  forms a $2$d TQFT over $\overline{\mbQ}_l$. 
Moreover, 
the corresponding Frobenius algebra $(\mcV, \eta)$ is semisimple.
  \end{intthm}
\SSP

The $2$d TQFT asserted above  will play an important role in studying enumerative geometry of do'pers.
In fact, each  $\vec{\pmb{\rho}} := (\pmb{\rho}_i)_{i=1}^r \in \Delta^{\times r}$ determines
  the closed substack $\mfO \mfp^{^\mr{Zzz...}}_{G,g,r, \vec{\pmb{\rho}}}$ (cf. (\ref{ee800}))    of $\mfO \mfp^{^\mr{Zzz...}}_{G,g,r}$ classifying $G$-do'pers of radii $\vec{\pmb{\rho}}$.
According to Theorem \ref{T10}, (i),   the value $\Lambda_{G, g,r}  (\bigotimes_{i=1}^r e_{\pmb{\rho}_i})$ 
 coincides with
the generic degree of $\mfO \mfp^{^\mr{Zzz...}}_{G,g,r, \vec{\pmb{\rho}}}/\overline{\mfM}_{g,r}$.
In particular, 
 it 
 is exactly  equal to  $\frac{1}{ |Z|}$ times 
the number of  (isomorphism classes of) $G$-do'pers of radii $\vec{\pmb{\rho}}$ on a sufficiently general curve classified by $\overline{\mfM}_{g,r}$.
 An explicit description of 
 the ring structure (e.g., the set of characters) of the Frobenius algebra $(\mcV, \eta)$ allows us to  compute these values.

In ~\cite[Proposition 7.33]{Wak5}, it was already proved that the collection $\{ \Lambda_{G, g,r}\}_{g, r}$ together form  a {\it pseudo-fusion rule}, i.e., 
a somewhat weaker variant of  a fusion rule, as well as a $2$d TQFT.
The notion of a pseudo-fusion rule provides  just enough framework
 to compute the generic degrees, but  at the time of writing ~\cite{Wak5} some of the conditions for it to be 2d TQFT remained unknown.
After proving  Theorem \ref{TBB}, one can obtain 
the same 
 calculations as obtained in {\it loc.\,cit.} by simply applying the general theory of 2d TQFTs.

For example,  
  the structure  of the Frobenius algebra $(\mcV, \eta)$ for $G = \mr{PGL}_2$  can be completely understood by comparing with the fusion ring of the $\mfs \mfl_2 (\mbC)$ Wess-Zumino-Witten model  (cf. \S\,\ref{SSjj5d}).
 In fact, under the natural  identification $\Delta = \{ 0, 1, \cdots, \frac{p-3}{2} \}$ (cf. (\ref{ER52})),  
 we obtain the following equality   defined  for each collection of radii $(n_i)_{i=1}^r \in \Delta^{\times r}$:
 \begin{align} \label{ee80ee2}
 \Lambda_{\mr{PGL}_2, g,r} (\bigotimes_{i=1}^r e_{n_i}) & \left(= \mr{deg} (\mfO \mfp^{^\mr{Zzz...}}_{\mr{PGL}_2, g,r, (n_i)_{i=1}^r}/\overline{\mfM}_{g,r}) \right) 
  = \frac{p^{g-1}}{2^{2g-1}}  \sum_{j =1}^{p-1} \frac{ \prod_{i=1}^r \sin \left( \frac{(2n_i +1) j \pi}{p}\right)}{\sin^{2g-2+r} \left(\frac{j \pi}{p} \right)}. 
 \end{align}
 (cf. ~\cite[Theorem 7.41]{Wak5}).

\LSP
\subsection*{0.6} \label{S1016}
In the final  section of the present paper, we discuss an analogue for  $
\Lambda_{G}$ of the Witten-Kontsevich theorem (cf. ~\cite{Wit1},  ~\cite{K1}, ~\cite{KaLa}, ~\cite{KiLi}, ~\cite{Mirz}, and ~\cite{Okou}), which is 
one of the landmark results concerning the intersection theory on the moduli 
stacks  of pointed stable curves   over $\mbC$.
The well-known Witten-Kontsevich theorem  asserts an equivalence of the intersection theory of  psi  classes on 
that space and 
 the Hermitian matrix model of $2$-dimensional gravity.
This implies that the partition function of the trivial CohFT is a tau function 
of the  KdV hierarchy, in other words, it satisfies a certain series  of partial differential equations.

The  partition function that we  deal with 
  is defined as follows (cf. Definition \ref{DFF061}):
\begin{align} \label{ee200}
 &  Z_{G} 
  :=  \mr{exp} \left(
 \sum_{g,r \geq 0}
  \frac{\hbar^{2g-2}}{r!} \hspace{-2mm}\sum_{\genfrac{.}{.}{0pt}{}{d_1, \cdots, d_r \geq 0,}{\pmb{\rho}_1, \cdots, \pmb{\rho}_r \in \Delta}} 
\hspace{-2mm}\left(\int_{[\mfO \mfp_{G, g,r}^{^\mr{Zzz...}}]^\mr{vir}}  \prod_{i=1}^r \mr{ev}_i^{*} (e_{\pmb{\rho}_i})\widehat{\psi}_i^{d_i} \right)
 \prod_{i=1}^r t_{d_i, \,  \pmb{\rho}_i} \right) 
  \\
 & \hspace{5mm} \left(= \mr{exp} \left(
 \sum_{g,r \geq 0}
   \frac{\hbar^{2g-2}}{r!} \hspace{-2mm}\sum_{\genfrac{.}{.}{0pt}{}{d_1, \cdots, d_r \geq 0,}{\pmb{\rho}_1, \cdots, \pmb{\rho}_r \in \Delta}} 
\hspace{-2mm}\left(\int_{\overline{\mfM}_{g,r}}\hspace{-2mm} \Lambda_{G, g,r} (\bigotimes_{i=1}^r e_{\pmb{\rho}_i}) \prod_{i=1}^r \psi_i^{d_i} \right)
 \prod_{i=1}^r t_{d_i, \,  \pmb{\rho}_i} \right)\right) \notag \\
&  \hspace{9mm}  \in  \overline{\mbQ}_l ((\hbar)) [[ \{t_{d,\, \pmb{\rho}} \}_{d \geq 0, \,  \pmb{\rho} \in \Delta}]], \notag
\end{align}
 where $\widehat{\psi}_i$ and $\psi_i$ ($i =1, \cdots, r$) denote  the $i$-th psi classes on $\mfO \mfp^{^\mr{Zzz...}}_{G, g,r}$ and $\overline{\mfM}_{g,r} $ respectively,  and the equality in the parenthesis follows from Proposition \ref{T02}.
Also,  let 
 $L_n$ ($n \geq  -1$) denote 
 a differential operator 
\begin{align} \label{ee505}
L_n := &  - \frac{(2n+3)!!}{2^{n+1}} \frac{\partial}{\partial t_{n+1,  \, \pmb{\varepsilon}}} + \sum_{i=0}
^\infty \frac{(2i+2n+1)!!}{(2i-1)!!2^{n+1}} \left( \sum_{\pmb{\rho} \in \Delta}^{} t_{i, \, \pmb{\rho}} \frac{\partial}{\partial t_{i+n, \, \pmb{\rho}}}\right) \\
& + \frac{ |Z| \hbar^2}{2} \sum_{i=0}^{n-1}\frac{(2i+1)!! (2n-2i-1)!!}{2^{n+1}} \left(\sum_{\pmb{\rho} \in \Delta}\frac{\partial^2}{\partial t_{i, \, \pmb{\rho}} \partial t_{n-1-i, \, \pmb{\rho}^\veebar}} \right) \notag \\
& + \delta_{n, -1} \frac{\hbar^{-2}}{2 |Z|} \left( \sum_{\pmb{\rho} \in \Delta}  t_{0, \, \pmb{\rho}}t_{0, \, \pmb{\rho}^\veebar}\right) + \delta_{n, 0} \frac{|\Delta|}{16}. \notag 
\end{align}
Then, 
we  prove the equality  
\begin{align} \label{ee865}
L_n Z_{G} =0
\end{align}
  for every  $n \geq -1$ (cf. Theorem \ref{TT01}).
That is to say, the  partition function  of $G$-do'pers  turns out to be  a solution to infinitely many partial differential equations $L_n y=0$ ($n \geq -1$).
In particular, this result gives  nontrivial constraints on the psi classes $\widehat{\psi}_i$ on $\mfO \mfp^{^\mr{Zzz...}}_{G, g,r}$.

\LSP
 \subsection*{Acknowledgements} 
 The author cannot  express enough his sincere and deep gratitude to Professors Shinichi Mochizuki and Kirti Joshi.
Without their philosophical viewpoints, theoretical insights, and endless creativity, my study of mathematics   would have remained ``{\it dormant}".
Also,  special thanks go to the moduli stack  of dormant $G$-opers $\mfO \mfp^{^\mr{Zzz...}}_{G,g,r}$, who has guided him to  the beautiful world of mathematics.
The author was partially  supported by the Grant-in-Aid for Scientific Research (KAKENHI No.\,18K13385, 21K13770).

\LSP
\subsection*{Notation and Conventions}
Let us introduce some  notation and conventions used in  the present paper.
Throughout the present paper, all schemes and  algebraic stacks are assumed to be locally noetherian.
We  fix a perfect field $k$, and 
denote by $\mfS \mfc \mfh_{/k}$ the category of schemes over $k$.

We use the word {\it stack} to mean algebraic stack in the sense of the appendix in ~\cite{Vis1}.
Let $\mfX$ be a Deligne-Mumford stack.
By a sheaf on $\mfX$, we mean, unless otherwise stated,  a sheaf  on the {\it small \'{e}tale site} of $\mfX$.
In particular, one obtains the {\it structure sheaf} $\mcO_\mfX$ of $\mfX$. 
Under the assumption that $\mfX$ is of finite type over $k$,
we denote by $| \mfX |$ the coarse moduli space of $\mfX$, which has  a natural  projection $\mr{coa}_\mfX : \mfX \migi | \mfX |$.

Each morphism $f: \mfX \migi \mfY$ of Deligne-Mumford stacks of finite type  over $k$
induces naturally a morphism 
$| f | :  | \mfX | \migi | \mfY |$
between their respective   coarse moduli spaces.
If, moreover, both $\mfX$ and $\mfY$ are integral and  $f$ is  separated and dominant, then
the {\it degree} $\mr{deg}(\mfX/ \mfY)$ of $\mfX$ over $\mfY$  are defined (cf. ~\cite[(1.15), Definition]{Vis1}).
Notice that the value $\mr{deg}(\mfX/ \mfY)$ may not be an integer in general unless $f$ is representable.
One may generalize naturally the notion of degree to the case where $\mfX$ is a disjoint union of finite number of integral stacks.
For each $n \in \mbZ_{\geq 0}$, denote by $A_n (\mfX)_\mbQ$ (cf. ~\cite[(3.4), Definition]{Vis1}) the rational Chow group of cycles of dimension
$n$ on $\mfX$ modulo rational equivalence tensored with $\mbQ$.

Basic references for the notion and properties  of a {\it log scheme} (or more generally, a {\it log stack}) are ~\cite{ILL}, ~\cite{KaFu}, and ~\cite{KaKa}.
For a log  stack indicated, say,  by $\mfY^\mr{log}$, we shall write $\mfY$ for the underlying stack of $\mfY^\mr{log}$.
For a morphism $\mfY^\mr{log} \migi \mfT^\mr{log}$ of fs log Deligne-Mumford stacks,
let us write $\mcT_{\mfY^\mr{log}/\mfT^\mr{log}}$  for the sheaf of logarithmic derivations of $\mfY^\mr{log}$ over $\mfT^\mr{log}$.
Also,   write $\Omega_{\mfY^\mr{log}/\mfT^\mr{log}}$ for its dual $\mcT_{\mfY^\mr{log}/\mfT^\mr{log}}^\vee$, i.e.,  the sheaf of logarithmic
  $1$-forms
   of $\mfY^\mr{log}$ over $\mfT^\mr{log}$.

For each positive integer $l$, we denote by $\mu_l$ the group of $l$-th roots of unity in an  algebraic closure of $k$.
If $G$ is  an  algebraic group  over $k$, then let  $\mcB G$ denote the classifying stack of $G$,  
which is defined as the quotient stack $\mcB G := [\mr{pt} / G]$
   for   the trivial action of $G$ on $\mr{pt} := \mr{Spec} (k)$. 

  Given a right $G$-bundle $\mcE$ on a $k$-stack $\mfY$ in the \'{e}tale topology  and  a $k$-vector space $\mfh$  equipped with a left $G$-action, we shall write $\mfh_\mcE$ for the vector bundle  on $\mfY$ associated with the relative affine space $\mcE \times^G \mfh$ ($: = (\mcE \times_k \mfh)/G$).
Denote by  $\mbG_m$ the multiplicative group over $k$.

Let $S$ be a scheme. By a {\bf nodal curve} over $S$, we mean a flat morphism of schemes   $f : X \migi S$  whose geometric fibers are connected and  reduced $1$-dimensional schemes with at most nodal  points as singularities. 
For simplicity, we often  write $X$ instead of $f : X \migi S$ when  denoting  a nodal curve over $S$.

Let $r$ be  a nonnegative integer.
An {\bf $r$-pointed nodal curve} over $S$  is defined to be a collection $\msX := (f : X \migi S, \{\sigma_{X, i} : S \migi X\}_{i=1}^r)$ consisting of a nodal curve $f : X \migi S$ over $S$  and $r$ sections $\sigma_{X, i} : S \migi X$ ($i =1, \cdots, r$) such that $\mr{Im} (\sigma_{X, i})$  lies, for any $i$,  in the smooth locus of $X$ (relative to $S$) and  that 
$\mr{Im}(\sigma_{X, i}) \cap \mr{Im}(\sigma_{X, j}) = \emptyset$
for any pair $(i, j)$ with $i \neq j$.

\vspace{10mm}
\section{Extended spin structures on  twisted  curves} \label{S01} 
\SSP

This section deals with   {\it extended spin structures}, 
 which are  some sort of natural generalization of  spin structures (in the classical sense)  defined on a twisted curve.
   The notion of an extended spin structure will be used to describe the gap between  a faithful twisted $G$-oper 
      and  the associated  $G_\mr{ad}$-oper, where $G_\mr{ad}$ denotes the adjoint group of $G$.
   In fact, we will see  (cf. Theorem \ref{P0012}) that the moduli stack of faithful twisted $G$-opers 
   is isomorphic to the product of the moduli stacks of $G_\mr{ad}$-opers and  of  certain extended spin structures.
   In this section, we prove  some properties of the moduli stack of extended spin structures (cf., e.g., Theorem \ref{P08f}, Propositions \ref{P08}, \ref{Pw081}, and \ref{P020}).

\LSP
\subsection{The moduli stack of pointed stable curves} \label{SS01}
To begin with, we shall introduce some notation concerning pointed curves and the  moduli stack classifying them.
Given a pair of nonnegative integers $(g, r)$ with $2g-2+r >0$,
we shall write  
$\overline{\mfM}_{g,r, \mbZ}$
for  the moduli stack of $r$-pointed stable curves of genus $g$, and write 
\begin{align}
\overline{\mfM}_{g,r} := \overline{\mfM}_{g,r, \mbZ} \times_{\mbZ} k.
\end{align}
Namely,  $\overline{\mfM}_{g,r}$ classifies the pointed stable curves
\begin{align}
\msX := (f: X \migi S, \{ \sigma_{X, i} : S \migi X \}_{i=1}^r)
\end{align}
 consisting  of a proper 
 nodal
 curve $f : X \migi S$ of genus $g$ over a $k$-scheme $S$   and $r$ marked points $\sigma_{X, i} : S \migi X$ ($i =1, \cdots, r$) satisfying certain conditions (cf. ~\cite[Definition 1.1]{Kn2}).
Recall (cf. ~\cite[Corollary 2.6 and Theorem 2.7]{Kn2}, ~\cite[\S\,5]{DM}) that $\overline{\mfM}_{g,r}$ may be represented by a geometrically connected, proper, and smooth Deligne-Mumford stack over $k$ of dimension $3g-3+r$.
Denote by $\mfM_{g,r}$ the dense open substack of $\overline{\mfM}_{g,r}$ classifying smooth curves.

Given  nonnegative integers $g_1$, $g_2$, $r_1$, $r_2$ with $2g_i-1+r_i >0$ ($i =1, 2$),
we shall write
\begin{align} \label{e072}
\Phi_\mr{tree} : \overline{\mfM}_{g_1, r_1 +1} \times_k \overline{\mfM}_{g_2, r_2 +1} \migi \overline{\mfM}_{g_1 +g_2, r_1 +r_2}
\end{align}
for the gluing map obtained by attaching the respective last marked points 
of curves classified by $\overline{\mfM}_{g_1, r_1 +1}$ and $\overline{\mfM}_{g_2, r_2 +1}$
to form a node.

Similarly, given  nonnegative integers $g$, $r$ with $2g+r >0$,
we shall write
\begin{align} \label{e071}
\Phi_\mr{loop} : \overline{\mfM}_{g, r+2} \migi \overline{\mfM}_{g+1, r}
\end{align}
for the gluing map 
 obtained  by attaching the last two  marked points of each curve classified by $\overline{\mfM}_{g, r+2}$ to form a node.

Finally, if $2g -2+r >0$, we define 
\begin{align} \label{e070}
\Phi_{\mr{tail}} : \overline{\mfM}_{g, r+1} \migi \overline{\mfM}_{g,r}
\end{align}
to be the morphism obtained by  forgetting the last marked point and successively contracting any resulting unstable components of each curve classified by $\overline{\mfM}_{g, r+1}$.

\LSP
\subsection{Twisted curves} \label{SSg01}
Next, 
let us recall  the notion of a pointed twisted curve and the log structure equipped with it.
This notion  will be used to define  a twisted $G$-oper.
Here, recall 
 that a Deligne-Mumford stack $\mfX$ over $k$  is called  {\bf tame} if, for any algebraically closed field $\overline{k}$ over $k$
 and any morphism $\overline{x} : \mr{Spec}(\overline{k}) \migi \mfX$,  the stabilizer group $\mr{Stab}_\mfX (\overline{x})$ of $\overline{x}$ has order invertible in $k$.

\SSP
\bde \label{D01}  Let $T^\mr{log}$ be an fs log scheme over $k$.
 A {\bf stacky  log curve} over $T^\mr{log}$ is an fs log Deligne-Mumford stack $\mfY^\mr{log}$ together with a log smooth integrable morphism $f^\mr{log} : \mfY^\mr{log} \migi T^\mr{log}$ such that the geometric fibers of the underlying morphism $f : \mfY \migi T$ of stacks  are reduced, connected, and $1$-dimensional.
 (In particular, both $\mcT_{\mfY^\mr{log}/T^\mr{log}}$ and $\Omega_{\mfY^\mr{log}/T^\mr{log}}$ are line bundles on $\mfY$.)
 \ede

\SSP
\bde[cf. ~\cite{O1}, Definition 1.2] \label{D011}
  Let $S$ be a $k$-scheme.
\begin{itemize}
\item[(i)]
A {\bf (balanced) twisted  curve} over $S$ is 
a proper   flat morphism $f: \mfX \migi S$ of  tame Deligne-Mumford stacks over $k$
 satisfying the following conditions:
\begin{itemize}
\item
The geometric  fibers of $f : \mfX \migi S$ are purely $1$-dimensional and 
 connected with at most nodal singularities.
\item
If $| \mfX |^\mr{sm}$ denotes the open subscheme  of $| \mfX |$  where $|f| : | \mfX | \migi S$ is smooth, then
the inverse image $\mfX \times_{|\mfX |} | \mfX |^\mr{sm}$ ($\subseteq \mfX$) coincides with the open substack of $\mfX$ where $f$ is smooth.
\item[$\bullet$]
For any geometric point $\overline{s} \migi S$, the map $\mr{coa}_\mfX \times \mr{id}_{\overline{s}} : \mfX \times_{S} \overline{s} \migi | \mfX | \times_{S} \overline{s}$ is an isomorphism over some dense open subscheme  of $| \mfX | \times_{S} \overline{s}$.
\item[$\bullet$]
Consider a geometric point $\overline{x} \migi | \mfX |$ mapping to a node.
Note that 
there exist
an affine open neighborhood  $T \ (:= \mr{Spec} (R))  \subseteq S$ of $| f | (\overline{x})$, 
an affine \'{e}tale neighborhood $\mr{Spec} (A) \migi |f|^{-1} (T)$ ($\subseteq | \mfX |$) of $\overline{x}$,  and an  \'{e}tale morphism
\begin{align} \label{WW130}
\mr{Spec} (A) \migi \mr{Spec} (R [s_0, t_0]/(s_0t_0 -u_0))
\end{align}
over $R$ for some $u_0 \in R$.
Then,   the pull-back $\mfX \times_{| \mfX |} \mr{Spec} (A)$ is isomorphic to the quotient stack
\begin{align} \label{Ww131}
[\mr{Spec} (A [ s_1, t_1]/(s_1 t_1 -u_1, s_1^l - s_0, t_1^l - t_0))/\mu_l]
\end{align}
for some   $u_1 \in R$ and some positive integer $l$ invertible in $k$
such that 
$\zeta \in \mu_l$ acts by $(s_1, t_1) \mapsto (\zeta s_1, \zeta^{-1} t_1)$.
\end{itemize}
\item[(ii)]
Let $g$ be a nonnegative integer.
We shall say that a twisted curve $f : \mfX \migi S$ is  {\bf of genus $g$}
if the genus of every   fiber of the proper nodal curve $|f| : | \mfX | \migi S$
 coincides with  $g$.
\end{itemize}
  \ede
\SSP

\bde[cf. ~\cite{AV}, Definition 4.1.2] \label{D012} 
Let $g$ and $r$ be  nonnegative integers and $S$ a $k$-scheme.
\begin{itemize}
\item[(i)]
An {\bf $r$-pointed twisted curve (of genus $g$)}  over $S$ is a collection of data
\begin{align}
\msX^\mr{tw} := (f : \mfX \migi S, \{ \sigma_{\mfX, i} : \mfS_i \migi \mfX \}_{i=1}^r)
\end{align}
consisting of a twisted curve $f : \mfX \migi S$ (of genus $g$)  and disjoint closed substacks
$\sigma_{\mfX, i} : \mfS_i \migi \mfX$ ($i =1, \cdots, r$) of $\mfX$ 
   satisfying the following conditions:
\begin{itemize}
\item[$\bullet$]
Each $\mr{Im} (\sigma_{\mfX, i})$ is contained in the smooth locus 
in $\mfX$ (relative to $S$).
\item[$\bullet$]
Each  $\mfS_i$ is  \'{e}tale gerbe over $S$.
\item[$\bullet$]
If $\mfX^{\mr{gen}}$ denotes the complement of the union of $\mr{Im}(\sigma_{\mfX, i})$ ($i =1, \cdots, r$) in the smooth locus
in $\mfX$ (relative to $S$),
then $\mfX^{\mr{gen}}$ may be represented by  a scheme.
\end{itemize}
\item[(ii)]
Let $\msX^\mr{tw}_j := (f_j : \mfX_j \migi S_j, \{ \sigma_{\mfX_j, i} : \mfS_{j, i} \migi \mfX_j \}_{i =1}^r)$ ($j=1,2$) be $r$-pointed  twisted curves.
A {\bf morphism of $r$-pointed twisted curves} from $\msX^\mr{tw}_1$ to $\msX^\mr{tw}_2$ is a pair of morphisms 
\begin{align}
(\alpha_S : S_1 \migi S_2, \alpha_\mfX : \mfX_1 \migi \mfX_2)
\end{align}
 such  that $\alpha^{-1}_\mfX (\mfS_{2, i}) = \mfS_{1, i}$ for any $i =1, \cdots, r$,  and moreover,   the  square diagram
\begin{align}
\begin{CD}
\mfX_1 @> f_1 >> S_1
\\
@V \alpha_\mfX VV @VV \alpha_S V
\\
\mfX_2 @>>f_2 > S_2
\end{CD}
\end{align}
is commutative and cartesian.
\item[(iii)]
Suppose that $2g-2+r >0$.
A {\bf twisted stable curve of type $(g,r)$} over $S$  is an $r$-pointed twisted curve $\msX^\mr{tw} := (f : \mfX \migi S, \{ \sigma_{\mfX, i} : \mfS_i \migi \mfX \}_{i=1}^r)$ over $S$  whose    coarse moduli space 
\begin{align} \label{ee230}
| \msX^\mr{tw} | := (| f | : | \mfX | \migi S, \{ |\sigma_{\mfX, i} | : | \mfS_i | \migi | \mfX | \}_{i=1}^r)
\end{align}
   forms an $r$-pointed stable curve of genus $g$.
 \end{itemize}
  \ede
\SSP

\bde
 \label{D012} 
 Let $r$ be 
a nonnegative integer and $S$ a $k$-scheme.
\begin{itemize}
\item[(i)]
Let 
 $\msX := (f: X \migi S, \{ \sigma_{X, i} \}_{i=1}^r)$ be an $r$-pointed nodal curve over $S$.
  A {\bf twistification} of $\msX$ is a pair
  \begin{align}
 (\msX^\mr{tw}, \gamma)
 \end{align}
consisting of an $r$-pointed  twisted  curve $\msX^\mr{tw} := (f : \mfX \migi S, \{ \sigma_{\mfX, i} \}_{i=1}^r)$ and 
an $S$-morphism $\gamma : \mfX \migi X$
such that the induced morphism
  $|\gamma | : | \mfX| \migi X$ is an isomorphism and 
the equality $|\gamma \circ \sigma_{\mfX, i} | = \sigma_{X, i}$ holds  for every  $i = 1, \cdots, r$.
\item[(ii)]
Let  $\msX$ 
 be an $r$-pointed nodal  curve  over $S$ and
let $(\msX^\mr{tw}_j, \gamma_j)$ ($j=1,2$) be  twistifications  of $\msX$.
An {\bf isomorphism of twistifications} from $(\msX^\mr{tw}_1, \gamma_1)$  to $(\msX^\mr{tw}_2, \gamma_2)$ is an isomorphism of $r$-pointed twisted curves  from $\msX^\mr{tw}_1$ to $\msX^\mr{tw}_2$ 
 compatible with $\gamma_1$ and $\gamma_2$.
\end{itemize}
  \ede
\SSP

Let $r$ be a nonnegative integer and $\msX^\mr{tw} := (f : \mfX \migi S, \{ \sigma_{\mfX, i} \}_{i=1}^r)$  an $r$-pointed twisted curve.
There exist   canonical  log structures on $\mfX$ and $S$ obtained by  the log structures described in 
~\cite[Theorem 3.6]{O1}
and the closed substacks $\sigma_{\mfX, i}$; we denote the resulting 
 log stacks  by $\mfX^\mr{log}$ and $S^\mr{log}$ respectively.
The  morphism $f : \mfX \migi S$  extends to a log smooth morphism  between  log stacks  $f^\mr{log} : \mfX^{\mr{log}} \migi S^{\mr{log}}$, by which $\mfX^{\mr{log}}$ specifies a stacky  log curve over $S^{\mr{log}}$.
If $\mfX$ is smooth, then $S^{\mr{log}} = S$.

\LSP
\subsection{Extended spin structures.} \label{y034}
In the rest of this section, we assume that $\mr{char} (k) \neq 2$.
Also, assume that we are given  a pair $(Z, \delta)$ consisting of a finite abelian group $Z$ (often regarded as  a finite group scheme over $k$) which has  order  invertible in $k$ and a morphism  of groups $\delta : (\{ \pm 1 \} =) \ \mu_2 \migi Z$.
Write 
\begin{align} \label{ee880}
\overline{\delta} : \mu_{\overline{2}} \ (:= \mu_2 / \mr{Ker} (\delta)) \migiincl  Z,
\end{align}
 where $\overline{2}$ is either $1$ or $2$,  for the injection  induced naturally by $\delta$. 
Denote by 
\begin{align} \label{ee280}
\widehat{Z}_\delta
\end{align}
 the cofiber product of $\delta$ and the natural inclusion $\mu_2 \migiincl \mbG_m$, which fits into the following morphism of short exact sequences of algebraic groups:
\begin{align} \label{ee028}
\begin{CD}
1 @>>> \mu_2 @> \mr{incl.} >> \mbG_m @>  x \mapsto x^2  >> \mbG_m @>>> 1 
\\
@. @VV \delta V @VV \delta_{\mbG_m}V @VV \mr{id}_{\mbG_m} V @.
\\
 1 @>>> Z @>>> \widehat{Z}_\delta @> \nu >> \mbG_m @>>> 1.
 \end{CD}
\end{align}

\SSP
\bde \label{y035dd}\begin{itemize}
\item[(i)]
Let $S^\mr{log}$ be an fs log scheme over $k$ and $\mfU^\mr{log}$ a stacky log curve over $S^\mr{log}$.
An  {\bf extended $(Z, \delta)$-spin structure} (or, a {\bf $(Z, \delta)$-structure}, for short) on $\mfU^\mr{log}/S^\mr{log}$ is 
a $\widehat{Z}_{\delta}$-bundle
\begin{align}
\pi_\mcG : \mcG \migi \mfU
\end{align}
on 
$\mfU$ in the \'{e}tale topology
 whose  classifying  morphism $\mfU \migi \mcB \widehat{Z}_\delta$ is representable and  fits into  the following  $1$-commutative diagram:
\begin{align} \label{ee030}
\xymatrix{
\mfU \ar[rr]^{[\Omega_{\mfU^\mr{log}/S^\mr{log}}]}
 \ar[rd]& & \mcB \mbG_m   \\
&\mcB \widehat{Z}_\delta,  \ar[ru]_{\mcB \nu}&
}
\end{align}
where the upper horizontal morphism $[\Omega_{\mfU^\mr{log}/S^\mr{log}}]$ denotes the classifying morphism  of the line bundle $\Omega_{\mfU^{\mr{log}}/S^{\mr{log}}}$ and $\mcB \nu$ denotes the morphism induced naturally by $\nu$ (cf.  (\ref{ee028})). 
\item[(ii)]
Let $r$ be a nonnegative integer, $S$ a $k$-scheme, and  $\msX := (X, \{ \sigma_{X, i} \}_{i=1}^r)$   an $r$-pointed nodal curve over $S$.
A {\bf $(Z, \delta)$-structure} on $\msX$ is a collection of data
\begin{align}
(\msX^\mr{tw}, \gamma, \pi_\mcG : \mcG \migi \mfX),
\end{align}
consisting of a twistification $(\msX^\mr{tw}, \gamma)$ of $\msX$ and a $(Z, \delta)$-structure $\pi_\mcG : \mcG \migi \mfX$ on $\mfX^\mr{log}/S^\mr{log}$.
\end{itemize}
\ede
\SSP

\bde \label{y03f5D} 
 In (i)-(iii) below,  let $r$ be a nonnegative integer.
\begin{itemize}
\item[(i)]
Let 
  $S$ be  a $k$-scheme. An {\bf $r$-pointed  $(Z, \delta)$-spin curve}
  over $S$ is a collection of data
\begin{align}
\mbX := (\msX, \msX^\mr{tw}, \gamma, \pi_\mcG : \mcG \migi \mfX),
\end{align}
where $\msX$ denotes   an $r$-pointed nodal curve
 over $S$ and 
$(\msX^\mr{tw}, \gamma, \pi_\mcG : \mcG \migi \mfX)$ denotes a $(Z, \delta)$-structure on $\msX$.
\item[(ii)]
For each $j \in \{1,2\}$, let $S_j$ be a $k$-scheme and     $\mbX_j:= (\msX_j, \msX^\mr{tw}_j, \gamma_j, \pi_{\mcG, j} : \mcG_{j} \migi \mfX_j)$ (where $\msX^\mr{tw}_j := (f_j : \mfX_j \migi S_j, \{ \sigma_{\mfX_j,  i} \}_{i=1}^r)$)
   an $r$-pointed $(Z, \delta)$-spin curve
  over $S_j$.
 A {\bf $1$-morphism} (or just a {\bf morphism}) {\bf of $r$-pointed $(Z, \delta)$-spin curves} from $\mbX_1$ to $\mbX_2$ is a triple   of morphisms
 \begin{align}
\alpha := (\alpha_S, \alpha_\mfX, \alpha_\mcG)
 \end{align}
which makes 
the following diagram  $1$-commutative:
 \begin{align}
 \begin{CD}
 \mcG_{1} @> \pi_{\mcG, 1} >> \mfX_1 @> f_1 >> S_1
 \\
 @VV \alpha_\mcG V @VV \alpha_\mfX V @VV \alpha_S V
 \\ 
  \mcG_{2} @>> \pi_{\mcG, 2} > \mfX_2 @>> f_2 > S_2,
 \end{CD}
 \end{align}
 where 
 \begin{itemize}
 \item[$\bullet$]
 the right-hand square forms a morphism of $r$-pointed twisted  curves;
  \item[$\bullet$]
 the left-hand square is cartesian,  and $\alpha_\mcG$ is compatible with the respective $\widehat{Z}_\delta$-actions    on $\mcG_{1}$ and $\mcG_{2}$.
 \end{itemize}
 In particular, by taking coarse moduli spaces, one may associate, to each such morphism $\alpha : \mbX_1 \migi \mbX_2$,
  a morphism $\alpha_X : \msX_1 \migi \msX_2$ between the  underlying  $r$-pointed nodal curves.
 \item[(iii)]
 Let $S_j$,  $\mbX_j$ ($j = 1,2$) be as in (ii)  and  
  $\alpha_{l} := (\alpha_{S, l}, \alpha_{\mfX, l}, \alpha_{\mcG, l})$ ($l = 1,2$) be   morphisms $\mbX_1 \migi \mbX_2$ of $r$-pointed $(Z, \delta)$-spin curves.
 A {\bf $2$-morphism} from $\alpha_1$ to $\alpha_2$ is a triple of natural transformations 
 \begin{align}
 \mfa := (\alpha_{S,1} \stackrel{\mfa_S}{\Rightarrow} \alpha_{S,2}, \alpha_{\mfX,1} \stackrel{\mfa_\mfX}{\Rightarrow} \alpha_{\mfX,2}, \alpha_{\mcG,1} \stackrel{\mfa_\mcG}{\Rightarrow} \alpha_{\mcG,2})
 \end{align}
   compatible with each other (hence, $\mfa_S$ coincides with the identity natural transformation).
\item[(iv)]
Let $(g, r)$ be a pair of nonnegative integers with $2g-2+r >0$ and $S$ a $k$-scheme.
A {\bf pointed stable $(Z, \delta)$-spin curve of type $(g,r)$} over $S$ is an $r$-pointed $(Z, \delta)$-spin curve  over $S$ whose underlying pointed nodal curve defines    a pointed stable curve of type $(g,r)$.
\end{itemize}
\ede
\SSP

By  Definition \ref{y03f5D} above, 
$r$-pointed $(Z, \delta)$-spin curves 
 form a $2$-category.
 It is verified that 
  this $2$-category
  is equivalent to
 the ($1$-)category fibered in groupoids  over $\mfS \mfc \mfh_{/k}$ 
 whose fiber over 
  $S \in \mr{Ob} (\mfS \mfc \mfh_{/k})$ forms 
  the groupoid classifying  $r$-pointed $(Z,\delta)$-spin curves over $S$ and $2$-isomorphism classes of ($1$-)morphisms between them.
Indeed,  since all $2$-morphisms are invertible, ~\cite[Lemma 4.2.3]{AV} implies that any $1$-morphism in that  $2$-category does not  have nontrivial automorphisms.
Thus, for each pair of nonnegative integers $(g,r)$ with $2g-2+r >0$, we obtain the  ($1$-)category
\begin{align} \label{ee240}
 \mfS \mfp_{Z, \delta,  g,r}.
\end{align}
consisting of  pointed stable  $(Z,\delta)$-spin curves   of type $(g,r)$ over $k$-schemes.
The assignment $\mbX := (\msX, \msX^\mr{tw}, \gamma, \pi_\mcG : \mcG \migi \mfX)$ $\mapsto \msX$ determines a  functor
\begin{align} \label{eq009}
\pi_{g,r}^{\mfS \mfp} :  \mfS \mfp_{Z, \delta,  g,r} \migi \overline{\mfM}_{g,r}.
\end{align}

\SSP
\begin{rema} \label{RRR037}
Let $r$ be a positive integer
 and  $\msX$ (resp., $\msX^\mr{tw}$)  an $r$-pointed nodal curve (resp., an $r$-pointed twisted  curve) over a $k$-scheme $S$.
 In the subsequent  discussion, we shall refer
 to each  $2$-isomorphism class of a $1$-isomorphism of pointed $(Z, \delta)$-spin curves inducing the identity morphism of $\msX$ (resp., $\msX^\mr{tw}$) as an {\bf isomorphism of $(Z, \delta)$-structures} on  $\msX$ (resp., $\msX^\mr{tw}$). 
 In particular, we obtain the groupoid  of $(Z, \delta)$-structures  on  $\msX$ (resp., $\msX^\mr{tw}$).
\end{rema}

\SSP
\begin{rema} \label{R037} 
Certain  special cases of ``$(Z, \delta)$" may be immediately understood or  found in the previous works, as explained   as follows.
\begin{itemize}
\item[(i)]
First, let us consider the case where  $(Z, \delta) = (\mu_2, \mr{id}_{\mu_2})$.
The notion of a pointed stable $(\mu_2, \mr{id}_{\mu_2})$-spin curve is evidently equivalent to  the notion of 
a {\it twisted $2$-spin curve} in the sense of   ~\cite[\S\,1.4]{AJ1}.
In particular,  $\mfS \mfp_{\mu_2, \mr{id}_{\mu_2}, g,r}$ is equivalent  to the category $\mcB_{g,r} (\mbG_m, \omega_{\mr{log}}^{1/2})$
described in {\it loc.\,cit.}
For a general  $(Z, \delta)$, 
changing the structure group of the underlying bundles 
 via $\delta_{\mbG_m}$
gives an assignment from 
   each  twisted $2$-spin curve (i.e., a pointed stable  $(\mu_2, \mr{id}_{\mu_2})$-spin curve) to  a pointed stable  $(Z, \delta)$-spin curve.
This assignment determines a functor  
\begin{align} \label{ee701}
\mfS \mfp_{\mu_2, \mr{id}_{\mu_2}, g,r} \migi \mfS \mfp_{Z, \delta, g,r}
\end{align}
over $\overline{\mfM}_{g,r}$.
Hence, since $\mfS \mfp_{\mu_2, \mr{id}_{\mu_2}, g,r}$ is nonempty (cf. ~\cite[Corollary 4.11 and Proposition 4.19]{Chi}),  $\mfS \mfp_{Z, \delta, g,r}$ turns out  to be  nonempty.
\item[(ii)]
Next, 
 assume that  $Z = \{ 1 \}$ (hence $\delta$ is identical to  the zero map $0$ and $\widehat{Z}_\delta \cong \mbG_m$).
Then, there exists exactly  one $(\{1\}, 0)$-structure on each pointed stable curve $\msX/S$ given by the $\mbG_m$-bundle corresponding to $\Omega_{X^\mr{log}/S^\mr{log}}$.
In particular, 
$\pi_{g,r}^{\mfS \mfp}$
is an equivalence of categories.

\item[(iii)]
More generally, let us consider the case where $Z$ is arbitrary but $\delta$ coincides with  the zero map.
Then, 
there exists a natural isomorphism
$\widehat{Z}_\delta \isom   \mbG \times_k Z$ and the surjection $\nu$ (cf. (\ref{ee028})) may be identified   (relative to this isomorphism) with the first projection $\mbG_m \times_k Z \migi \mbG_m$.
It follows immedately  that $\mfS \mfp_{Z, \delta, g,r}$ is naturally isomorphic to the moduli stack $\overline{\mfM}_{g,r} (\mcB G)$ 
of
  twisted stable maps into $\mcB G$, studied by D. Abramovich, T. Graber, A. Vistoli  et al. (cf.  ~\cite{AGV0}, ~\cite{AGV}, ~\cite{ACV}, and  ~\cite{JK}). 
According to ~ \cite[Theorems 2.1.7 and  3.0.2]{ACV}, $\overline{\mfM}_{g,r} (\mcB G) \left(\cong \mfS \mfp_{Z, \delta, g,r} \right)$ may be represented by a proper smooth Deligne-Mumford stack with projective coarse moduli space.
By applying an argument similar to the argument in the proof of this result, one may obtain Theorem \ref{P08f} below, which may be thought of as  its generalization to the case of an arbitrary $(Z, \delta)$.
\end{itemize}
\end{rema}
\SSP

\bt \label{P08f}
 Let $(g, r)$ be a pair of nonnegative integers with $2g-2+r >0$.
 Then, 
 $\mfS \mfp_{Z, \delta, g,r}$ may be represented by a nonempty proper smooth   Deligne-Mumford stack over $k$ admitting a projective coarse moduli space. 
 The forgetting morphism  
 $\pi_{g,r}^{\mfS \mfp} :  \mfS \mfp_{Z, \delta,  g,r} \migi \overline{\mfM}_{g,r}$
  is finite and flat.
  Moreover, its restriction $\mfS \mfp_{Z, \delta,  g,r} \times_{\overline{\mfM}_{g,r}} \mfM_{g,r} \migi \mfM_{g,r}$ is \'{e}tale.
   \et
\begin{proof}
Since {\it the relative cotangent complex 
  of  $\mcB \nu : \mcB \widehat{Z}_\delta \migi \mcB \mbG_m$ is verified to be trivial},
the first and second assertions follow from  the arguments in ~\cite[\S\,1.5, \S\,2.1, and  \S\,2.2]{AJ1} (or the argument in the proof of ~\cite[Theorem 3.0.2]{ACV}), where $\mcB_{g,n} (\mbG_m, \omega^{1/r}_{\mr{log}})$ and  $\kappa_r : \mcB \mbG_m \migi \mcB \mbG_m$  in {\it loc.\,cit.} are  replaced with $\mfS \mfp_{Z, \delta, g,r}$ and  $\mcB \nu : \mcB \widehat{Z}_\delta \migi \mcB \mbG_m$ respectively. 
(The nonemptiness follows from the discussion in  Remark \ref{R037} (i) above.)
Moreover, the last assertion follows from ~\cite[Theorem 1.8]{O1}  and  the italicized assertion described above, which implies that 
deformations and obstructions of a pointed stable $(Z, \delta)$-spin curve are identical to those of the underlying twistification.
\end{proof}

\LSP
\subsection{Radii of extended spin structures.
} \label{y034gg}
Given a  Deligne-Mumford stack $\mfX$ of finite type  over $k$, 
we have the  {\it stack of  cyclotomic gerbes}
\begin{align}
\overline{\mcI}_{\mu} (\mfX)
\end{align}
 in $\mfX$, as described in  ~\cite[Definition 3.3.6]{AGV}.
By definition, $\overline{\mcI}_{\mu} (\mfX)$ is the disjoint union $\coprod_{l \geq 1} \overline{\mcI}_{\mu_l} (\mfX)$, where  $ \overline{\mcI}_{\mu_l} (\mfX)$ (for each positive integer $l$) denotes the category fibered in groupoids  over $\mfS \mfc \mfh_{/k}$ whose fiber over $S \in \mr{Ob} (\mfS \mfc \mfh_{/k})$ is the groupoid classifying  pairs $(\mfS, \phi)$ consisting of  a gerbe  $\mfS$ over $S$ banded by $\mu_l$  (hence $| \mfS | \cong S$) and a representable  morphism $\phi : \mfS \migi \mfX$ over $k$.

Let us consider the stack of  cyclotomic gerbes  $\overline{\mcI}_\mu (\mcB Z)$ in the classifying stack $\mcB Z$.
Denote by 
\begin{align} \label{ee702}
\mr{Inj} (\mu, Z)
\end{align}
 the set of injective morphisms of groups $\mu_l \migi Z$, where $l$ is some positive integer.
It is a finite set because  $Z$ is finite.
Given  each element  $\kappa : \mu_l \migi Z$ of $\mr{Inj} (\mu, Z)$,
we obtain an open and closed substack
\begin{align} \label{ee035}
\overline{\mcI}_\mu (\mcB Z)_\kappa
\end{align}
of $\overline{\mcI}_\mu (\mcB Z)$ classifying representable morphisms $\phi : \mfS \migi \mcB Z$ 
which are, \'{e}tale  locally on $| \mfS |$, identified with
the composite $| \mfS | \times_k \mcB \mu_l \migi \mcB Z$ of the second projection $| \mfS | \times_k \mcB \mu_l \migi \mcB \mu_l$ and the morphism $\mcB \kappa : \mcB \mu_l \migi \mcB Z$
 arising from  $\kappa$.

The stack $\overline{\mcI}_\mu (\mcB Z)$ decomposes into the disjoint union
\begin{align} \label{ee034}
\overline{\mcI}_\mu (\mcB Z) = \coprod_{\kappa \in \mr{Inj} (\mu, Z)} \overline{\mcI}_\mu (\mcB Z)_\kappa.
\end{align}
It gives  rise to   a decomposition 
\begin{align} \label{ee032}
\overline{\mcI}_\mu (\mcB Z)^{\times r} =  \coprod_{(\kappa_i )_{i=1}^r  \in \mr{Inj}(\mu, Z)^{\times r}} 
\prod_{i=1}^r  \overline{\mcI}_\mu (\mcB Z)_{\kappa_i}
\end{align}
(where $(-)^{\times r}$ denotes the product of $r$ copies of $(-)$).

\SSP
\bpr \label{Lw001}
 For each element $\kappa : \mu_l \migiincl Z$ of $\mr{Inj}(\mu, Z)$, there exists a canonical isomorphism of $k$-stacks
 \begin{align} \label{Ww900}
 \overline{\mcI}_\mu (\mcB Z)_\kappa \isom  \mcB \mr{Coker}(\kappa).
 \end{align}
  \epr
\begin{proof}
Denote by $\mcI_{\mu_l} (\mcB Z)$
  the $k$-stack classifying representable morphisms from $\mcB \mu_l$ to $\mcB Z$.
It  follows from  ~\cite[Proposition 3.4.1]{AGV} that
$\overline{\mcI}_{\mu_l} (\mcB Z)$ is canonically isomorphic to the rigidification (cf. ~\cite[Definition 5.1.9]{ACV}) of  $\mcI_{\mu_l} (\mcB Z)$ along $\mu_l$.
The representable morphisms $\mcB \mu_l \migi \mcB Z$ 
inducing, via rigidification, a morphism in $\overline{\mcI}_\mu (\mcB Z)_\kappa \left(\subseteq \overline{\mcI}_{\mu_l} (\mcB Z) \right)$
 form
a substack of $\mcI_{\mu_l} (\mcB Z)$,  which  is isomorphic to $\mcB Z$. 
Hence, 
$\overline{\mcI}_\mu (\mcB Z)_\kappa$ is isomorphic to the rigidification of  $\mcB Z$ along $\mr{Im}(\kappa) \left(\subseteq Z \right)$, namely, isomorphic to $\mcB \mr{Coker}(\kappa)$. 
\end{proof}
\SSP

Now,
let $(g, r)$ be a pair of nonnegative integers with  $2g-2+r >0$, $r >0$, and
 $\msX^\mr{tw} := (\mfX, \{ \sigma_{\mfX, i} : \mfS_i \migi \mfX \}_{i=1}^r)$
  an $r$-pointed twisted curve of genus $g$ over a $k$-scheme $S$.
  Also, let   $\mcG$  be a $(Z, \delta)$-structure on $\mfX^\mr{log}/S^\mr{log}$.
For each $i = 1, \cdots, r$,  
the restriction $\sigma^*_{\mfX, i}(\mcG)$ of $\mcG$ to $\mfS_i$ corresponds  
to a representable morphism $\mfS_i \migi \mcB \widehat{Z}_\delta$.
Notice that  the composite $\mfS_i \migi \mcB \mbG_m$  of this morphism and $\mcB \nu : \mcB \widehat{Z}_\delta \migi \mcB \mbG_m$ classifies the line bundle $\sigma^*_{\mfX, i} (\Omega_{\mfX^\mr{log}/S^\mr{log}})$,  which is  canonically identified with  the trivial line bundle  $\mcO_{\mfS_i}$
via   the residue map.
Hence,  the morphism  $\mfS_i \migi \mcB \widehat{Z}_\delta$ factors through
$\mcB Z \migi \mcB \widehat{Z}_\delta$.
If $\phi_{\mcG, i} : \mfS_i \migi \mcB Z$ denotes the resulting morphism, then we obtain
an object
\begin{align} \label{ER1}
\widetilde{\kappa}_{\mcG, i} := (\mfS_i, \phi_{\mcG, i}) \in \mr{Ob} (\overline{\mcI}_\mu (\mcB Z) (S)).
\end{align}
One may find 
 a unique  element 
\begin{align} \label{ER2}
\kappa_{\mcG, i} \in \mr{Inj} (\mu, Z)
\end{align}
such that  $\widetilde{\kappa}_{\mcG, i}$ lies in  $\mr{Ob} (\overline{\mcI}_\mu (\mcB Z)_{\kappa_{\mcG, i}} (S))$.

\SSP
\bde \label{DD011} 
 We shall refer to $\kappa_{\mcG, i}$ as the {\bf radius} of the $(Z, \delta)$-structure $\mcG$ at the marked point $\sigma_{\mfX, i}$.
  \ede
\SSP

For $r$-pointed  $(Z, \delta)$-spin curves $\mbX_j :=(\msX_j, \msX^\mr{tw}_j, \gamma_j,  \mcG_{j})$
 ($j \in \{1,2 \}$),
each morphism $\mbX_1\migi \mbX_2$ 
    induces,
  in a natural way, a morphism $\widetilde{\kappa}_{\mcG_{1}, i} \migi \widetilde{\kappa}_{\mcG_{2}, i}$ in $\overline{\mcI}_\mu (\mcB Z)$.
 Hence,
 the assignment $\mbX \mapsto \widetilde{\kappa}_{\mcG, i}$ gives rise to  a morphism of $k$-stacks
 \begin{align} \label{ee036}
\mr{ev}_{i}^{\mfS \mfp} :  \mfS \mfp_{Z, \delta, g,r} \migi \overline{\mcI}_\mu (\mcB Z).
 \end{align}
Thus, we obtain  a morphism
\begin{align} \label{ER4}
\mr{ev}^{\mfS \mfp} := (\mr{ev}^{\mfS \mfp}_1, \cdots, \mr{ev}_r^{\mfS \mfp}) :  \mfS \mfp_{Z, \delta, g,r} \migi 
\overline{\mcI}_\mu (\mcB Z)^{\times r}.
\end{align}
Each  $\vec{\kappa} := (\kappa_i)_{i=1}^r \in \mr{Inj} (\mu, Z)^{\times r}$ determines  an open and closed substack
\begin{align} \label{ee039}
\mfS \mfp_{Z, \delta, g,r, \vec{\kappa}} := (\mr{ev}^{\mfS \mfp})^{-1} (\prod_{i=1}^r \overline{\mcI}_{\mu} (\mcB Z)_{\kappa_i})
\end{align}
of $\mfS \mfp_{Z, \delta, g,r, k}$.
That is to say, $\mfS \mfp_{Z, \delta, g,r, \vec{\kappa}}$ classifies pointed stable $(Z, \delta)$-spin curves of type $(g,r)$ whose  $(Z, \delta)$-structure is of radii $\vec{\kappa}$.

\LSP
\subsection{Gluing $(Z, \delta)$-structures.} \label{y079}
In this subsection, we observe
  certain factorization properties of the moduli stack of pointed stable  $(Z, \delta)$-spin curves according to the gluing maps $\Phi_\mr{tree}$, $\Phi_\mr{loop}$  (cf. (\ref{e072}), (\ref{e071})).

Let $\mfS$ be 
 a gerbe  over a $k$-scheme $S$ banded by $\mu_l$ ($l >0$).
 One may  change the banding of the gerbe through the inversion automorphism $\zeta \mapsto \zeta^{-1}$ of $\mu_l$.
 Denote by  $\mfS^{\veebar}$  the resulting gerbe 
 and  by  $\mr{inv}_\mfS : \mfS \isom \mfS^{\veebar}$
 the isomorphism   over $S$ arising from the inversion.
The  assignment $(\mfS, \phi) \mapsto (\mfS^\veebar, \phi^\veebar)$ (where $\phi^\veebar :=  \phi \circ \mr{inv}_\mfS^{-1}$)
  induces an involution of $\overline{\mcI}_{\mu_l} (\mcB Z)$.
 By applying this argument to each piece $\overline{\mcI}_{\mu_l} (\mcB Z)$ ($l \geq 1$) of $\overline{\mcI}_\mu (\mcB Z)$ separately, we obtain an involution 
 \begin{align} \label{ee038}
 (-)^\veebar : \overline{\mcI}_\mu (\mcB Z) \isom \overline{\mcI}_\mu (\mcB Z)
 \end{align}
  of $\overline{\mcI}_\mu (\mcB Z)$.
  Each injective morphism  $\kappa : \mu_l \migiincl  Z$ in $\mr{Inj} (\mu, Z)$
  induces an injection  $\kappa^\veebar : \mu_l \migi Z$
    given by $\zeta \mapsto \kappa (\zeta^{-1})$.
  The involution $(-)^\vee$ obtained above
      restricts to an isomorphism $\overline{\mcI}_\mu (\mcB Z)_\kappa \isom \overline{\mcI}_\mu (\mcB Z)_{\kappa^\veebar}$.

Let $S$ be as above.
For each $j = 1,2$, let $(g_j, r_j)$ be a pair of nonnegative integers with $2g_j -1 + r_j >0$, 
$\msX_j$ an $(r_j+1)$-pointed stable  curve of genus $g_j$ over $S$, and
  $(\msX^\mr{tw}_j, \gamma_j)$ 
  (where $\msX^\mr{tw}_j := (\mfX_j, \{ \sigma_{j, i} : \mfS_{j, i} \migi \mfX_j \}_{i=1}^{r_j + 1})$)
  a twistification  of $\msX_j$.
  Denote by $\msX_\mr{tree}$ the $(r_1 + r_2)$-pointed stable curve of  genus $(g_1 + g_2)$ obtained  by attaching the respective last marked points of $\msX_1$ and $\msX_2$ to form a node.
Suppose that we are given an isomorphism $\epsilon : \mfS_{1, r_1 +1} \isom \mfS_{2, r_2 +1}^\veebar$
over $S$ compatible with the bands.
Then, $(\msX^\mr{tw}_1, \gamma_1)$ and $(\msX^\mr{tw}_2, \gamma_2)$ may be glued together, by means of $\epsilon$ (cf. the {\it balanced-ness} described in the fourth condition in Definition \ref{D011}, (i)),  to 
a twistification 
\begin{align}
(\msX^\mr{tw}_\mr{tree}, \gamma_\mr{tree})
\end{align}
(where $\msX^\mr{tw}_\mr{tree}:= (\mfX_{\mr{tree}}, \{ \sigma_{\mr{tree}, i}\}_{i=1}^{r_1 + r_2})$)  of  $\msX_\mr{tree}$.
 In particular, we obtain closed immersions 
 \begin{align}
 \mr{imm}_j : \mfX_j \migi \mfX_\mr{tree}
 \end{align}
  ($j =1,2$),
 which fit into the following $1$-commutative diagram:
 \begin{align} \label{ee041}
 \xymatrix{
 \mfS_{1, r_1+1} \ar[rr]^{\epsilon} \ar[rd]_{\mr{imm}_1 \circ \sigma_{1, r_1 +1}}& & \mfS^{\veebar}_{2, r_2 +1} \ar[ld]^{ \ \ \ \mr{imm}_2 \circ \sigma_{2, r_2 +1} \circ \mr{inv}^{-1}_{\mfS_{2, r_2+1}}}
 \\ & \mfX_{\mr{tree}}. & 
 }
 \end{align}

It is immediately verified that 
the morphism
\begin{align} \label{ee705}
\mr{imm}_j^* (\Omega_{\mfX^\mr{log}_\mr{tree}/S^\mr{log}}) \migi \Omega_{\mfX^\mr{log}_j/S^\mr{log}}
\end{align}
induced by
$\mr{imm}_j$ is  an isomorphism.
That is to say, there exists a canonical $2$-isomorphism
\begin{align} \label{ee706}
[\Omega_{\mfX^\mr{log}_\mr{tree}/S^\mr{log}}]\circ \mr{imm}_j \stackrel{\sim}{\Rightarrow} [\Omega_{\mfX^\mr{log}_j/S^\mr{log}}]
\end{align}
(cf. (\ref{ee030})). 
The  following lemma may be immediately verified.

\SSP
\ble \label{L01} 
 Let us keep the above notation.
 \begin{itemize}
 \item[(i)]
 Let  $\pi_\mcG : \mcG \migi \mfX_\mr{tree}$
  be a $(Z, \delta)$-structure on $\mfX_\mr{tree}^\mr{log}/S^\mr{log}$.
 Then, for each $j = 1,2$,   the pull-back  $\pi_\mcG |_{\mcG_j} : \mcG_{j}  \left(:= \mr{imm}_j^*(\mcG)\right) \migi \mfX_j$ of $\pi_\mcG$ by  $\mr{imm}_j$
 forms, via  (\ref{ee706}),  a $(Z, \delta)$-structure on $\mfX_j^\mr{log}/S^\mr{log}$.
 Moreover, by restricting $\mcG_{j}$ to $\sigma_{j, r_j +1}$ and applying the commutativity of   diagram (\ref{ee041}), we obtain  an isomorphism
 \begin{align} \label{ee040}
\epsilon_{\mcG} := (\epsilon, \epsilon_\phi) :   \widetilde{\kappa}_{\mcG_{1}, r_1+1} 
\ \Big(= (\mfS_{1, r_1 +1}, \phi_{\mcG_{1}, r_1+1})\Big) 
\isom  
\widetilde{\kappa}^\veebar_{\mcG_2, r_2+1} 
\left(= (\mfS_{2, r_2 +1}^\veebar, \phi^\veebar_{\mcG_{2}, r_2+1})\right) 
  \end{align}
 in $\overline{\mcI}_\mu (\mcB Z)$,  where $\epsilon$ is as in (\ref{ee041}) and $\epsilon_\phi$ denotes some   $2$-isomorphism $\phi_{\mcG_1, r_1+1} \stackrel{\sim}{\Rightarrow}  \phi_{\mcG_2, r_2 +1}^\veebar 
 \circ \epsilon$.
  \item[(ii)]
 Conversely, for each $j  \in \{1,2 \}$, 
 let $\pi_{\mcG, j} : \mcG_{j} \migi \mfX_j$  be a  $(Z, \delta)$-structure on $\mfX_j^\mr{log}/S^\mr{log}$.
 Also, assume that we are given 
 an isomorphism  of the form $(\epsilon, \epsilon_\phi) :  \widetilde{\kappa}_{\mcG_{1}, r_1+1}   \isom 
 \widetilde{\kappa}_{\mcG_{2}, r_2+1}^\veebar$.
 Then, there exists a unique (up to isomorphism) $(Z, \delta)$-structure $\pi_{\mcG, \mr{tree}} : \mcG_\mr{tree} \migi \mfX_\mr{tree}$ on $\mfX_\mr{tree}^\mr{log}/S^\mr{log}$ 
 together with isomorphisms  of $(Z, \delta)$-structures $\mr{imm}^*_j (\mcG_\mr{tree}) \isom \mcG_j$ ($j =1,2$) via which the equality $(\epsilon, \epsilon_\phi) = \epsilon_{\mcG_{\mr{tree}}}$  (cf. (\ref{ee040}) above) holds.
 \end{itemize}
   \ele
\SSP

By applying the above lemma, we obtain the following proposition.

\SSP
\bpr \label{P08} 
Let  $g_1$, $g_2$, $r_1$, and $r_2$ be nonnegative integers with $2g_j -1 +r_j > 0$ ($j =1, 2$).
Also,  let 
$\vec{\kappa}_1 \in \mr{Inj} (\mu, Z)^{\times r_1}$ and  $\vec{\kappa}_2 \in \mr{Inj} (\mu, Z)^{\times r_2}$.
We shall write $g := g_1 + g_2$, $r := r_1 +r_2$.
\begin{itemize}
\item[(i)]
Let $\kappa : \mu_l \migiincl Z$ ($l \geq 1$) be an element of $\mr{Inj} (\mu, Z)$.
We shall write
\begin{align}\label{WW100}
\mfS \mfp_{\mr{tree}, \kappa} :=  \mfS \mfp_{Z, \delta, g_1,r_1+1, (\vec{\kappa}_1, \kappa)}
\times_{\mr{ev}^{\mfS \mfp}_{r_1 +1}, \overline{\mcI}_\mu (\mcB Z)_{\kappa}, (-)^\veebar \circ \mr{ev}^{\mfS \mfp}_{r_2 +1}} 
\mfS \mfp_{Z, \delta, g_2,r_2+1, (\vec{\kappa}_2, \kappa^{\veebar})}.
\end{align}
Then,  
  the assignment
\begin{align} \label{Ww1001}
((\msX_1, \msX_1^\mr{tw}, \gamma_1, \mcG_1), (\msX_2, \msX_2^\mr{tw}, \gamma_2, \mcG_2)) \mapsto (\msX_\mr{tree}, \msX_{\mr{tree}}^\mr{tw}, \gamma_{\mr{tree}}, \mcG_{\mr{tree}})
\end{align}
resulting from  Lemma \ref{L01}
  determines  a morphism
\begin{align} \label{Ww1002}
\Phi^{\mfS \mfp}_{\mr{tree}, \kappa}
: \mfS \mfp_{\mr{tree},  \kappa}  \migi \mfS \mfp_{Z, \delta, g,r, (\vec{\kappa}_1, \vec{\kappa}_2)},
\end{align}
which makes the following square diagram commute:
\begin{align} \label{Ww1003}
\begin{CD} 
\mfS \mfp_{\mr{tree},  \kappa}
@>  \Phi^{\mfS \mfp}_{\mr{tree}, \kappa} >> 
\mfS \mfp_{Z, \delta, g,r, (\vec{\kappa}_1, \vec{\kappa}_2)} 
\\
@V  \pi^{\mfS \mfp}_{g_1,r_1} \times \pi^{\mfS \mfp}_{g_2,r_2} VV @V \pi^{\mfS \mfp}_{g, r} VV
\\
\overline{\mfM}_{g_1, r_1+1} \times_k \overline{\mfM}_{g_2, r_2+1} @>> \Phi_{\mr{tree}}> \overline{\mfM}_{g, r}.
\end{CD}
\end{align}
Moreover, the ramification index of $\pi^{\mfS \mfp}_{g, r}$ at the image of $\Phi^{\mfS \mfp}_{\mr{tree}, \kappa}$ coincides with $l$.
\item[(ii)]
 We shall write
 \begin{align} \label{Ww101}
 \mfW_\mr{tree} := (\overline{\mfM}_{g_1, r_1+1} \times_k \overline{\mfM}_{g_2, r_2+1}) \times_{\Phi_\mr{tree}, \overline{\mfM}_{g,r}, \pi^{\mfS \mfp}_{g,r}} \mfS \mfp_{Z, \delta, g,r, (\vec{\kappa}_1, \vec{\kappa}_2)}
 \end{align}
  and write $\mfW^\mr{red}_\mr{tree}$ for the reduced stack  associated to $\mfW_\mr{tree}$.
  Then, the morphism
  \begin{align} \label{WW102}
  \coprod_{\kappa \in \mr{Inj}(\mu, Z)} \mfS \mfp_{\mr{tree}, \kappa} \migi \mfW^\mr{red}_\mr{tree}
  \end{align}
induced by  $\coprod_{\kappa \in \mr{Inj}(\mu, Z)} \Phi^{\mfS \mfp}_{\mr{tree}, \kappa}$ (because of the commutativity of (\ref{Ww1003}))  is an isomorphism.
 \item[(iii)]
 The following equality holds:
\begin{align} \label{Ww103}
\mr{deg}(\mfS \mfp_{Z, \delta, g,r, (\vec{\kappa}_1, \vec{\kappa}_2)}/\overline{\mfM}_{g,r}) 
= |Z| \cdot \sum_{\kappa \in \mr{Inj}(\mu, Z)} \prod_{j=1}^2 \mr{deg}(\mfS \mfp_{Z, \delta, g_j,r_j+1, (\vec{\kappa}_j, \kappa)}/\overline{\mfM}_{g_j, r_j+1}).
\end{align}
\end{itemize}
   \epr
\begin{proof}
First, let us prove assertion (i).
The nontrivial portion is the last statement.
Let $q$ be a geometric point in $\mfS \mfp_{Z, \delta, g, r, (\vec{\kappa}_1, \vec{\kappa}_2)}$ determined by  the image via  $\Phi_{\mr{tree}, \kappa}^{\mfS \mfp}$ of an irreducible component of $\mfS \mfp_{\mr{tree}, \kappa}$.
Denote by $\overline{q} \in \overline{\mfM}_{g,r}$ the image of $q$ via $\pi_{g,r}^{\mfS \mfp}$.
Also, denote by $\msX$ the pointed stable curve classified by $\overline{q}$ 
and by $(\msX^\mr{tw}, \gamma, \mcG)$ the
$(Z, \delta)$-structure on $\msX$ classified by $q$.
 Since the relative cotangent complex of $\mcB \nu : \mcB \widehat{Z}_\delta \migi \mcB \mbG_m$ is trivial, the deformations of 
 the 
  $(Z, \delta)$-structure  $(\msX^\mr{tw}, \gamma, \mcG_{Z, \delta})$
 are identical to those of the underlying twistification
 $(\msX^\mr{tw}, \gamma)$
  (cf. the proof in Theorem \ref{P08f}).
Let $D := \mr{Spec}(R)$ (for some $k$-algebra $R$) be  the versal deformation space of  $\overline{q} \in \overline{\mfM}_{g,r}$.
Then, by ~\cite[Theorem 4.5]{Chi} (or the discussion in the proof of Theorem 4.4 in {\it loc.\,cit}.), 
the deformation space of the twistification $(\msX^\mr{tw}, \gamma)$  may be represented  by 
the quotient stack $\widetilde{D} := [\mr{Spec}(R[z]/(z^l-w))/\mu_l]$ over $D$, where $w  \left(\in R\right)$ denotes a function defining $\overline{q}$ and the $\mu_l$-action on $\mr{Spec}(R [z]/(z^l-w))$ is given by $(\zeta, z) \mapsto \zeta  \cdot z$ (for each $\zeta \in \mu_l$).
(To be precise, we apply the {\it pointed curve} version of the result in {\it loc.\,cit.}, which may be formulated and proved similarly.)
In particular,  the ramification index of $\widetilde{D}/D$ at 
the point classifying  $(\msX^\mr{tw}, \gamma)$  coincides with $l$.
This completes the proof of assertion (i).

Next, assertion (ii) follows from the local description (\ref{Ww131}) of 
a twisted curve at a node.
Indeed, if $R$ is a reduced ring,
then the \'{e}tale stack  over $\mr{Spec}(R[s_0, t_0]/(s_0t_0))$ of the form (\ref{Ww131}) must satisfy the equality  $u_1 =0$.
This  stack  may be obtained by gluing together the two stacks $[\mr{Spec}(A[s_1]/(s_1^l-s_0))/\mu_l]$ and $[\mr{Spec}(A[t_1]/(t_1^l-t_0))/\mu_l]$.
This implies 
that the pointed $(Z, \delta)$-spin curves classified by  $\mfW^\mr{red}_\mr{tree}$ come from $\coprod_{\kappa \in \mr{Inj}(\mu, Z)} \mfS \mfp_{\mr{tree},\kappa}$. Hence, (\ref{WW102}) turns out to be an isomorphism.

Finally, let us consider assertion (iii).
Observe that the  following equalities hold:
\begin{align} \label{Ww120}
& \ \mr{deg}(\mfS \mfp_{\mr{tree},\kappa}/(\overline{\mfM}_{g_1, r_1+1} \times_k \overline{\mfM}_{g_2, r_2+1})) \\
= & \ |\mr{Coker}(\kappa)| \cdot \prod_{j=1}^2 \mr{deg}(\mfS \mfp_{Z, \delta, g_j, r_j+1, (\vec{\kappa}_j, \kappa)}/\overline{\mfM}_{g_j, r_j+1}) \notag \\
= & \ \frac{|Z|}{l} \cdot \prod_{j=1}^2 \mr{deg}(\mfS \mfp_{Z, \delta, g_j, r_j+1, (\vec{\kappa}_j, \kappa)}/\overline{\mfM}_{g_j, r_j+1}). \notag
\end{align}
Given an element $\kappa : \mu_l \migiincl Z$ of $\mr{Inj}(\mu, Z)$, we shall write
$\mfW_{\mr{tree}, \kappa}$ for the open and closed substack of $\mfW_{\mr{tree}}$ determined by the image via (\ref{WW102}) of 
$\mfS \mfp_{\mr{tree},  \kappa} \left(\subseteq \coprod_{\kappa \in \mr{Inj}(\mu, Z)} \mfS \mfp_{\mr{tree},\kappa} \right)$.
Also, write $\mfW_{\mr{tree}, \kappa}^\mr{red}$ for the reduced stack associated to $\mfW_{\mr{tree}, \kappa}$ (hence, $\mfW_{\mr{tree}}^\mr{red} = \coprod_{\kappa \in \mr{Inj}(\mu, Z)} \mcN_{\mr{tree}, \kappa}^\mr{red}$).
By  the last assertion of  (ii), we have
\begin{align} \label{Ww300}
\mr{deg} (\mfW_{\mr{tree}, \kappa}/ (\overline{\mfM}_{g_1, r_1+1} \times_k \overline{\mfM}_{g_2, r_2 +1})) = l \cdot 
\mr{deg} (\mfW_{\mr{tree}, \kappa}^\mr{red}/ (\overline{\mfM}_{g_1, r_1+1} \times_k \overline{\mfM}_{g_2, r_2 +1})). 
\end{align}
Thus, the following sequence of equalities holds:
\begin{align}
& \ \ \  \ \mr{deg}(\mfS \mfp_{Z, \delta, g, r, (\vec{\kappa}_1, \vec{\kappa}_2)}/\overline{\mfM}_{g,r}) \\
& =
\mr{deg}(\mfW_\mr{tree}/ (\overline{\mfM}_{g_1, r_1+1} \times_k \overline{\mfM}_{g_2, r_2 +1}))
\notag \\
& = \sum_{\kappa \in \mr{Inj}(\mu, Z)} \mr{deg}(\mfW_{\mr{tree}, \kappa}/ (\overline{\mfM}_{g_1, r_1+1} \times_k \overline{\mfM}_{g_2, r_2 +1})) \notag  \\
& \stackrel{(\ref{Ww300})}{=} \sum_{\kappa \in \mr{Inj}(\mu, Z)} l \cdot\mr{deg}(\mfW_{\mr{tree}, \kappa}^\mr{red}/ (\overline{\mfM}_{g_1, r_1+1} \times_k \overline{\mfM}_{g_2, r_2 +1})) \notag \\
& \stackrel{\text{(ii)}}{=}
 \sum_{\kappa \in \mr{Inj}(\mu, Z)} 
l \cdot \mr{deg}(\mfS \mfp_{\mr{tree},  \kappa}/(\overline{\mfM}_{g_1, r_1+1} \times_k \overline{\mfM}_{g_2, r_2 +1})) 
\notag \\
& \stackrel{(\ref{Ww120})}{=}  |Z| \cdot \sum_{\kappa \in \mr{Inj}(\mu, Z)} \prod_{j=1}^2 \mr{deg}(\mfS \mfp_{Z, \delta, g_j,r_j+1, (\vec{\kappa}_j, \kappa)}/\overline{\mfM}_{g_j, r_j+1}), \notag
\end{align}
where the first equality follows from the fact that $\pi_{g,r}^{\mfS \mfp}$ is finite and flat (cf. Theorem \ref{P08f}).
This completes the proof of assertion (iii).
\end{proof}
\SSP

Next, let $(g, r)$ be a pair of  nonnegative integers with $2g+r >0$, $\msX$ an $(r+2)$-pointed stable  curve of genus $g$ over $S$, and $(\msX^\mr{tw}, \gamma)$ (where $\msX^\mr{tw} := (\mfX, \{ \sigma_{\mfX, i} : \mfS_i \migi \mfX \}_{i=1}^{r+2})$) a twistification  of $\msX$.
Denote by $\msX_\mr{loop}$ the $r$-pointed stable curve of genus $g$ obtained by attaching the last two marked points of $\msX$ to form a node.
Suppose that we are given an isomorphism
$\epsilon : \mfS_{r+1} \isom \mfS_{r +2}^{\veebar}$ over $S$ compatible with the bands.
 By attaching the last two marked points of $\msX^\mr{tw}$ along $\epsilon$, 
we obtain
a twistification  $(\msX^\mr{tw}_{\mr{loop}}, \gamma_{\mr{loop}})$ of 
$\msX_\mr{loop}$.

There exists 
  a natural bijective correspondence $\mcG \mapsto \mcG_\mr{loop}$  between    the  set of isomorphism classes of  $(Z, \delta)$-structures $\mcG$ on $\mfX^\mr{log}/S^\mr{log}$ admitting an isomorphism
  $\widetilde{\alpha}_{\mcG_{}, r+1} \isom \widetilde{\alpha}_{\mcG_{}, r+2}^\veebar$ compatible with  $\epsilon$  and the set of isomorphism classes of $(Z, \delta)$-structures $\mcG_\mr{loop}$ on $\mfX_\mr{loop}^\mr{log}/S^\mr{log}$.
 Thus, we obtain  the following propositions.

\SSP
\bpr \label{Pw081} 
 Let $(g, r)$ be a pair nonnegative integers with $2g + r >0$, and let $\vec{\kappa} \in \mr{Inj} (\mu, Z)^{\times r}$.
\begin{itemize}
 \item[(i)]
Let $\kappa : \mu_l \migiincl Z$ ($l \geq 1$) be an element of $\mr{Inj} (\mu, Z)$.
We shall write
\begin{align} \label{Ww109}
\mfS \mfp_{\mr{loop}, \kappa} := \mr{Ker}\left(\mfS \mfp_{Z, \delta, g,r+2, (\vec{\kappa}, \kappa, \kappa^{\veebar})} \stackrel{\mr{ev}^{\mfS \mfp}_{r+1}, (-)^\veebar \circ  \mr{ev}^{\mfS \mfp}_{r+2}}{\rightrightarrows} \overline{\mcI}_\mu (\mcB Z)_{\kappa}\right).
\end{align}
 Then, 
  the assignment  
 \begin{align} \label{Ww1004}
 (\msX, \msX^\mr{tw}, \gamma, \mcG) \mapsto (\msX_\mr{loop}, \msX^\mr{tw}_{\mr{loop}}, \gamma_{\mr{loop}}, \mcG_{\mr{loop}})
 \end{align}
resulting from the above discussion 
  determines a morphism
 \begin{align} \label{Ww1005}
 \Phi^{\mfS \mfp}_{\mr{loop}, \kappa} :  \mfS \mfp_{\mr{loop},  \kappa} \migi \mfS \mfp_{Z, \delta, g+1,r, \vec{\kappa}, k},
 \end{align}
and the following square diagram is  commutative:
 \begin{align} \label{Ww1006}
 \begin{CD}
 \mfS \mfp_{\mr{loop},  \kappa}
    @> 
     \Phi_{\mr{loop}, \kappa}^{\mfS \mfp} >> \mfS \mfp_{Z, \delta, g+1,r, \vec{\kappa}} 
 \\
 @V 
  \pi^{\mfS \mfp}_{g,r+2} VV @V \pi^{\mfS \mfp}_{g+1, r} VV
 \\
 \overline{\mfM}_{g, r+2} @>> \Phi_{\mr{loop}}> \overline{\mfM}_{g+1, r}.
\end{CD}
\end{align}
Moreover, the ramification index of $\pi^{\mfS \mfp}_{g+1, r}$ at the image of $\Phi_{\mr{loop}, \kappa}^{\mfS \mfp}$ coincides with $l$.
\item[(ii)]
 We shall write
 \begin{align} \label{Ww110}
 \mfW_\mr{loop} := \overline{\mfM}_{g, r+2} \times_{\Phi_\mr{loop}, \overline{\mfM}_{g+1,r}, \pi^{\mfS \mfp}_{g+1,r}} \mfS \mfp_{Z, \delta, g+1,r, \vec{\kappa}}
 \end{align}
  and write $\mfW^\mr{red}_\mr{loop}$ for the reduced stack  associated to $\mfW_\mr{loop}$.
  Then, the morphism
  \begin{align} \label{WW111}
  \coprod_{\kappa \in \mr{Inj}(\mu, Z)} 
  \mfS \mfp_{\mr{loop},  \kappa}  \migi \mfW^\mr{red}_\mr{loop}
  \end{align}
induced by  $\coprod_{\kappa \in \mr{Inj}(\mu, Z)} \Phi^{\mfS \mfp}_{\mr{loop}, \kappa}$ (because of the commutativity of (\ref{Ww1006}))  is an isomorphism.

\item[(iii)]
 The following equality holds:
\begin{align} \label{Ww112}
\mr{deg}(\mfS \mfp_{Z, \delta, g+1,r, \vec{\kappa}}/\overline{\mfM}_{g+1,r}) 
= |Z| \cdot \sum_{\kappa \in \mr{Inj}(\mu, Z)}\mr{deg}(\mfS \mfp_{Z, \delta, g, r+2, (\vec{\kappa}, \kappa, \kappa^\veebar)}/\overline{\mfM}_{g, r+2}).
\end{align}
\end{itemize}
   \epr
\begin{proof}
The assertions follow from arguments similar to the arguments in the proof of  Proposition \ref{P08}.
\end{proof}

\LSP
\subsection{Forgetting tails.} \label{y040}
Next, we observe  the  factorization of $\mfS \mfp_{Z, \delta, g,r, \vec{\kappa}}$ according to the forgetting-tails map $\Phi_\mr{tail}$ (cf. (\ref{e070})).
Let $r$ be a nonnegative integer, 
$S$ a $k$-scheme,  and 
$\msX := (X, \{ \sigma_{X, i} \}_{i=1}^{r+1})$  an $(r+1)$-pointed {\it smooth} curve over $S$. 
Write $\msX_\mr{tail} := (X, \{ \sigma_{X, i} \}_{i=1}^{r})$, i.e.,
the $r$-pointed curve
 obtained from $\msX$ by forgetting the last marked point. 
Moreover, let $(\msX_\mr{tail}^\mr{tw}, \gamma_\mr{tail})$ (where $\msX^\mr{tw}_\mr{tail} := (\mfX_\mr{tail}, \{ \sigma_{\mfX_\mr{tail}, i} \}_{i=1}^r)$) be  a twistification of $\msX_\mr{tail}$.

\SSP
\ble \label{Pfg08} 
There exists a unique (up to isomorphism) twistification $(\msX_{+ \delta}^\mr{tw}, \gamma_{+ \delta})$ of $\msX$ which is isomorphic to $(\msX_\mr{tail}^\mr{tw}, \gamma_\mr{tail})$ when restricted to $X \setminus \mr{Im} (\sigma_{X, r+1})$ and such that the stabilizer at any  geometric point in $\mfX_{+ \delta}$ (:= the underlying curve of $\msX_{+ \delta}^\mr{tw}$)  lying over  $\mr{Im} (\sigma_{X, r+1})$ ($\subseteq X$)  has order $2$. 
   \ele
\begin{proof}
Let us construct the desired 
twistification. 
Consider the category $X [\mr{Im}(\sigma_{X, r+1})/2]$ (cf. ~\cite[Definition 2.2]{Chi}) consisting of collections  $(S, M, j, s)$, where
\begin{itemize}
\item[$\bullet$]
$S$ is an $X$-scheme $S \migi X$;
\item[$\bullet$]
$M$ is a line bundle on $S$;
\item[$\bullet$]
$j$ is an isomorphism between $M^{\otimes 2}$ and the pull-back of $\mcO_X (\mr{Im}(\sigma_{X, r+1}))$ on $S$;
\item[$\bullet$]
$s$ is a global section of $M$ such that $j (s^{\otimes 2})$ equals the tautological section of $\mcO_X (\mr{Im}(\sigma_{X, r+1}))$ vanishing along $\sigma_{X, r+1}$.
\end{itemize}
The morphisms in this category are  defined in an obvious way.
$X [\mr{Im}(\sigma_{X, r+1})/2]$  may be represented by a  Deligne-Mumford stack over $k$ with coarse moduli space $X$.
The projection $\mr{coa} : X [\mr{Im}(\sigma_{X, r+1})/2] \migi X$ (i.e., $\mr{coa} := \mr{coa}_{X [\mr{Im}(\sigma_{X, r+1})/2]}$) is an isomorphism over $X \setminus \mr{Im} (\sigma_{X, r+1})$.
Also, $X [\mr{Im}(\sigma_{X, r+1})/2]$ is equipped with a tautological line bundle $\mcM$ and an isomorphism $\mcM^{\otimes 2} \isom \mr{coa}^*(\mcO_X (\mr{Im}(\sigma_{X, r+1})))$.
Let us write   $\mfX_{+\delta} := \mfX_{\mr{tail}} \times_X X [\mr{Im}(\sigma_{X, r+1})/2]$ and write $\xi : \mfX_{+ \delta} \migi \mfX_\mr{tail}$ for the projection to the first  factor.
There exists  a unique closed immersion $\sigma_{\mfX_{+\delta}, r+1} : \mfS_{r+1} \migi \mfX_{+\delta}$ (where $\mfS_{r+1}$ is an \'{e}tale gerbe over $S$ banded by $\mu_2$) inducing the closed immersion $\sigma_{X, r+1}$ between their respective  coarse moduli spaces.
Since  $\xi$ restricts to an isomorphism $\mfX_{+\delta} \setminus \mr{Im}(\sigma_{\mfX_{+\delta}, r+1}) \isom \mfX_{\mr{tail}} \setminus \mr{Im} (\xi \circ \sigma_{\mfX_{+ \delta}, r+1})$,
each $\sigma_{X, i}$ ($i=1, \cdots, r$) may be regarded as a marked point $\sigma_{\mfX_{+\delta}, i}$ of $\mfX_{+\delta}$.
By setting $\msX^\mr{tw}_{+\delta} := (\mfX_{+\delta}, \{ \sigma_{\mfX_{+\delta}, i} \}_{i=1}^r)$ and $\gamma_{+ \delta} := \gamma_\mr{tail} \circ \xi$,
we obtain  a pair  $(\msX^\mr{tw}_{+\delta}, \gamma_{+ \delta})$  forming  the desired  twistification of $\msX$.

Next, let us consider the uniqueness assertion.
Let $(\msX^\mr{tw}, \gamma)$ (where $\msX^\mr{tw} := (\mfX, \{ \sigma_{\mfX, i} \}_{i=1}^{r+1})$) be a twistification of $\msX$ satisfying the required  properties.
The line bundle $\mcO_{\mfX} (\mr{Im}(\sigma_{\mfX, r+1}))^{\otimes 2}$ on $\mfX$ descends to 
$\mcO_{X} (\mr{Im} (\sigma_{X, r+ 1}))$.
Hence,   the line bundle  $\mcO_{\mfX} (\mr{Im}(\sigma_{\mfX, r+1}))$
 and  its tautological section vanishing along $\sigma_{\mfX, r+1}$ induce a morphism $\mfX \migi X [\mr{Im}(\sigma_{X, r+1})/2]$. This morphism extends naturally  to a morphism $\mfX \migi \mfX_{+ \delta}$ ($= \mfX_{\mr{tail}} \times_X X [\mr{Im}(\sigma_{X, r+1})/2]$); it specifies, by construction,  an isomorphism of twistifications $(\msX^\mr{tw}, \gamma) \isom (\msX_{+ \delta}^\mr{tw}, \gamma_{+ \delta})$.
This completes the proof of the uniqueness assertion. 
\end{proof}
\SSP

 Next, let  $\mcG_\mr{tail}$ be  a $(Z, \delta)$-structure on $\mfX_\mr{tail}^\mr{log}/S$.
By using this,
we shall  construct a $(Z, \delta)$-structure on $\mfX_{+\delta}^\mr{log}/S$.
Observe that 
\begin{align}
\Omega_{\mfX^\mr{log}_{+\delta}/S} & \cong \xi^*(\Omega_{\mfX_\mr{tail}^\mr{log}/S^\mr{log}}) \otimes (\mr{coa}_{\mfX_\mr{tail}} \circ \xi)^*(\mcO_X (\mr{Im}(\sigma_{X, r+1}))) \\
& \cong \xi^*(\Omega_{\mfX_\mr{tail}^\mr{log}/S}) \otimes (\mr{coa} \circ \mr{coa}')^* (\mcO_X (\mr{Im}(\sigma_{X, r+1})))  \notag \\
& \cong \xi^*(\Omega_{\mfX_\mr{tail}^\mr{log}/S}) \otimes  \mr{coa}'^*(\mcM)^{\otimes 2}, \notag
\end{align}
where $\mr{coa}'$ denotes the projection $\mfX_{+\delta} \migi X [\mr{Im}(\sigma_{X, r+1})/2]$.
Hence, by twisting $\mcG_\mr{tail}$ by the $\mbG_m$-bundle corresponding to $\mr{coa}'^*(\mcM)$,
we obtain a $\widehat{Z}_\delta$-bundle $\mcG_{+ \delta}$,
 which 
 forms a $(Z, \delta)$-structure on $\mfX_{+\delta}^\mr{log}/S$.
Since the nontrivial automorphism of $\mfS_{r+1}$ over $S$ arises from the automorphism of $\mcM$ given by multiplication by $(-1) \in \mu_2$, the commutativity of the left-hand square in (\ref{ee028}) implies that  the
 radius  of $\mcG_{+\delta}$ at $\sigma_{\mfX_{+\delta}, r+1}$ coincides with $\overline{\delta} : \mu_{\overline{2}} \migiincl Z$ ($\subseteq \widehat{Z}_\delta$).
Thus, 
we have obtained a $(Z, \delta)$-structure
\begin{align} \label{ee730}
(\msX_{+\delta}^\mr{tw}, \gamma_{+\delta}, \mcG_{+\delta})
\end{align}
on $\msX$ with $\kappa_{\mcG_{+\delta}, r+1} = \overline{\delta}$.

Conversely, suppose that we are given a $(Z, \delta)$-structure $(\msX^\mr{tw}, \gamma, \mcG)$ on $\msX$ 
with $\kappa_{\mcG, r+1} = \overline{\delta}$.
One may find  
  a unique (up to isomorphism)
twistification $(\msX_{-\delta}^\mr{tw}, \gamma_{-\delta})$ (where $\msX^\mr{tw}_{-\delta} := (\mfX_{-\delta}, \{ \sigma_{\mfX_{-\delta}, i} \}_{i=1}^r)$) of $\msX_\mr{tail}$
such that there exists an isomorphism $\mfX \isom \mfX_{-\delta} \times_X X [\mr{Im} (\sigma_{X, r+1})/2]$ whose composite with $\gamma_{-\delta} \times \mr{cos} : \mfX_{-\delta} \times_X X [\mr{Im} (\sigma_{X, r+1})/2] \migi X$ coincides with $\gamma$.
If $\mcG \cdot \mcM$ denotes the twist of $\mcG$ by the pull-back of $\mcM$ via the projection $\mfX \migi X [\mr{Im} (\sigma_{X, r+1})/2]$, 
then $2 \overline{\delta} = 0 \ (\in \mr{Inj}(\mu, Z))$ implies that  the automorphism group of $\mfS_{r+1}$ acts on $\sigma_{\mfX, r+1}^* (\mcG \cdot \mcM)$  trivially.
Hence, $\mcG \cdot \mcM$ comes
from
   a (unique up to isomorphism) $(Z, \delta)$-structure $\mcG_{- \delta}$ on $\mfX_{-\delta}^\mr{log}/S$ via pull-back by the projection $\mfX \migi \mfX_{-\delta}$.
Consequently, we have obtained  the following proposition.

\SSP
\bpr \label{P020} 
 \begin{itemize}
 \item[(i)]
 The assignments $(\msX^\mr{tw}_\mr{tail}, \gamma_\mr{tail},  \mcG) \mapsto (\msX^\mr{tw}_{+\delta}, \gamma_{+\delta},  \mcG_{+ \delta})$ and  $(\msX^\mr{tw}, \gamma, \mcG) \mapsto (\msX_{-\delta}^\mr{tw}, \gamma_{-\delta},  \mcG_{- \delta})$ constructed  above  determine  an equivalence of categories between the groupoid  of $(Z, \delta)$-structures on $\msX_\mr{tail}$ and the groupoid  of $(Z, \delta)$-structures on $\msX$ whose radius at $\sigma_{\mfX, r+1}$ coincides with  $\overline{\delta}$.
\item[(ii)]
   Let $(g, r)$ be  a pair of nonnegative integers with
    $2g -1+r >0$, and let $\vec{\kappa} \in \mr{Inj} (\mu, Z)^{\times r}$.
 Then, the assignment 
 \begin{align}
 (\msX, (\msX_\mr{tail}, \msX_\mr{tail}^\mr{tw}, \gamma_\mr{tail}, \mcG)) \mapsto  (\msX, (\msX, \msX^\mr{tw}, \gamma, \mcG_{+ \delta}))
 \end{align}
    determines an isomorphism
 \begin{align}
\mfM_{g, r+1} \times_{\overline{\mfM}_{g,r}} \mfS \mfp_{Z, \delta, g,r, \vec{\kappa}}  \isom \mfM_{g, r+1} \times_{\overline{\mfM}_{g,r+1}} \mfS \mfp_{Z, \delta, g,r+1, (\vec{\kappa}, \overline{\delta})} 
 \end{align}
over $\mfM_{g, r+1}$.
\item[(iii)]
The following equality holds:
\begin{align}
\mr{deg}(\mfS \mfp_{Z, \delta, g,r+1, (\vec{\kappa}, \overline{\delta})}/\overline{\mfM}_{g,r+1}) = \mr{deg}(\mfS \mfp_{Z, \delta, g,r, \vec{\kappa}}/\overline{\mfM}_{g,r}).
\end{align}
\end{itemize}
  \epr

\LSP
\subsection{Automorphisms of $(Z, \delta)$-structures} \label{Whhh034}
In this subsection, suppose that $k$ is algebraically closed.
Let $V$ be a finite set and $\{ \msX_v \}_{v \in V}$  a collection of pointed  {\it irreducible} stable curves over $k$ indexed by $V$, where each $\msX_v$ is $r_v$-pointed  and of  genus $g_r$ with $2g_v -2 +r_v >0$.
Also, 
let  $\msX := (X, \{ \sigma_{X, i} \}_{i=1}^r)$ be a pointed  stable curve over $k$ obtained  from the $\msX_v$'s  by attaching some of their marked points to form nodes (cf. ~\cite[\S\,7.3]{Wak5}).
In particular, $V$ may be identified with 
the set of vertices of the dual graph $\Gamma$  of $\msX$ (cf.  ~\cite[\S\,1.5]{JKV}).
Denote by $E$ the set of  edges of $\Gamma$. 
Given an element $e$ of $E$, 
we  denote by $e_1$, $e_2$  the two half-edges belonging to $e$ (i.e., the branches of the node corresponding to $e$) and  
  by $v_{e, j}$ (for each $j \in \{1, 2\}$) the element of $V$
at which $e_j$ is attached.
Moreover, for any  $e \in E$ and $j \in \{1, 2\}$, 
 denote by $i_{e,j}$ the unique element of $\{1, \cdots, r_{v_{e, j}} \}$ such that the $i_{e, j}$-th marked point  of $\msX_{v_{e, j}}$ corresponds to  $e_{j}$.

Suppose further that, for each $v \in V$, we are given
 a $(Z, \delta)$-structure $(\msX_v^\mr{tw}, \gamma_v, \mcG_{v})$ on $\msX_v$ 
 admitting 
 an isomorphism  $\widetilde{\kappa}_{\mcG_{v_{e,1}}, i_{e, 1}} \isom \widetilde{\kappa}_{\mcG_{v_{e, 2}}, i_{e, 2}}^\veebar$ for any $e \in E$.
 In particular, we have  the equality 
 \begin{align} \label{Ww911}
 \kappa_e := \kappa_{\mcG_{v_{e,1}}, i_{e, 1}} = \kappa_{\mcG_{v_{e,2}}, i_{e, 2}}^\veebar \in \mr{Inj}(\mu, Z)
 \end{align}
for every $e \in E$.
Then, 
one may  apply Lemma \ref{L01}, (ii),
successively  to obtain  a $(Z, \delta)$-structure $(\msX^\mr{tw}, \gamma, \mcG)$ on $\msX$.
In particular, $(\msX^\mr{tw}, \gamma, \mcG)$ becomes $(\msX_v^\mr{tw}, \gamma_v, \mcG_v)$ ($v \in V$) after pull-back to $\msX_v$. 
Let us write
\begin{align}
\mbX := (\msX, \msX^\mr{tw}, \gamma, \mcG), \hspace{5mm} \mbX_v := (\msX_v, \msX^\mr{tw}_v, \gamma_v, \mcG_v)
\end{align} 
($v \in V$).
Then, we obtain the following proposition, which may be found 
 in ~\cite[Proposition 1.18]{JKV} (and the comment following that proposition) when $(Z, \delta) = (\mu_2, \mr{id}_{\mu_2})$.

\SSP
\bpr \label{L010}
Let us keep the above notation.
Denote by $\mr{Aut}_{\msX} (\mbX)$ (resp., $\mr{Aut}_{\msX_v} (\mbX_v)$ for each $v \in V$) the automorphism group of $\mbX$ (resp., $\mbX_v$) inducing the identity of $\msX$ (resp., $\msX_v$).
\begin{itemize}
\item[(i)]
For each $v \in V$, 
$\mr{Aut}_{\msX_v} (\mbX_v)$ is canonically  isomorphic to $Z$.
\item[(ii)]
The group $\mr{Aut}_{\msX} (\mbX)$
fits into the following exact sequence:
\begin{align}
1 \longmigi \mr{Aut}_{\msX}(\mbX) \xrightarrow{ \ \alpha \ } \prod_{v \in V} \mr{Aut}_{\msX_v}(\mbX_v) \xrightarrow{ \ \beta \ } \prod_{e \in E_\mr{nl}} \mr{Coker}(\kappa_e),
\end{align}
\end{itemize}
where  
\begin{itemize}
\item
$\alpha$ denotes the morphism induced  by  pull-back to  $\msX_v$ ($v \in V$);
\item
$E_\mr{nl}$ denotes the set of edges which do not start and end at the same vertex (i.e., non-loops);
\item
 $\beta$ maps 
each element $(\zeta_v)_{v \in V} \in \prod_{v \in V} \mr{Aut}_{\msX_v}(\mbX_v)  \left( = \prod_{v \in V} Z \right)$ to
$(\zeta_{i_{e, 1}}\zeta_{i_{e, 2}}^{-1} \  \mr{mod}  \ \mr{Im}(\kappa_e))_{e \in E_\mr{nl}}$. 
\end{itemize}
   \epr
\begin{proof}
Assertion (i) follows from the fact that 
any element of  $\mr{Aut}_{\msX_v} (\mbX_v)$
may be expressed uniquely  as 
 the automorphism given by  translation  by an element  of $Z$.
Assertion (ii) follows immediately from Proposition \ref{Lw001} and  Lemma \ref{L01}.
\end{proof}

\LSP
\subsection{$(Z, \delta)$-structures on the  $2$-pointed projective line} \label{yhhh034}
Let $S$ be  a $k$-scheme, and 
let  $\mbP_S := \mcP roj (\mcO_S [x, y])$ denote the projective line over $S$, i.e., the moduli space classifying ratios $[x : y]$ (for $x$, $y \in \mcO_S$ with $(x, y) \neq (0, 0)$).
Denote by  $\sigma_{\mbP, 1}$
  and  $\sigma_{\mbP, 2}$
   the marked points of $\mbP_S$ determined by the values $0$ ($= [0 : 1]$) and   $\infty$ ($= [1: 0]$) respectively.
In particular, we have a $2$-pointed curve $\msP'_S := (\mbP_S,\{ \sigma_{\mbP, 1}, \sigma_{\mbP, 2}\})$ over $S$.
The $S$-scheme $\mbP_S$ has two open subschemes $U_1 := \mbP_S \setminus \mr{Im} (\sigma_{\mbP, 2}) = \mcS pec (\mcO_S [s_1])$ (where $s_1 := x/y$) and  $U_2 :=\mbP_S \setminus  \mr{Im} (\sigma_{\mbP, 1}) = \mcS pec (\mcO_S [s_2])$ (where $s_2 := y/x$).

\SSP
\bpr \label{P08ff}
\begin{itemize}
\item[(i)]
Let $(\msP^\mr{tw}, \gamma, \mcG)$ be a $(Z, \delta)$-structure on $\msP_S$.
Then, the equality $\kappa_{\mcG, 1} = \kappa^\veebar_{\mcG, 2}$ ($\in \mr{Inj} (\mu, Z)$)  holds.
\item[(ii)]
For each $\kappa \in \mr{Inj} (\mu, Z)$, there exists a unique (up to isomorphism) $2$-pointed $(Z, \delta)$-spin curve $(\msP^\mr{tw}, \gamma, \mcG)$ over $S$ whose underlying pointed curve is $\msP'_S$ and which satisfies the equalities $\kappa = \kappa_{\mcG, 1} = \kappa_{\mcG, 2}^\veebar$. 
 \end{itemize}
   \epr
\begin{proof}
Since assertions (i) and (ii) are  of local nature with respect to the \'{e}tale topology of $S$,
 one may assume that 
$\mu_l$ is contained in $\Gamma (S, \mcO_S)$ for any $l$ dividing $|Z|$.

Before proceeding,
let us introduce some notation.
Given a $k$-scheme $V$, we denote by  $Z^\mr{triv}_V$ the trivial $Z$-bundle on $V$.
For each $j \in \{1,2\}$ and a positive integer $l$, we shall define $U_{j, l}$ to be the $U_j$-scheme
$\mcS pec (\mcO_S [s_j^{1/l}])$, equipped with  the $\mu_l$-action given by $s_j^{1/l} \mapsto \zeta s_j^{1/l}$ ($\zeta \in \mu_l$);
given  each $\zeta \in \mu_l$, we denote the corresponding automorphism of $U_{j, l}$ by $\alpha_{j, l, \zeta}$.
Denote by $\pi_{j, l} : U_{j, l} \migi [U_{j, l}/ \mu_{l}]$ the natural projection.
Finally, given $\xi \in Z$, we denote by $\beta_\xi : Z \isom Z$ the translation by $\xi$.

Now, we shall  prove assertion (i).
To this end,  one may assume, without loss of generality, that $S = \mr{Spec} (k)$  and $k$ is algebraically closed.
Let  $(\msP^\mr{tw}, \gamma, \mcG)$ (where $\msP^\mr{tw} := (\mfP, \{ \sigma_{\mfP, 1}, \sigma_{\mfP,2} \})$)  be
  a $(Z, \delta)$-structure
 on $\msP_k := \msP_S$.
Since 
\begin{align}
\Omega_{\mfP^\mr{log}/k} \cong \gamma^* (\Omega_{\mbP_k^\mr{log}/k}) \cong \gamma^* (\mcO_{\mbP_k}) \cong \mcO_\mfP,
\end{align} 
the classifying morphism $\mfP \migi \mcB \widehat{Z}_\delta$ of $\mcG$ factors through $\mcB Z \migi \mcB \widehat{Z}_\delta$.
That is to say, 
$\mcG$ may be regarded  as a $Z$-bundle on $\mfP$.
Denote by $l_j$ (for each $j =1,2$) the order of the stabilizer in $\mfP$ at $\sigma_{\mfP, j}$.
Then, there exists an isomorphism $\iota_{j, l_j} : [U_{j, l_j}
/\mu_{l_j}] \isom \mfP \setminus \mr{Im}(\sigma_{\mfP, j})$ of stacks over $U_j$.
Since we have assumed that $Z$ has order invertible in $k$ and $k$ is 
 algebraically closed,
$(\iota_{j, l_j} \circ \pi_{j, l_j})^*(\mcG)$ turns out to be trivial.
Let us identify $(\iota_{j, l_j} \circ \pi_{j, l_j})^*(\mcG)$ with $Z^\mr{triv}_{U_{j, l_j}}$ by a  fixed  isomorphism.
 For each $\zeta \in \mu_{l_j}$,
 the equality $(\iota_{j, l_j} \circ \pi_{j, l_j}) \circ \alpha_{j, l_j, \zeta} = (\iota_{j, l_j} \circ \pi_{j, l_j})$ 
 implies that  $\alpha_{j, l_j, \zeta}$ induces  
 an automorphism of $(\iota_{j, l_j} \circ \pi_{j, l_j})^*(\mcG)$ over $U_{j, l_j}$.
 It follows from the definition of $\kappa_{\mcG, j}$ that 
 this automorphism may be expressed  as
 $\alpha_{j, l_j, \zeta} \times \beta_{\kappa_{\mcG, l_j}(\zeta)}$  under the identification  $(\iota_{j, l_j} \circ \pi_{j, l_j})^*(\mcG) = Z^\mr{triv}_{U_{j, l_j}}$ ($= U_{j, l_j} \times_k Z$).

Next, let us consider the commutative diagram
\begin{align} \label{ee274}
\vcenter{\xymatrix{ & U \ar[ld]_{\upsilon_1} \ar[rd]^{\upsilon_2} & \\
U_{1, l_1} \ar[rd]& & 
U_{2, l_2}
 \ar[ld] \\ & \mfP, &}}
\end{align}
where $U := \mr{Spec}(k [s_1^{\pm 1/L}]) = \mr{Spec}(k [s_2^{\pm 1/L}])$  ($L := \mr{lcm} (l_1, l_2)$, $s_1^{1/L} = s_2^{-1/ L}$).
If $\tau_j : \mu_L \migi \mu_{l_j}$  denotes the surjection given by  $\zeta \mapsto \zeta^{L /l_j}$, 
then, for each $\zeta \in \mu_L$,   the following equality of  automorphisms of $Z^\mr{triv}_U$ ($= \upsilon_j^* (Z^\mr{triv}_{U_{j, l_j}})$)  holds:
\begin{align} \label{ee270}
\upsilon^*_j (\alpha_{j, l_j, \tau_j(\zeta)} \times \beta_{\kappa_{\mcG, l_j}(\zeta)}) = (\alpha_{j, L, \zeta} |_U) \times \beta_{(\kappa_{\mcG, l_j}\circ \tau_j) (\zeta)}.
\end{align}
The commutativity of (\ref{ee274}) gives rise to
an isomorphism 
\begin{align} \label{ee271}
(Z^\mr{triv}_U =) \ \upsilon_1^* (Z^\mr{triv}_{U_{1, l_1}}) \isom \upsilon_2^* (Z^\mr{triv}_{U_{2, l_2}}) \ (= Z^\mr{triv}_U)
\end{align}
 over $U$ (given by $\beta_{\xi}$ for some $\xi \in Z$).
On the other hand, the equality
 $s_1^{1/L} = s_2^{-1/L}$ implies 
  $\alpha_{1, L, \zeta} |_U = \alpha_{2, L, \zeta^{-1}} |_U$.
 By passing to (\ref{ee271}) and taking account of  (\ref{ee270}), we obtain the equalities
\begin{align}
\beta_{(\kappa_{\mcG, l_1}\circ \tau_1) (\zeta)} = \left(\beta_{\xi}^{-1} \circ \beta_{(\kappa_{\mcG, l_2}\circ \tau_2) (\zeta^{-1})}\circ \beta_{\xi}=\right)  \beta_{(\kappa_{\mcG, l_2}\circ \tau_2) (\zeta^{-1})} = \beta_{(\kappa_{\mcG, l_2}\circ \tau_2)^\veebar (\zeta)}
\end{align}
(for all $\zeta \in \mu_L$).
 Thus, 
 the equality $\kappa_{\mcG, l_1}\circ\tau_1 = (\kappa_{\mcG, 2} \circ \tau_2)^\veebar$ of morphisms $\mu_L \migi Z$ holds. 
Since both $\kappa_{\mcG, l_1}$ and $\kappa_{\mcG, l_2}$ are surjective,  we have $l_1 = l_2$ and hence,  $\kappa_{\mcG, l_1} = \kappa_{\mcG, 2}^\veebar$.
This  completes the proof of assertion (i).

Next, let us consider assertion (ii). 
Let   $\kappa :\mu_l \migiincl Z$ be an element of $\mr{Inj} (\mu, Z)$.
Denote by $\widetilde{\kappa} : S \times_k \mcB \mu_l \migi \mcB Z$ the object of $\overline{\mcI}_{\mu} (\mcB Z)_\kappa$ defined as the composite of the second projection $S \times_k \mcB \mu_l \migi \mcB \mu_l$ and $\mcB \kappa : \mcB \mu_l \migi \mcB Z$.
Since
 \begin{align}
 [U_{1, l}/\mu_l] \times_{\mbP_S} (U_1 \cap U_2) \cong U_1 \cap U_2 \cong  [U_{2, l}/\mu_l] \times_{\mbP_S} (U_1 \cap U_2),
 \end{align}
the two stacks $[U_{1, l}/\mu_l]$ and  $[U_{2, l}/\mu_l]$ may be glued together to obtain  
 a twisted curve $\mfP_{\widetilde{\kappa}}$ over $S$  equipped with   two marked points $\sigma_{\mfP, 1}$, $\sigma_{\mfP, 2}$  over $\sigma_{\mbP, 1}$, $\sigma_{\mbP, 2}$ respectively.
Denote by $\msP_{\widetilde{\kappa}}^\mr{tw} := (\mfP_{\widetilde{\kappa}}, \{ \sigma_{\mfP, 1}, \sigma_{\mfP, 2} \})$ the resulting pointed twisted curve and by  $\gamma_{\widetilde{\kappa}} : \mfP_{\widetilde{\kappa}} \migi \mbP_S$  the natural projection, i.e., $\gamma_{\widetilde{\kappa}} := \mr{coa}_{\mfP_{\widetilde{\kappa}}}$.
One may find a unique (up to isomorphism) $Z$($= \widehat{Z}_\delta$)-bundle on
$\mfP_{\widetilde{\kappa}}$
whose restriction to $\sigma_{\mfP, 1}$ and   $\sigma_{\mfP, 2}$ are  classified by $\widetilde{\kappa}$ and $\widetilde{\kappa}^\veebar$ respectively.
Indeed, 
for each $j = 1,2$, consider the trivial $Z$-bundle $Z^{\mr{triv}}_{U_{j, l}}$ on $U_{l, j}$.
Let us   equip $Z^{\mr{triv}}_{U_{1, l}}$ (resp., $Z^{\mr{triv}}_{U_{2, l}}$)
with
the $\mu_l$-action given by    $\alpha_{1, l, \zeta} \times \beta_{\kappa (\zeta)}$ (resp., $\alpha_{2, l, \zeta^{-1}} \times \beta_{\kappa^\veebar (\zeta)}$) for each $\zeta \in \mu_l$.
These  $\mu_l$-actions on $Z^{\mr{triv}}_{U_{1, l}}$ and $Z^{\mr{triv}}_{U_{2, l}}$
are compatible when restricted to $U := U_{1, l} \times_{\mbP_S} U_{2, l}$.
By means of these actions and the identity morphism   ($Z^\mr{triv}_U =$) $Z^\mr{triv}_{U_{1, l}} |_{U} \isom  Z^\mr{triv}_{U_{2, l}} |_{U}$ ($=Z^\mr{triv}_U$) of $Z^\mr{triv}_U$, 
the bundles $Z^{\mr{triv}}_{U_{1, l}}$, $Z^{\mr{triv}}_{U_{2, l}}$
may be glued together to obtain  the desired  $Z$-bundle $\mcG_{\widetilde{\kappa}}$.
The injectivity of $\kappa$ implies that the classifying morphism $\mfP_{\widetilde{\kappa}} \migi \mcB Z$ of $\mcG_{\widetilde{\kappa}}$ is representable.
Thus, we obtain a $(Z, \delta)$-structure $(\msP^\mr{tw}_{\widetilde{\kappa}}, \gamma_{\widetilde{\kappa}}, \mcG_{\widetilde{\kappa}})$ on $\msP_S$ satisfying the required conditions.
This implies  the validity of the existence portion.
The uniqueness portion follows immediately from
the above construction of $(\msP^\mr{tw}_{\widetilde{\kappa}}, \gamma_{\widetilde{\kappa}}, \mcG_{\widetilde{\kappa}})$.
This completes the proof of assertion (ii). 
\end{proof}
\SSP

\bco \label{C08} 
We have
\begin{align}
 \mfS \mfp_{Z, \delta, 0,3, (\kappa_1, \kappa_2, \overline{\delta})} \cong \begin{cases}
 \mcB Z & \text{if $\kappa_1 = \kappa_2^\veebar$;}
 \\ 
 \emptyset & \text{if otherwise}.
\end{cases}
\end{align}
\eco
\begin{proof}
The assertion follows from   Propositions \ref{P020}, (i), \ref{L010}, (i),  and  \ref{P08ff}.
\end{proof}

\vspace{10mm}
\section{Twisted opers on  pointed stable curves}  \label{y032}\SSP

In this section, we define 
the notion of a {\it (faithful) twisted $G$-oper} on a pointed stable curve.
We construct the moduli space, denoted by $\mfO \mfp_{G,g,r, (\vec{\pmb{\rho}})}$,  classifying 
pointed stable curves together  with a faithful twisted  $G$-oper (of prescribed radii $\vec{\pmb{\rho}}$).
The main result of this section asserts (cf. Theorems \ref{P0012} and \ref{P04}) that
this moduli space may be represented by a smooth Deligne-Mumford stack which is flat  over 
$\overline{\mfM}_{g,r}$ of constant relative dimension.
Also, we 
study 
 the  faithful twisted $G$-opers on the projective line with two or three marked points (cf. Proposition \ref{Pff01} and Corollary \ref{cff01}).

\LSP
\subsection{Algebraic groups and Lie algebras.} \label{y033}
First, we introduce some notation concerning algebraic groups and Lie algebras (cf. ~\cite[\S\S\,2.1-2.2]{Wak5}).
Let  $(G, T)$ be  a split semisimple algebraic group over $k$, where $T$ denotes a maximal torus $T$ of $G$.
In this section, {\it we shall assume  that $\mr{char} (k) =0$ or $\mr{char} (k) = p$ for some prime $p$ satisfying the  condition $(*)_{G}$ described in Introduction}.
Let us  fix  a Borel subgroup $B$ defined over $k$ containing $T$.
In particular, we obtain a natural surjection $B \migisurj (B/[B, B] \cong) \ T$.
Denote by $\mfg$, $\mft$,  and  $\mfb$ the Lie algebras of $G$,  $T$,  and $B$ respectively (hence, $\mft \subseteq \mfb \subseteq \mfg$).

For each character  $\beta$ of $T$, we write   
\begin{align}
\mfg^\beta := \big\{ x \in \mfg \ \big| \ \text{$\mr{ad}(t)(x) =  \beta (t) \cdot x$ for all  $t\in T$}  \big\}.
\end{align}
Let $\Gamma$ denote the set of simple roots in $B$ with respect to $T$.
For each $\alpha \in \Gamma$, we fix a generator $x_\alpha$ of $\mfg^\alpha$. 
 Write $p_{1} := \sum_{\alpha \in \Gamma} x_\alpha$ ($\in \mfg$) and $\check{\rho} := \sum_{\alpha \in \Gamma} \check{\omega}_\alpha$ ($\in \mft$),
where $\check{\omega}_\alpha$ (for each $\alpha \in \Gamma$) denotes the fundamental coweight of $\alpha$, regarded   as an element of $\mft$ via differentiation.
There exists  a unique  collection $(y_\alpha)_{\alpha \in \Gamma}$, where $y_\alpha$ is a generator of $\mfg^{-\alpha}$, such that if we write $p_{-1}  := \sum_{\alpha \in \Gamma} y_\alpha$,
then the set 
 $\{  p_{-1}, 2 \check{\rho}, p_1  \}$ forms an $\mfs \mfl_2$-triple.

Finally, recall a canonical 
 decreasing filtration $\{ \mfg^j \}_{j \in \mbZ}$
   on $\mfg$
 such that $\mfg^0 = \mfb$, $\mfg^0/ \mfg^1= \bigoplus_{\alpha \in \Gamma} \mfg^\alpha$, and $[\mfg^{j_1}, \mfg^{j_2}] \subseteq \mfg^{j_1 +j_2}$ for $j_1$, $j_2 \in \mbZ$.

\LSP
\subsection{Twisted $G$-opers on stacky log curves} \label{y034}
Let 
$S^\mr{log}$ be an fs log scheme (or, more generally, an fs log stack) over $k$,  $\mfU^\mr{log}$ a stacky  log curve over $S^\mr{log}$,  
and
$\pi : \mcE \migi \mfU$ a right $G$-bundle on $\mfU$.
By pulling-back the log structure of $\mfU^\mr{log}$  via $\pi$, one may obtain a log structure on $\mcE$; we denote  the resulting log stack by   $\mcE^\mr{log}$.
The $G$-action on $\mcE$ carries a  $G$-action on  the direct image   $\pi_{*}(\mcT_{\mcE^\mr{log}/S^\mr{log}})$ of  $\mcT_{\mcE^\mr{log}/S^\mr{log}}$.
 Denote by  $\widetilde{\mcT}_{\mcE^\mr{log}/S^\mr{log}}$
  the subsheaf of $G$-invariant sections  of $\pi_{*}(\mcT_{\mcE^\mr{log}/S^\mr{log}})$.
The differential of  $\pi$  gives rises to  a short exact sequence
\begin{align}
 0 \longmigi  \mfg_{\mcE} \longmigi  \widetilde{\mcT}_{\mcE^\mr{log}/S^\mr{log}} \stackrel{d^\mr{log}_{\mcE}}{\longmigi}  \mcT_{\mfU^\mr{log}/S^\mr{log}} \longmigi 0\label{Ex0}
 \end{align}
of $\mcO_\mfU$-modules.

\SSP
\bde \label{DdD011} An {\bf $S^\mr{log}$-connection} on $\mcE$ is  an $\mcO_\mfU$-linear morphism $\nabla  : \mcT_{\mfU^\mr{log}/S^\mr{log}} \migi \widetilde{\mcT}_{\mcE^\mr{log}/S^\mr{log}}$ with $d^\mr{log}_{\mcE} \circ \nabla = \mr{id}_{\mcT_{\mfU^\mr{log}/S^\mr{log}}}$.
  \ede
\SSP

Now, suppose that we are given  a right $B$-bundle
$\pi_B : \mcE_B  \migi \mfU$   on $\mfU$. 
Denote by
$\pi_G : \left(\mcE_B \times^B G =:\right) \mcE_G \migi \mfU$ the $G$-bundle on $\mfU$ obtained by  change of structure group via the inclusion $B \migiincl G$.
The natural morphism $\mcE_B \migi \mcE_G$ yields 
a canonical isomorphism $\mfg_{\mcE_B} \isom \mfg_{\mcE_G}$ and moreover
a morphism between short exact sequences:
\begin{equation} \begin{CD}
0 @>>> \mfb_{ \mcE_B} @>>> \widetilde{\mcT}_{\mcE_B^\mr{log}/S^\mr{log}} @> d_{\mcE_B}^\mr{log} >> \mcT_{\mfU^\mr{log}/S^\mr{log}} @>>> 0
\\
@. @VV \iota_{\mfg/\mfb} V  @VV \widetilde{\iota}_{\mfg/\mfb} V @VV \mr{id} V @.
\\
0 @>>>  \mfg_{\mcE_G} @>>> \widetilde{\mcT}_{\mcE_G^\mr{log}/S^\mr{log}} @>> d_{\mcE_G}^\mr{log} > \mcT_{\mfU^\mr{log}/S^\mr{log}} @>>> 0,
\end{CD}  \label{extension} \end{equation}
 where the upper and lower horizontal sequences are (\ref{Ex0}) applied to $\mcE_B$ and $\mcE_G$ respectively.
Since $\mfg^j  \left(\subseteq \mfg\right)$  is closed under the adjoint action of $B$,
  one obtains   vector bundles $\mfg^j_{\mcE_B}$ ($j \in \mbZ$) associated with $\mcE_B  \times^B \mfg^j$.
  The collection $\{\mfg^j_{\mcE_B}\}_{j \leq 0}$ defines 
     a decreasing filtration on $\mfg_{\mcE_B}  \left(\cong \mfg_{\mcE_G }\right)$.
On the other hand,  diagram (\ref{extension})  induces a composite isomorphism
\begin{equation}  \mfg_{\mcE_B} / \mfg^0_{\mcE_B} \isom \mfg_{\mcE_G}/\iota_{\mfg/\mfb} (\mfb_{\mcE_B}) \isom \widetilde{\mcT}_{\mcE_G^\mr{log}/S^\mr{log}} /\widetilde{\iota}_{\mfg/\mfb}( \widetilde{\mcT}_{\mcE_B^\mr{log}/S^\mr{log}}).\end{equation}
The filtration $\{ \mfg^j_{\mcE_B} \}_{j \leq 0}$
  carries, via this composite isomorphism,  a decreasing filtration
$\{ \widetilde{\mcT}_{\mcE_G^\mr{log}/S^\mr{log}}^j \}_{j \leq 0}$
  on $\widetilde{\mcT}_{\mcE_G^\mr{log}/S^\mr{log}}$ in such a way that  
$\widetilde{\mcT}_{\mcE_G^\mr{log}/S^\mr{log}}^0 = \widetilde{\iota}_{\mfg/\mfb}(\widetilde{\mcT}_{\mcE_B^\mr{log}/S^\mr{log}})$ and  the resulting morphism 
$\mfg^{j-1}_{\mcE_B}/ \mfg^{j}_{\mcE_B} \migi \widetilde{\mcT}_{\mcE_G^\mr{log}/S^\mr{log}}^{j-1}/ \widetilde{\mcT}_{\mcE_G^\mr{log}/S^\mr{log}}^{j}$ is an isomorphism.
Since each $\mfg^{-\alpha}$ ($\alpha \in \Gamma$) is closed under the $B$-action  defined by the composite $B \migisurj T \xrightarrow{\mr{adj.\,rep.}} \mr{Aut} (\mfg^{-\alpha})$, the canonical  decomposition 
 $\mfg^{-1}/\mfg^0 = \bigoplus_{\alpha \in \Gamma} \mfg^{-\alpha}$
 gives  rise to  a  decomposition 
\begin{equation} \label{decom97} \widetilde{\mcT}_{\mcE_G^\mr{log}/S^\mr{log}}^{-1} / \widetilde{\mcT}^0_{\mcE_G^\mr{log}/S^\mr{log}} = \bigoplus_{\alpha \in \Gamma}    \mfg^{-\alpha}_{\mcE_B} .\end{equation}

\SSP
\bde \label{Dy0351} \begin{itemize}
\item[(i)]
  A {\bf twisted $G$-oper}
on  $\mfU^\mr{log}/S^\mr{log}$ is a pair
 \begin{equation}
  \msE^\spadesuit := (\pi_B: \mcE_B \migi \mfU, \nabla : \mcT_{\mfU^\mr{log}/S^\mr{log}} \migi \widetilde{\mcT}_{\mcE_G^\mr{log}/S^\mr{log}})\end{equation}
  consisting of a $B$-bundle $\mcE_B$ on $\mfU$ and an $S^\mr{log}$-connection $\nabla$  
  on the  $G$-bundle $\pi_G : \mcE_G \migi \mfU$  induced by $\mcE_B$ satisfying the following two conditions:
\begin{itemize}
\item[$\bullet$]
 $\nabla (\mcT_{\mfU^\mr{log}/S^\mr{log}}) \subseteq  \widetilde{\mcT}_{\mcE_G^\mr{log}/S^\mr{log}}^{-1}$,
\item[$\bullet$]
For any $\alpha \in \Gamma$, the composite
\begin{equation}
\label{isomoper}
 \mcT_{\mfU^\mr{log}/S^\mr{log}} \xrightarrow{\nabla} \widetilde{\mcT}_{\mcE_G^\mr{log}/S^\mr{log}}^{-1} \migisurj \widetilde{\mcT}_{\mcE^\mr{log}_G/S^\mr{log}}^{-1} /\widetilde{\mcT}^0_{\mcE_G^\mr{log}/S^\mr{log}} \migisurj   \mfg^{-\alpha}_{\mcE_B}
  \end{equation}
is an isomorphism, where the third arrow denotes the natural projection with respect to   decomposition (\ref{decom97}).
\end{itemize}
If $\mfU^\mr{log}/S^\mr{log} = \mfX^\mr{log}/S^\mr{log}$
for a
 pointed twisted curve  $\msX^\mr{tw} := (\mfX/S, \{ \sigma_{\mfX, i} \}_{i})$,
 then we shall refer to any twisted  $G$-oper on $\mfX^\mr{log}/S^\mr{log}$ as a {\bf twisted $G$-oper on $\msX^\mr{tw}$}.

If  $\msX$ is  a pointed nodal curve, then
we define a {\bf twisted $G$-oper} on $\msX$ to be a collection of data
\begin{align} \label{ee100}
\mbE^\spadesuit := (\msX^\mr{tw}, \gamma, \msE^\spadesuit),
\end{align}
consisting of a twistification  $(\msX^\mr{tw}, \gamma)$ of $\msX$ and a twisted $G$-oper $\msE^\spadesuit$ on $\msX^\mr{tw}$.
\item[(ii)]
For each $j  \in \{1,2\}$, let $S_j$ be a $k$-scheme and $\mbX_j^{\spadesuit} := (\msX_j, \msX_j^\mr{tw}, \gamma_j, \msE_j^\spadesuit)$
  a collection   consisting of a pointed nodal curve $\msX_j$ and a twisted $G$-oper $(\msX_j^\mr{tw}, \gamma_j, \msE_j^\spadesuit)$ on it, where  $\msX^\mr{tw}_j := (f_j : \mfX_j \migi S_j, \{ \sigma_{\mfX_j, i} \}_{i})$ and  $\msE_j^\spadesuit := (\pi_{B, j} : \mcE_{B, j} \migi \mfX_j, \nabla_{j})$.
A {\bf $1$-morphism} (or just a {\bf morphism})  from 
$\mbX_1^{\spadesuit}$ to $\mbX_2^{\spadesuit}$
is a triple
\begin{align}
\alpha^{\spadesuit} := (\alpha_S^\spadesuit, \alpha_\mfX^\spadesuit, \alpha_\mcE^\spadesuit)
\end{align}
consisting of morphisms  of $k$-stacks
which make the following diagram  $1$-commutative:
\begin{align}
\begin{CD}
\mcE_{B, 1} @> \pi_{B, 1} >> \mfX_1 @> f_1 >> S_1  
\\
@VV \alpha_\mcE^\spadesuit V @VV \alpha_\mfX^\spadesuit V @VV \alpha_S^\spadesuit V
\\
\mcE_{B, 2} @>> \pi_{B, 2} > \mfX_2 @>> f_2 > S_2,
\end{CD}
\end{align}
where
\begin{itemize}
\item[$\bullet$]
the right-hand square diagram forms a morphism of  pointed twisted curves (cf. Definition \ref{D012}, (ii));
\item[$\bullet$]
the left-hand square is cartesian, and $\alpha_\mcE^\spadesuit$ is compatible with the respective $B$-actions of $\mcE_{B,1}$ and $\mcE_{B,2}$;
\item[$\bullet$]
the  morphism $\mcE_{G,1} \ (:=\mcE_{B, 1} \times^B G) \migi \mcE_{G,2} \ (:= \mcE_{B, 2} \times^B G)$ induced by $\alpha_\mcE^\spadesuit$ is compatible with the respective connections $\nabla_{1}$,  $\nabla_{2}$. 
\end{itemize}
In particular, 
one may associate, to  such a morphism $\alpha^{\spadesuit}$, a morphism $\alpha^\spadesuit_\msX : \msX_1 \migi \msX_2$ between the underlying pointed nodal curves.
\item[(iii)]
Let $\mbX_j^{\spadesuit}$ ($j =1,2$) be as in (ii) and $\alpha_l^{\spadesuit} := (\alpha_{S, l}^\spadesuit, \alpha_{\mfX, l}^\spadesuit, \alpha_{\mcE, l}^{\spadesuit})$ ($l=1,2$) morphisms from  $\mbX_1^{\spadesuit}$ to $\mbX_2^{\spadesuit}$.
A {\bf $2$-morphism} from $\alpha_1^{\spadesuit}$ to $\alpha_2^{\spadesuit}$
is a triple of natural transformations
\begin{align}
\mfa^{\spadesuit} := (\alpha_{S,1}^\spadesuit \stackrel{\mfa_S^\spadesuit}{\Rightarrow} \alpha_{S, 2}^\spadesuit, \alpha_{\mfX, 1} \stackrel{\mfa_\mfX^\spadesuit}{\Rightarrow} \alpha_{\mfX, 2}^\spadesuit, \alpha_{\mcE, 1}^\spadesuit \stackrel{\mfa_\mcE^\spadesuit}{\Rightarrow} \alpha_{\mcE, 2}^\spadesuit)
\end{align}
 compatible with each other (hence, $\mfa^\spadesuit_S$ coincides with the identity natural transformation).
\end{itemize}
  \ede

\LSP
\subsection{Faithful twisted $G$-opers and their moduli} \label{y03f4}
Next, we shall introduce the notion of a {\it faithful} twisted $G$-oper.
As discussed below, the faithfulness should be imposed on twisted $G$-opers to  construct a suitable moduli space classifying them.
Let $\mfU^\mr{log}/S^\mr{log}$ be as before.

\SSP
\bde \label{Dy035}
 We shall say that  
 a twisted $G$-oper $\msE^\spadesuit := (\mcE_B, \nabla)$  on $\mfU^\mr{log}/S^\mr{log}$
  is {\bf faithful} if 
 the classifying morphism
 $\mfU \migi \mcB T$ of the $T$-bundle  $\mcE_B \times^B T$ is representable.
Also, let  $\mbE^\spadesuit := (\msX^\mr{tw}, \gamma, \msE^\spadesuit)$ be a twisted $G$-oper on a pointed nodal curve.
Then, we shall say that $\mbE^\spadesuit$ is {\bf faithful} if $\msE^\spadesuit$ is faithful.
  \ede
\SSP

We shall  examine the  relationship with the notion  of an extended spin structure discussed in \S\,\ref{S01}.
Denote by $Z$ the center of $G$.
Because of  the assumption 
imposed at the beginning of \S\,\ref{y033}, $G$ is finite and has order invertible in $k$.
Write $G_\mr{ad} := G/Z$ (i.e., the {\it adjoint group} of $G$) and $T_\mr{ad} := T/Z$.
Let   $\lambda : \mbG_m \migi T_{\mr{ad}}$ (cf. ~\cite[\S\,3.4.1]{BD1})  be the morphism
 determined by the condition   that
for any 
$\alpha \in \Gamma$,
$\lambda (t)$
 acts on $\mfg^{\alpha}$ (via the adjoint representation $(T_\mr{ad} \subseteq) \ G_\mr{ad} \migi \mr{GL} (\mfg)$) as multiplication by
$t$.
Then, one may find a unique morphism  
$\lambda^\sharp : \mbG_m \migi T$ (cf. ~\cite[Eq.\,(54)]{BD1}) such that $\lambda (t)^2 = \lambda^\sharp (t)$ mod $Z$ for any $t \in \mbG_m$.
The morphism $\lambda^\sharp$ restricts to a  morphism $\delta^\sharp : \mu_2 \migi Z$ ($\subseteq T$).
The following square diagram  is verified to be commutative:
\begin{align} \label{ee881}
\begin{CD}
\mu_2 @> \text{incl.} >> \mbG_m
\\
@V \delta^\sharp VV @VV v \mapsto (\lambda^\sharp (v), v^2) V
\\
Z @>> z \mapsto (z, e) >  T \times_{T_\mr{ad}, \lambda} \mbG_m,
\end{CD}
\end{align}
 where $e$ denotes the unit of $\mbG_m$. 
Hence, this diagram  determines  a morphism 
\begin{align} \label{ee286}
\widehat{Z}_{\delta^\sharp} \ (\cong Z \times^{\delta^\sharp, \mu_2}\mbG_m) \migi T \times_{T_\mr{ad}, \lambda} \mbG_m.
\end{align}
 This morphism  is an isomorphism since it fits into the following morphism of short exact sequences:
\begin{align}
\begin{CD}
0 @>>> Z @>>> \widehat{Z}_{\delta^\sharp} @> \nu >> \mbG_m @>>>0
\\
@. @VV \mr{id}_Z V @VV (\ref{ee286}) V  @VV \mr{id}_{\mbG_m}V @.
\\
0 @>>> Z @>>z \mapsto (z, e)> T \times_{T_\mr{ad}, \lambda} \mbG_m @>> (h, g) \mapsto g > \mbG_m @>>> 0, 
\end{CD}
\end{align}
where the upper horizontal sequence is the lower horizontal sequence in (\ref{ee028}).

Now, let $\msE^\spadesuit := (\mcE_B, \nabla)$ be a twisted $G$-oper on $\mfU^\mr{log}/S^\mr{log}$.
If $(\Omega_{\mfU^\mr{log}/S^\mr{log}})^\times$ denotes the $\mbG_m$-bundle on $\mfU$ corresponding to the line bundle $\Omega_{\mfU^\mr{log}/S^\mr{log}}$, then
 the $T_\mr{ad}$-bundle  $\mcE_B \times^B T_{\mr{ad}}$ induced  from $\mcE_B$ via  change of structure group by the composite $B \migisurj T \migisurj T_\mr{ad}$ is isomorphic to $(\Omega_{\mfU^\mr{log}/S^\mr{log}})^\times \times^{\mbG_m, \lambda} T_\mr{ad}$ (cf. the discussion in ~\cite[\S\,3.4.1]{BD1}).
 Hence, by passing to  (\ref{ee286}), one may use the $T$-bundle $\mcE_B \times^B T$ and the $\mbG_m$-bundle $(\Omega_{\mfU^\mr{log}/S^\mr{log}})^\times$  to obtain   
 a $\widehat{Z}_{\delta^\sharp}$-bundle 
 \begin{align} \label{ee720}
 \pi_{\mcG, \msE^\spadesuit} := \mcG_{\msE^\spadesuit} \migi \mfU.
 \end{align}
 By definition, $\mcG_{\msE^\spadesuit}  \times^{\widehat{Z}_{\delta^\sharp}, \nu} \mbG_m$ is isomorphic to  $(\Omega_{\mfU^\mr{log}/S^\mr{log}})^\times$.

\SSP
\bpr \label{PPPP08}
The twisted $G$-oper  $\msE^\spadesuit$ is faithful if and only if  $\mcG_{\msE^\spadesuit}$ forms a $(Z, \delta^\sharp)$-structure on $\mfU^\mr{log}/S^\mr{log}$.
 \epr
\begin{proof}
Let us consider the composite
\begin{align} \label{ee601}
\widehat{Z}_{\delta^\sharp} \xrightarrow{(\ref{ee286})} T \times_{T_\mr{ad}, \lambda} \mbG_m \longmigi T,
\end{align}
where the second arrow denotes the first projection.
The composite of the classifying morphism $[\mcG_{\msE^\spadesuit}] : \mfU \migi \widehat{Z}_{\delta^\sharp}$ of $\mcG_{\msE^\spadesuit}$ and the morphism  $\mcB \widehat{Z}_{\delta^\sharp} \migi \mcB T$ induced by (\ref{ee601})   coincides with the classifying morphism $[\mcE_B \times^B T] : \mfU \migi \mcB T$ of $\mcE_B \times^B T$.
 Since (\ref{ee601}) is a closed immersion, it follows from  ~\cite[Lemma 4.4.3]{AV}  that $[\mcG_{\msE^\spadesuit}]$ is representable  if and only if $[\mcE_B \times^B T]$ is representable.
 This completes the proof of Proposition \ref{PPPP08}.
\end{proof}
\SSP

\begin{rema} \label{ppp037} 
Let us consider the case where $G$ is of adjoint type, i.e., $G = G_\mr{ad}$.
The definition of a faithful twisted $G$-oper
  introduced  above may be  identified with  the classical definition of a $\mfg$-oper  in the sense of ~\cite[Definition 2.1, (i)]{Wak5}.
Indeed, let 
$\msX$ be a pointed nodal curve and $\mbE^\spadesuit := (\msX^\mr{tw}, \gamma, \msE^\spadesuit)$ 
 a twisted $G$-oper on $\msX$.
The representability of the classifying morphism  $\mfX \migi \widehat{Z}_{\delta^\sharp}$ of $\mcG_{\msE^\spadesuit}$ together with  the assumption $Z = \{ 1 \}$ imply (cf. ~\cite[Lemma 4.4.3]{AV}) that 
the stabilizers of the nodes and the marked points  of $\msX^\mr{tw}$ are trivial.
Hence, 
the morphism $\gamma : \mfX \migi X$
 must be  an isomorphism and $\msE^\spadesuit$ specifies  a $\mfg$-oper  on $\msX$ in the classical sense.
 Conversely, given a $\mfg$-oper $\msE^\spadesuit$ on $\msX$, we obtain a faithful twisted $G$-oper $(\msX, \mr{id}_X, \msE^\spadesuit)$ on $\msX$.
 In this way, we shall not distinguish between faithful twisted $G$-opers (in the case of $G = G_\mr{ad}$) and $\mfg$-opers.
 \end{rema}
\SSP

By Definition \ref{Dy0351}, (i)-(iii),
the collections $(\msX, \mbE^\spadesuit)$ of  a pointed nodal  curve $\msX$  and  a twisted $G$-oper $\mbE^\spadesuit$ on $\msX$ 
 form a $2$-category.
Just as in the case of $\mfS \mfp_{Z, \delta, g,r}$ (cf. the discussion following Definition \ref{y03f5D}), 
this $2$-category
specifies a category fibered in groupoids over $\mfS \mfc \mfh_{/k}$.

Given   a pair of nonnegative integers $(g, r)$  with $2g-2+r >0$,
   we shall denote by  
\begin{align}
\mfO \mfp_{G, g,r}
\end{align}
the category of pairs $(\msX, \mbE^\spadesuit)$ consisting of an $r$-pointed stable curve   $\msX$ of  genus $g$ over a $k$-scheme and a faithful twisted $G$-oper $\mbE^\spadesuit$ on $\msX$.
The assignments $(\msX,  \mbE^\spadesuit) \mapsto \msX$ and $(\msX,  \mbE^\spadesuit) \mapsto (\msX, \mcG_{\msE^\spadesuit})$ (where $\mbE^\spadesuit := (\msX^\mr{tw}, \gamma, \msE^\spadesuit)$) determine  functors
\begin{align} \label{ee077}
\mfO \mfp_{G, g,r} \migi \overline{\mfM}_{g,r} \hspace{5mm} \text{and} \hspace{5mm}\mfO \mfp_{G, g,r} \migi \mfS \mfp_{Z, \delta^\sharp, g,r}
\end{align}
respectively.
Also, 
by  change of structure group via $G \migisurj G_\mr{ad}$, we obtain a functor
\begin{align} \label{e051}
\mr{op}_{\mr{ad}} : \mfO \mfp_{G,g,r} \migi \mfO \mfp_{G_{\mr{ad}},g,r}
\end{align}
over $\overline{\mfM}_{g,r}$.

\SSP
\begin{rema} \label{gRRR037} Let 
$\msX$ be a pointed stable curve.
 We shall refer
 to each  $2$-isomorphism class of a $1$-isomorphism defined in Definition \ref{Dy0351}, (ii)
  inducing the identity morphism of $\msX$
   as an {\bf isomorphism of faithful twisted $G$-opers} on $\msX$.
 In this way, we obtain the groupoid of faithful twisted $G$-opers  on $\msX$.
\end{rema}

\SSP
\bt \label{P0012}
\begin{itemize}
\item[(i)]
The functor 
 \begin{align} \label{e052}
 \mfO \mfp_{G,g,r}  \isom \mfO \mfp_{G_{\mr{ad}},g,r} \times_{\overline{\mfM}_{g,r}} \mfS \mfp_{Z, \delta^\sharp, g,r}
 \end{align}
induced by (\ref{e051}) and  the second morphism in (\ref{ee077}) is an equivalence of categories over $\overline{\mfM}_{g,r}$.
  \item[(ii)]
 $\mfO \mfp_{G, g,r}$  may be represented by a nonempty  smooth  Deligne-Mumford stack over $k$ which is flat  over $\overline{\mfM}_{g,r}$ of relative dimension 
 \begin{align}
 (g-1) \cdot  \mr{dim} (G) + \frac{r}{2} \cdot (\mr{dim} (G) + \mr{rk} (G)),
 \end{align}
  where $\mr{dim}(G)$ and  $\mr{rk} (G)$ denote the dimension  and rank of $G$ respectively.
  Also, the morphism $\mr{op}_{\mr{ad}}$ is finite, flat, and generically \'{e}tale.
 \end{itemize}
  \et
\begin{proof}
We shall 
consider assertion  (i).
 Let  
 $\msX$ be an $r$-pointed stable curve of genus $g$ over a $k$-scheme $S$.
 Suppose that we are given a $G_\mr{ad}$-oper $\msE_\mr{ad}^\spadesuit := (\mcE_{B_\mr{ad}}, \nabla)$ on $\msX$ and a $(Z, \delta^\sharp)$-structure $(\msX^\mr{tw}, \gamma, \mcG)$ on $\msX$.
 In particular, $\mbX := (\msX, \msX^\mr{tw}, \gamma, \mcG)$ forms a pointed stable $(Z, \delta^\sharp)$-spin curve.
Both the $B_\mr{ad}$-bundle  $\gamma^*(\mcE_{B_\mr{ad}})$ on $\mfX$ and the $T$-bundle $\mcG \times^{\widehat{Z}_{\delta^\sharp}} T$ obtained from $\mcG$ via change of structure group by the composite of (\ref{ee286}) and the first projection $T \times_{T_\mr{ad}} \mbG_m \migi T$  induce the same $T_\mr{ad}$-bundle $(\Omega_{\mfX^\mr{log}/S^\mr{log}})^{\times} \times^{\mbG_m, \lambda} T_\mr{ad} \left(\cong (\gamma^*(\Omega_{X^\mr{log}/S^\mr{log}}))^{\times} \times^{\mbG_m, \lambda} T_\mr{ad}\right)$.
Since $B \cong B_\mr{ad} \times_{T_\mr{ad}} T$, these bundles give rise to   a $B$-bundle $\mcE_B$ on $\mfX$.
The natural morphism $\mcE_G  \left(:= \mcE_B \times^B G\right)\migi \mcE_{G_\mr{ad}}$ induces an isomorphism
$\widetilde{\mcT}_{\mcE_G^\mr{log}/S^\mr{log}} \isom \widetilde{\mcT}_{\mcE_{G_\mr{ad}}^\mr{log}/S^\mr{log}}$.
 By this isomorphism,  $\nabla$ may be regarded  as an $S^\mr{log}$-connection on $\mcE_G$.
Hence, the collection $\mbE^\spadesuit := (\msX^\mr{tw}, \gamma,  (\mcE_B, \nabla))$
  forms a faithful twisted $G$-oper on $\msX$.
One may verify that the resulting assignment 
$((\msX, \msE_\mr{ad}^\spadesuit), \mbX) \mapsto (\msX, \mbE^\spadesuit)$
 determines a functor
\begin{align}
\mfO \mfp_{G_{\mr{ad}},g,r} \times_{\overline{\mfM}_{g,r}} \mfS \mfp_{Z, \delta^\sharp, g,r} \migi  \mfO \mfp_{G,g,r},
\end{align}
and that  it specifies  the  inverse of (\ref{e052}).
This  complete the proof of  assertion  (i).

Assertion (ii) follows from Theorem  \ref{P08f} and ~\cite[Theorem A]{Wak5}.
\end{proof}

\LSP
\subsection{Radii of twisted $G$-opers.} \label{y034}
Let us write  $\mfc :=  \mfg \ooalign{$/$ \cr $\,/$}\hspace{-0.5mm}G$, i.e., the GIT quotient of $\mfg$ by the adjoint action of $G$.
Also,  write $\chi : \mfg \migi \mfc$ for the natural projection.
The $k$-scheme $\mfc$ has 
 the involution
$\rho \mapsto \rho^\veebar$   arising from the  automorphism of $\mfg$ given by  multiplication by $(-1)$.
This involution induces  an involution $(-)^\veebar$ on $\mfc (S) \times \mr{Inj} (\mu, Z)$ defined  by
$(\rho, \kappa)^\veebar := (\rho^\veebar, \kappa^\veebar)$.
Let us  write
\begin{align} \label{ee180}
\varepsilon := \chi (- \check{\rho}) \in \mfc (k), \hspace{10mm} \pmb{\varepsilon} := (\varepsilon, \overline{\delta}^\sharp) \in \mfc (k) \times \mr{Inj} (\mu, Z).
\end{align}
Then, the equalities   $\varepsilon^\veebar = \varepsilon$ and  $\pmb{\varepsilon}^\veebar = \pmb{\varepsilon}$ hold.

  The  quotient stack $[\mfc / Z]$ induced by  the trivial $Z$-action on $\mfc$   is canonically  isomorphic to $\mfc \times_k \mcB Z$.
We identify the stack  of cyclotomic gerbes $\overline{\mcI}_\mu ([\mfc /Z])$ in $[\mfc/Z]$ with $\mfc \times_k \overline{\mcI}_\mu (\mcB Z)$ by the  composite of natural  isomorphisms
\begin{align} \label{e0020}
\overline{\mcI}_\mu ([\mfc/Z])   \isom \overline{\mcI}_\mu (\mfc \times_k \mcB Z) 
 \isom \mfc \times_k \overline{\mcI}_\mu (\mcB Z).
\end{align}
Moreover, according to  decomposition (\ref{ee034}),
$\overline{\mcI}_\mu ([\mfc /Z])$ may be identified with $\coprod_{\kappa \in \mr{Inj}(\mu, Z)} \mfc \times_k \overline{\mcI}_\mu (\mcB Z)_\kappa$.
Given $\pmb{\rho} := (\rho, \kappa) \in \mfc (k) \times \mr{Inj} (\mu, Z)$,
we obtain a closed substack
\begin{align}
\overline{\mcI}_\mu ([\mfc / Z])_{\pmb{\rho}} 
\end{align}
of $\overline{\mcI}_{\mu} ([\mfc / Z])$ corresponding, via (\ref{e0020}), to the close immersion  
$\overline{\mcI}_\mu (\mcB Z)_\kappa \ (= \mr{Spec}(k)$ $\times_k \overline{\mcI}_\mu (\mcB Z)_\kappa) \migiincl \mfc \times_k \overline{\mcI}_\mu (\mcB Z)$  defined as the product of  $\rho$ and the inclusion  $\overline{\mcI}_\mu (\mcB Z)_\kappa$ $\migiincl \overline{\mcI}_\mu (\mcB Z)$.

Let 
$r$ be a positive integer, 
 $\msX := (X, \{ \sigma_{X, i} \}_{i=1}^r)$  an $r$-pointed nodal curve  over a $k$-scheme $S$,  and 
$\mbE^\spadesuit := (\msX^\mr{tw}, \gamma, \msE^\spadesuit)$ a faithful twisted $G$-oper  on $\msX$.
Let  
$(\msE^\spadesuit_\mr{ad}, \mcG)$ be  the pair of a faithful twisted $G_\mr{ad}$-oper  on $\msX$ and a $(Z, \delta^\sharp)$-structure on $\msX^\mr{tw}$ corresponding to  
$\mbE^\spadesuit$ via  equivalence of categories (\ref{e052}).

Recall from ~\cite[Definition 2.32]{Wak5} that  the radius  of $\msE^\spadesuit_\mr{ad}= (\mcE_{B, \mr{ad}}, \nabla)$ at each $\sigma_{X, i}$ ($i =1, \cdots, r$) is defined as a certain element of $\mfc (S)$, which we shall denote by $\rho_{\msE^\spadesuit_{\mr{ad}}, i}$.
That is to say,
if $\mcE_{G, \mr{ad}}$ denotes the $G_\mr{ad}$-bundle $\mcE_{B, \mr{ad}} \times^{B_\mr{ad}} G_\mr{ad}$ and $\mu_i^{ \nabla}$ denotes the {\it monodromy operator} of $\nabla$ at $\sigma_{X, i}$ in the sense of ~\cite[Definition 1.46]{Wak5},
then
$\rho_{\msE^\spadesuit_{\mr{ad}}, i}$ is defined as the image  of 
the pair $(\sigma^*_{X, i} (\mcE_{G, \mr{ad}}), \mu_i^{\nabla})$ via the natural projection $[\mfg/G] \migi \mfc$.
Here, the quotient stack $[\mfg / G]$  represents the functor which, to any $k$-scheme $T$,  assigns the groupoid of pairs $(\mcF, \mu)$ consisting of a $G$-bundle on $T$ and $\mu \in \Gamma (T, \mfg_T)$. 
Let us write
\begin{align} \label{ee600}
\pmb{\rho}_{\mbE^\spadesuit, i} &:= (\rho_{\mcE_\mr{ad}^\spadesuit, i}, \kappa_{\mcG, i}) \in \mfc (S)  \times \mr{Inj} (\mu, Z),\\
 \widetilde{\pmb{\rho}}_{\mbE^\spadesuit, i} &:= (\rho_{\mcE_\mr{ad}^\spadesuit, i}, \widetilde{\kappa}_{\mcG, i}) \in   \mr{Ob} 
(\overline{\mcI}_\mu ([\mfc /Z]) (S))  \left(\stackrel{(\ref{e0020})}{\cong} \mfc (S) \times \mr{Ob} (\overline{\mcI}_\mu (\mcB Z) (S)) \right). \notag 
\end{align}
We shall refer to $\pmb{\rho}_{\mbE^\spadesuit, i}$ as the {\bf radius} of $\mbE^\spadesuit$ at $\sigma_{X, i}$.
In particular, if $G$ is of adjoint type, then the notion of radius coincides with the classical definition discussed in {\it loc.\,cit.}

\SSP
\bde \label{y035} 
Let $\msX$ be as above and $\vec{\pmb{\rho}} := (\pmb{\rho}_i)_{i=1}^r$ an element  of $(\mfc (S) \times \mr{Inj} (\mu, Z))^{\times r}$.
Then, we shall say that a faithful twisted $G$-oper $\mbE^\spadesuit$ on $\msX$ is {\bf of radii $\vec{\pmb{\rho}}$}
if $\pmb{\rho}_{\mbE^\spadesuit, i} = \pmb{\rho}_i$ for every  $i = 1, \cdots, r$.
  \ede
\SSP

For each $i = 1, \cdots, r$, the assignment $\mbE^\spadesuit \mapsto   \widetilde{\pmb{\rho}}_{\mbE^\spadesuit, i}$
determines
 a morphism
\begin{align} \label{ee101}
\mr{ev}^{\mfO \mfp}_{i} : \mfO \mfp_{G, g,r} \migi \overline{\mcI}_\mu([\mfc /Z]).
\end{align}
Moreover, the morphisms $\mr{ev}^{\mfO \mfp}_i$ determine  a morphism 
\begin{align} \label{ee102}
\mr{ev}^{\mfO \mfp} := (\mr{ev}^{\mfO \mfp}_1, \cdots, \mr{ev}^{\mfO \mfp}_r) : \mfO \mfp_{G, g,r} \migi  \overline{\mcI}_\mu ([\mfc /Z])^{\times r}.
\end{align}
Given an element  $\vec{\pmb{\rho}} := (\pmb{\rho}_{i})_{i=1}^r \in (\mfc (k) \times \mr{Inj} (\mu, Z))^{\times r}$,
  we obtain the closed substack
\begin{align}
\mfO \mfp_{G,g,r, \vec{\pmb{\rho}}} := (\mr{ev}^{\mfO \mfp})^{-1} (\prod_{i=1}^r \overline{\mcI}_\mu ([\mfc /Z])_{\pmb{\rho}_i}),
\end{align}
i.e., the substack classifying faithful twisted  $G$-opers of radii $\vec{\pmb{\rho}}$.

\SSP
\bt \label{P04}
 Let $\vec{\pmb{\rho}} := ((\rho_i, \kappa_i))_{i=1}^r \in (\mfc (k) \times\mr{Inj} (\mu, Z))^{\times r}$.
  Then,  isomorphism (\ref{e052}) restricts to  an isomorphism
 \begin{align} \label{e05}
\mfO \mfp_{G, g,r, \vec{\pmb{\rho}}}  \isom  \mfO \mfp_{G_{\mr{ad}}, g,r, (\rho_i)_{i=1}^r}  \times_{\overline{\mfM}_{g,r}} \mfS \mfp_{Z, \delta^\sharp, g,r, (\kappa_i)_{i=1}^r}.
 \end{align}
In particular,  $\mfO \mfp_{G, g,r, \vec{\pmb{\rho}}}$  may be represented by a smooth   Deligne-Mumford stack over $k$ which is flat over $\overline{\mfM}_{g,r}$ of relative dimension 
\begin{align}
(g-1) \cdot \mr{dim} (G) + \frac{r}{2} \cdot (\mr{dim} (G) -\mr{rk} (G)).
\end{align}
  \et
\begin{proof}
The assertion follows from the various definitions involved together with  Theorem \ref{P0012} and  ~\cite[Theorem A]{Wak5}.
\end{proof}

\LSP
\subsection{Forgetting tails.} \label{y0344}
Let $r$ be a  nonnegative integer, $\vec{\pmb{\rho}}$ an element of $(\mfc (S) \times \mr{Inj} (\mu, Z))^{\times r}$, 
 and  $\msX := (X, \{ \sigma_{X, i} \}_{i=1}^{r+1})$   an $(r+1)$-pointed {\it smooth} curve  
  over a $k$-scheme $S$ (hence $S^\mr{log} = S$).
We shall write $\msX_\mr{tail} := (X_\mr{tail}, \{ \sigma_{X_\mr{tail}, i} \}_{i=1}^r)$, i.e., the $r$-pointed curve obtained from $\msX$ by forgetting the last marked point.
In particular,  $X_\mr{tail} = X$ and $\sigma_{X_\mr{tail}, i} = \sigma_{X, i}$ for every  $i \in \{1, \cdots, r \}$.
Also, write 
\begin{align}
\mfO \mfp_{\msX_\mr{tail}, \vec{\pmb{\rho}}}  \left(\text{resp.,}  \ \mfO \mfp_{\msX, (\vec{\pmb{\rho}}, \,\pmb{\varepsilon})}\right)
\end{align}
for the groupoid of  faithful twisted $G$-opers on $\msX_\mr{tail}$ of radii $\vec{\pmb{\rho}}$
(resp., on $\msX$ of radii $(\vec{\pmb{\rho}}, \pmb{\varepsilon})$).
In what follows, we shall construct a canonical functor  $\mfO \mfp_{\msX_\mr{tail}, \vec{\pmb{\rho}}} \migi \mfO \mfp_{\msX, (\vec{\pmb{\rho}}, \,\pmb{\varepsilon})}$.

First, let us consider the case
where $G$ is of adjoint type (hence, 
$\mr{Inj} (\mu, Z) = \{1 \}$ 
and $\pmb{\varepsilon} =\varepsilon$).
Denote by ${^\dagger}\mcE_{B, \msX}$ (resp., ${^\dagger}\mcE_{B, \msX_\mr{tail}}$)
 the $B$-bundle ${^\dagger}\mcE_{\mbB, \hslash, U^\mr{log}/S^\mr{log}}$  introduced  in ~\cite[Eq.\,(210)]{Wak5},  where   the triple ``$(\mbB, \hslash, U^\mr{log}/S^\mr{log})$"  is taken to be  $(B, 1, X^\mr{log}/S)$ (resp., $(B, 1, X_\mr{tail}^\mr{log}/S)$).
Write ${^\dagger}\mcE_{G, \msX} := {^\dagger}\mcE_{B, \msX} \times^B G$ and ${^\dagger}\mcE_{G, \msX_\mr{tail}} := {^\dagger}\mcE_{B, \msX_\mr{tail}} \times^B G$.

Now, let
$\msE^\spadesuit_\mr{tail} := (\mcE_{B, \mr{tail}}, \nabla_{\mr{tail}})$ be a faithful twisted $G$-oper (i.e., a $\mfg$-oper in the sense of ~\cite[Definition 2.1]{Wak5}) on $\msX_\mr{tail}$ of radii $\vec{\rho} \in \mfc (S)^{\times r}$.
According to ~\cite[Proposition 2.19]{Wak5},
there exists a unique pair 
$({^\dagger}\msE_\mr{tail}^{\spadesuit}, \mr{nor}_{\msE^\spadesuit_\mr{tail}})$ consisting of a faithful twisted $G$-oper ${^\dagger}\msE_\mr{tail}^{\spadesuit} := ({^\dagger}\mcE_{B, \msX_\mr{tail}}, {^\dagger}\nabla_{\mr{tail}})$ on $\msX_\mr{tail}$ which is $\{ x_\alpha \}_\alpha$-normal (cf. ~\cite[Definition 2.14]{Wak5})   and an isomorphism of $G$-opers
$\mr{nor}_{\msE^\spadesuit_\mr{tail}} : \msE_\mr{tail}^{\spadesuit} \isom  {^\dagger}\msE_\mr{tail}^{\spadesuit}$.
By using  $\mr{nor}_{\msE^\spadesuit_\mr{tail}}$, we shall identify $\msE^\spadesuit_\mr{tail}$ with ${^\dagger}\msE^{\spadesuit}_\mr{tail}$, i.e., assume that $\msE^\spadesuit_\mr{tail}$ is $\{ x_\alpha \}_\alpha$-normal.

We shall   construct from $\nabla_{\mr{tail}}$ an $S^\mr{log} \left(=S\right)$-connection on ${^\dagger}\mcE_{G, \msX}$.
Let us take 
a pair $\mbU = (U, t)$ consisting of 
an open subscheme $U$ of $X$ with 
$U \cap \mr{Im} (\sigma_{X, r+1}) \neq \emptyset$, $U \cap \mr{Im} (\sigma_{X, i}) = \emptyset$ ($i =1, \cdots, r$) and an element   $t \in \Gamma (U, \mcO_X)$ defining   the closed subscheme $\sigma_{X, r+1}$ such that
  $d \mr{log} (t)$ generates $\Omega_{X^{\mr{log}}/S} |_U$.
Hence, $t \cdot d \mr{log} (t)$ ($= dt$) generates $\Omega_{X_\mr{tail}^\mr{log}/S} |_U$.
According to ~\cite[Eq.\,(211)]{Wak5}, 
  there exists a canonical trivialization
$\mr{triv}_{t \mbU} : {^\dagger}\mcE_{B, \msX_\mr{tail}} |_{U} \isom U \times_k B$ of the $B$-bundle  ${^\dagger}\mcE_{B, \msX_\mr{tail}} |_{U}$ arising from the pair $t\mbU := (U, dt)$.
After a change of  structure group by $B \migiincl G$,
the differential of this trivialization  specifies an $\mcO_U$-linear isomorphism
\begin{align} \label{ee189}
\widetilde{\mcT}_{{^\dagger}\mcE_{G, \msX_\mr{tail}}^\mr{log}/S} |_U \isom \mcT_{U/S} \oplus (\mcO_U \otimes_k \mfg).
\end{align} 
  For  any $h \in T$ and $\alpha \in \Gamma$, 
  denote by $h_\alpha$ 
    the automorphism of $\mfg^{\alpha}$
     given by multiplication by $\alpha (h) \in \mbG_m$.
 The assignment 
   $h \mapsto (h_\alpha)_{\alpha \in \Gamma}$ determines an isomorphism  
 \begin{align} \label{ee110}
 T \isom \prod_{\alpha \in \Gamma} \mr{GL} (\mfg^{\alpha}).
 \end{align}
Write $U^o:= U \setminus \mr{Im} (\sigma_{X, r+1})$.
Also, write $\iota_{U} : U^o  \migiincl  U$ for  the natural open immersion and $\mr{mult}_{\alpha, t}$ (where $\alpha \in \Gamma$) for
the automorphism of $\mcO_{U^o} \otimes_k  \mfg^\alpha$ given by multiplication by $\iota_U^*(t) \in \Gamma (U^o, \mcO_X^\times)$.
The translation by 
 the element of $T (U^o)$ corresponding to $(\mr{mult}_{\alpha, t})_{\alpha \in \Gamma} \in \prod_{\alpha \in \Gamma}\mr{GL} (\mfg^\alpha)$ via  (\ref{ee110}) 
 determines   an automorphism $\mr{mult}_\mbU$ of the trivial $B$-bundle $U^o \times_k B$ on $U^o$.
The composite
\begin{align}
 \mcT_{U^o/S}  &  \xrightarrow{\nabla_{\mr{tail}} |_{U^o}}  \widetilde{\mcT}_{{^\dagger}\mcE_{G, \msX_\mr{tail}}^\mr{log}/S} |_{U^o}
   \xrightarrow{(\ref{ee189})|_{U^o}} \mcT_{U^o/S} \oplus (\mcO_{U^o} \otimes_k \mfg) 
  \xrightarrow{d (\mr{mult}_\mbU)} \mcT_{U^o/S} \oplus (\mcO_{U^o} \otimes_k \mfg) 
\end{align}
 specifies an $S^\mr{log}$-connection on $U^o \times_k G$, and it
  extends uniquely to an $S^\mr{log}$-connection 
 \begin{align}
 \nabla_{\mbU} : \mcT_{X^\mr{log}/S} |_U \migi \mcT_{X^\mr{log}/S} |_U  \oplus (\mcO_U \otimes_k \mfg)
 \end{align}
  on $X^\mr{log} |_U \times_k G$.
 Moreover, the monodromy operator  $\mu^{\nabla_\mbU}_{r+1} 
 \in \mfg (U \times_{X, \sigma_{X, r+1}} S)$
   of $\nabla_\mbU$ at $\sigma_{X, r+1}$ 
    coincides with $- \check{\rho} \in \mfg$ (cf. ~\cite[Proposition 9.2.1]{Fr}).
 The pair $(U \times_k B, \nabla_\mbU)$ forms a faithful twisted $G$-oper on $X^\mr{log} |_U/S$ whose radius at $\sigma_{X, r+1}$ coincides with $\varepsilon$ ($= \chi (- \check{\rho})$).
 Hence, by ~\cite[Proposition 2.19]{Wak5},  there exists uniquely an  $\{ x_\alpha\}_\alpha$-normal  faithful twisted $G$-oper ${^\dagger}\mcE^{\spadesuit}_{\mbU} := ({^\dagger}\mcE_{B, \msX}|_U, {^\dagger}\nabla_\mbU)$ on $X^\mr{log} |_U$ isomorphic to $(U \times_k B, \nabla_\mbU)$.
On the other  hand, since there exists a natural identification ${^\dagger}\mcE_{G, \msX_\mr{tail}} |_{X \setminus \mr{Im} (\sigma_{X, r+1})} = {^\dagger}\mcE_{G, \msX} |_{X \setminus \mr{Im} (\sigma_{X, r+1})}$,
the restriction  $\nabla_{\mr{tail}} |_{X \setminus \mr{Im} (\sigma_{X, r+1})}$ of $\nabla_{\mr{tail}}$ to $X \setminus \mr{Im} (\sigma_{X, r+1})$  may be regarded as 
an $S^\mr{log}$-connection on ${^\dagger}\mcE_{B, \msX} |_{X \setminus \mr{Im} (\sigma_{X, r+1})}$.
Because of the equality  ${^\dagger}\nabla_\mbU |_{U^o} = \nabla_{\mr{tail}} |_{U^o}$,  the $S^\mr{log}$-connection  $\nabla_{\mr{tail}} |_{X \setminus \mr{Im} (\sigma_{X, r+1})}$  and the various $S^\mr{log}$-connections ${^\dagger}\nabla_\mbU$
 (where $\mbU$ ranges over the pairs $\mbU = (U, t)$ as above)
 may be glued together to obtain an $S^\mr{log}$-connection $\nabla_{+ \varepsilon}$ on ${^\dagger}\mcE_{G, \msX}$.
The resulting pair    
 \begin{align} \label{ee731}
 \msE^\spadesuit_{+ \varepsilon} := ({^\dagger}\mcE_{B, \msX}, \nabla_{+\varepsilon})
 \end{align}
   forms a faithful twisted $G$-oper on $\msX$ 
   of radii $(\vec{\rho}, \varepsilon) \in \mfc (S)^{\times (r+1)}$.
The assignment   $\msE_\mr{tail}^\spadesuit \mapsto \msE^\spadesuit_{+ \varepsilon}$  determines 
 a functor 
 $\mfO \mfp_{\msX_\mr{tail}, \vec{\rho}} \migi  \mfO \mfp_{\msX, (\vec{\rho}, \varepsilon)}$.

Next, let us remove the assumption that $G$ is of adjoint type.
Let $\mbE_\mr{tail}^\spadesuit$
   be  a faithful twisted  $G$-oper on $\msX_\mr{tail}$;
 it  corresponds, via (\ref{e052}), to 
  a pair $(\msE_\mr{tail}^\spadesuit, (\msX_\mr{tail}^\mr{tw}, \gamma_\mr{tail}, \mcG_{\mr{tail}}))$ consisting of a faithful twisted $G_\mr{ad}$-oper $\msE^\spadesuit_\mr{tail}$ on $\msX_\mr{tail}$ and a $(Z, \delta^\sharp)$-structure $(\msX_\mr{tail}^\mr{tw}, \gamma_\mr{tail}, \mcG_{\mr{tail}})$ on $\msX_\mr{tail}$.
According to  the above discussion and the discussion in \S\,\ref{y040},
these data induce a faithful twisted $G_\mr{ad}$-oper  $\msE^\spadesuit_{+ \varepsilon}$ on $\msX$ and a $(Z, \delta^\sharp)$-structure $(\msX_{+\delta^\sharp}^\mr{tw}, \gamma_{+ \delta^\sharp}, \mcG_{+ \delta^\sharp})$ on $\msX$ respectively.
The pair $(\msE^\spadesuit_{+ \varepsilon}, (\msX^\mr{tw}_{+ \delta^\sharp}, \gamma_{+ \delta^\sharp}, \mcG_{+\delta^\sharp}))$
corresponds, via (\ref{e052}) again,
to  a faithful twisted $G$-oper
\begin{align}
\mbE_{+ \pmb{\varepsilon}}^\spadesuit 
\end{align}
 on $\msX$.
The assignment $\mbE_\mr{tail}^\spadesuit \mapsto \mbE_{+ \pmb{\varepsilon}}^\spadesuit$
 determines a functor 
$\mfO \mfp_{\msX_\mr{tail}, \vec{\pmb{\rho}}} \migi \mfO \mfp_{\msX, (\vec{\pmb{\rho}}, \,\pmb{\varepsilon})}$.
We shall write
\begin{align} \label{WW1}
\mfO \mfp_{\msX, (\vec{\pmb{\rho}}, \,\pmb{\varepsilon})}^\mr{triv}
\end{align}
for the essential image of this functor.
The above discussion and Proposition \ref{P020} imply 
 the following proposition.

\SSP
\bpr \label{P301} 
 \begin{itemize}
 \item[(i)]
 Let $\msX$ and $\msX_\mr{tail}$ be as above.  Then, the assignment $\mbE_\mr{tail}^\spadesuit \mapsto \mbE^\spadesuit_{+ \pmb{\varepsilon}}$
   determines an equivalence of categories 
 \begin{align}
 \mfO \mfp_{\msX_\mr{tail}, \vec{\pmb{\rho}}} \isom \mfO \mfp_{\msX, (\vec{\pmb{\rho}}, \,\pmb{\varepsilon})}^\mr{triv}.
 \end{align}
 \item[(ii)]
Let $(g, r)$ be a pair of nonnegative integers
with $2g-1+r>0$, and
let   $\vec{\pmb{\rho}}  \in (\mfc (k) \times \mr{Inj} (\mu, Z))^{\times r}$.
 Then, 
 the assignment 
 \begin{align}
 (\msX, (\msX_\mr{tail}, \mbE_\mr{tail}^\spadesuit)) \mapsto (\msX, (\msX, \mbE^\spadesuit_{+ \pmb{\varepsilon}}))
 \end{align}
 determines a fully faithful functor
\begin{align}
\mfM_{g,r+1} \times_{\overline{\mfM}_{g,r}} \mfO \mfp_{G, g,r, \vec{\pmb{\rho}}}  \migi  \mfM_{g,r+1} \times_{\overline{\mfM}_{g,r+1}} \mfO \mfp_{G, g,r+1, (\vec{\pmb{\rho}}, \, \pmb{\varepsilon})}
\end{align}
over  $\mfM_{g,r+1}$.
 \end{itemize}
\epr

\LSP
\subsection{Faithful twisted $G$-opers on the $2$-pointed projective line.} \label{y03ff24}
In this subsection, we shall study the  faithful twisted  $G$-opers  on the $2$-pointed projective line.
Let us keep   the notation in \S\,\ref{yhhh034}.

\SSP
\bpr \label{Pff01} 
 \begin{itemize}
 \item[(i)]
 Let  $\mbE^\spadesuit$ be a faithful twisted $G$-oper  on $\msP'_S$.
 Then, 
  the equality $\pmb{\rho}_{\mbE^\spadesuit, 1} = \pmb{\rho}_{\mbE^\spadesuit, 2}^\veebar$ holds.
 \item[(ii)]
For each $\pmb{\rho} \in \mfc (S) \times \mr{Inj} (\mu, Z)$,
there exists a unique (up to isomorphism)
 faithful twisted $G$-oper $\mbE^\spadesuit$ on $\msP'_S$ with $\pmb{\rho} = \pmb{\rho}_{\mbE^\spadesuit, 1} = \pmb{\rho}_{\mbE^\spadesuit, 2}^\veebar$.
\end{itemize}
\epr
\begin{proof}
By Proposition \ref{P08ff} and Theorem  \ref{P04}, one may assume that $G = G_\mr{ad}$ (hence $\overline{\mcI}_\mu ([\mfc /Z])$ $= \mfc$).
First, we shall consider assertion (i).
Let $\msE^{\spadesuit} := ({^\dagger}\mcE_{B, \msP'_S}, \nabla)$ be an $\{ x_\alpha \}_\alpha$-normal faithful twisted $G$-oper on $\msP'_S$.
The global section $\displaystyle \frac{d s_1}{s_1} \left(= -\frac{ds_2}{s_2}\right) \in \Gamma (P_S, \Omega_{P_S^\mr{log}/S})$ generates globally $\Omega_{P^\mr{log}_S/S}$.
It follows that  the pair $(P_S, \frac{ds_1}{s_1})$  gives 
 a trivialization
${^\dagger}\mcE_{G, \msP'_S} \isom P_S \times_k G$ (cf. ~\cite[Eq.\,(211)]{Wak5}).
By this trivialization,
we have
\begin{align} \label{ER23}
\mr{Hom} (\mcT_{P^\mr{log}_S/S}, \widetilde{\mcT}_{{^\dagger}\mcE_{G, \msP'_S}^\mr{log}/S}) & \cong 
\mr{Hom} (\mcT_{P^\mr{log}_S/S}, \widetilde{\mcT}_{(P_S^\mr{log} \times G)/S}) \\
& \cong  \Gamma (P_S, \Omega_{P_S^\mr{log}/S}\otimes (\mcT_{P_S^\mr{log}/S}  \oplus (\mcO_{P_S} \otimes_k \mfg))) \notag\\
& \cong \Gamma (P_S, \mcO_{P_S} \oplus (\Omega_{P_S^\mr{log}/S} \otimes_k \mfg)). \notag
\end{align}
The $S^\mr{log}$-connection $\nabla'$ on $P_S \times_k G$ corresponding to $\nabla$ via (\ref{ER23}) may be expressed as  
\begin{align} \label{ee290}
\left(1,  \frac{ds_1}{s_1} \otimes v\right) 
 \left(= \left(1, -\frac{d s_2}{s_2} \otimes v\right)\right) 
\in \Gamma (P_S, \mcO_{P_S} \oplus (\Omega_{P_S^\mr{log}/S} \otimes \mfg))
\end{align}
 for some $v \in \mfg (S)$.
Hence, we obtain the following sequence of equalities: 
\begin{align}
\rho_{(P_S \times G, \nabla'), 1} & = \chi \left(\sigma_{\mbP, 1}^* \left(\frac{ds_1}{s_1} \otimes v\right)\right) \\
& = \chi \left(\sigma_{\mbP, 2}^* \left(-\frac{d s_2}{s_2} \otimes v\right)\right) \notag \\
& =  \chi \left(\sigma_{\mbP, 2}^* \left(\frac{d s_2}{s_2} \otimes v\right)\right)^\veebar \notag \\
&  = \rho^\veebar_{(P_S \times G, \nabla'),2}. \notag
\end{align}
This completes the proof of  assertion (i).

Next, we shall  consider assertion (ii).
Write $\mfg^{\mr{ad}(p_1)} := \{ x \in \mfg \  | \  \mr{ad} (p_1)(x) = 0\}$.
Since we have assumed that $\mr{char} (k) =0$ or $\mr{char} (k) =p$ for some prime $p$ satisfying the condition $(*)_{G}$,
the morphism $\mr{kos} : \mfg^{\mr{ad}(p_1)}\migi \mfc$ given by $s \mapsto \chi (p_{-1} +s)$ (for each $s \in \mfg^{\mr{ad}(p_1)}$) is an isomorphism (cf. ~\cite[Lemma 1.2.1]{Ngo}).
Hence, for each $\rho \in \mfc (S)$, there exists 
a unique   $v_0 \in \mfg^{\mr{ad} (p_1)} (S)$
with $\mr{kos}  (v_0) = \rho$.
Let   $\nabla_\rho$  denote the $S^\mr{log}$-connection on $P_S \times_k G$  of the form 
(\ref{ee290}),  where $v$ is taken to be $v:= p_{-1}+v_0$.
 Then, 
 the pair $\msE_\rho^\spadesuit :=  (P_S \times_k B, \nabla_\rho)$ forms a faithful twisted $G$-oper on $\msP'_S$ satisfying the required condition. 
This proves the existence portion of (ii).
The uniqueness portion follows immediately  from the construction of  $\msE_\rho^\spadesuit$.
Thus,   the proof of assertion (ii) is completed.
\end{proof}

\SSP
\bco \label{cff01} 
 Let $\msX$ be a $3$-pointed projective line over $k$.
 Then, we have
 \begin{align}
 \mfO  \mfp^\mr{triv}_{\msX, (\pmb{\rho}_1, \, \pmb{\rho}_2, \, \pmb{\varepsilon})} \cong 
 \begin{cases} \mcB Z& \text{if $\pmb{\rho}_1 = \pmb{\rho}_2^\veebar$;} \\ \emptyset & \text{if otherwise}. \end{cases}
 \end{align}
 \eco
\begin{proof}
The assertion follows from
Propositions  \ref{P301} and \ref{Pff01}.
\end{proof}

\vspace{10mm}
\section{Do'pers and their moduli} \label{SSS444}\SSP

In this section, we study 
 {\it $G$-do'pers} (=
 {\it dormant faithful twisted $G$-opers}) and their moduli, which are central objects in the present paper.
 We prove that the moduli stack of $G$-do'pers is a proper Deligne-Munford stack (cf. Theorem \ref{P05})  and 
satisfies certain factorization properties 
  according to clutching morphisms of the stacks $\overline{\mfM}_{g,r}$ for various pairs $(g,r)$ (cf. Propositions \ref{P0191}, \ref{P0192}, and \ref{P0193}).

Let us fix  a split semisimple algebraic group $(G, T)$ over $k$,
where $T$ denotes a maximal torus $T$ of $G$.
Also,  fix  a Borel subgroup of $G$ defined over $k$ containing $T$.
Denote by $\mfg$, $\mft$, and $\mfb$ the Lie algebras of $G$, $T$, and $B$ respectively.
We shall suppose that $\mr{char}(k) = p>0$ and $G$ satisfies  the condition $(*)_{G}$.

\LSP
\subsection{$G$-do'pers} \label{Sy034}
To begin with, we recall the definition of the {\it $p$-curvature} of a logarithmic  connection (cf., e.g., ~\cite[\S\,3.3]{Wak5}).
Let $S^\mr{log}$ be an fs log scheme over $k$, 
$\mfU^\mr{log}$ a stacky log curve over $S^\mr{log}$, 
and  $(\pi : \mcE \migi \mfU, \nabla)$  
a $G$-bundle  $\pi : \mcE \migi \mfU$ on $\mfU$ paired with an $S^\mr{log}$-connection $\nabla$ on $\mcE$.
If $\partial$ is a logarithmic derivation corresponding to a local section of $\mcT_{\mfU^\mr{log}/S^\mr{log}}$
(resp., $\widetilde{\mcT}_{\mcE^\mr{log}/S^\mr{log}} := (\pi_*(\mcT_{\mcE^\mr{log}/S^\mr{log}}))^G$),
then we shall denote by $\partial^{[p]}$ the $p$-th symbolic power of $\partial$ (i.e., ``$\partial \mapsto \partial^{(p)}$"  asserted  in  ~\cite[Proposition 1.2.1]{Og}),
which is also a logarithmic derivation corresponding to a local section  of $\mcT_{\mfU^\mr{log}/S^\mr{log}}$
(resp., $\widetilde{\mcT}_{\mcE^\mr{log}/S^\mr{log}}$).
Then, 
there exists  
 a unique $\mcO_{\mfU}$-linear morphism $\psi_{(\mcE, \nabla)} : \mcT_{\mfU^\mr{log}/S^\mr{log}}^{\otimes p} \migi \mfg_\mcE \ (\subseteq \widetilde{\mcT}_{\mcE^\mr{log}/S^\mr{log}})$ determined  by $\partial^{\otimes p} \mapsto \nabla (\partial)^{[p]}-\nabla (\partial^{[p]})$.
We shall refer $\psi_{(\mcE, \nabla)}$
 to as the {\bf $p$-curvature map} of $(\mcE, \nabla)$.

\SSP
\bde \label{DD035} 
\begin{itemize}
\item[(i)]
 We shall say that a twisted $G$-oper  
 $\msE^\spadesuit := (\mcE_B, \nabla)$
on $\mfU^\mr{log}/S^\mr{log}$ is {\bf dormant} if
 $\psi_{(\mcE_{B} \times^B G, \nabla)}$ vanishes identically on $\mfU$.
 
 \item[(ii)]
 Let $\msX$ be a pointed nodal curve.
A {\bf $G$-do'per} on $\msX$ is defined to be a  faithful twisted $G$-oper $\mbE^\spadesuit := (\msX^\mr{tw}, \gamma, \msE^\spadesuit)$  on $\msX$ such that $\msE^\spadesuit$ is dormant.
\end{itemize}
  \ede
\SSP

Let $(g, r)$ be a pair of nonnegative integers with $2g-2+r>0$. 
Write
\begin{align} \label{ee215}
\mfO \mfp^{^\mr{Zzz...}}_{G,g,r}
\end{align}
for the closed substack  of $\mfO \mfp_{G,g,r}$ classifying $G$-do'pers.
Also,  write 
\begin{align} \label{ee216}
\pi_{g,r}
 : \mfO \mfp^{^\mr{Zzz...}}_{G,g,r} \migi \overline{\mfM}_{g,r}, \hspace{10mm} \iota_{g,r}
  : \mfO \mfp^{^\mr{Zzz...}}_{G, g,r} \migi \mfO \mfp_{G, g,r}
\end{align}
for the forgetting morphism  and the natural  closed immersion 
 respectively.

\LSP
\subsection{Radii of  $G$-do'pers} \label{SS3}
Let $\mft_\mr{reg}$ denote the open subscheme of  $\mft$   classifying   regular  elements.
Since $G$ is a Chevalley group, all the groups $G$, $T$, and $B$ can be defined over the prime field $\mbF_p$.
It follows that the  $k$-scheme 
 $\mft_\mr{reg}$ comes  from a certain  $\mbF_p$-scheme
  $(\mft_\mr{reg})_{\mbF_p}$ via base-change.
In particular, the Frobenius twist 
$\mft_\mr{reg}^{(1)}$ is isomorphic to 
$\mft_\mr{reg}$ itself, and 
 the relative Frobenius morphism 
 $F_{\mft_\mr{reg}/k}$ may be regarded as a $k$-endomorphism of 
 $\mft_\mr{reg}$.
We shall write
 \begin{align} \label{ee500}
 \mft_{\mr{reg}}^F := \mr{Ker}\left( \mft_\mr{reg} \stackrel{F_{\mft_\mr{reg}/k}, \mr{id}}{\rightrightarrows} \mft_\mr{reg}\right).
 \end{align}
If $W$ denotes the Weyl group of $(G, T)$, then
the natural  $W$-action on $\mft$ restricts to a $W$-action on $ \mft_{\mr{reg}}^F$.
The composite  $\mft \migiincl \mfg \stackrel{\chi}{\migisurj}  \mfc :=  \mfg \ooalign{$/$ \cr $\,/$}\hspace{-0.5mm}G$ induces a closed immersion 
from the resulting quotient   $\mft^F_\mr{reg}/W$ to $\mfc$.
The finite set
\begin{align} \label{ee170}
(\mft_\mr{reg}^F /W) (k)  \left(=\mft^F_{\mr{reg}} (k)/W \right)
\end{align}
may be considered as a subset of $\mfc (k)$.
The $k$-scheme $\mft_\mr{reg}^F/W$
decomposes into 
the disjoint union 
$\coprod_{\rho \in (\mft_\mr{reg}^F /W) (k)} \mr{Spec}(k)_\rho$, where $\mr{Spec}(k)_\rho := \mr{Spec}(k)$.
So we occasionally identify   $\mft_\mr{reg}^F /W$ with the set of  its $k$-rational points in a natural manner.   
The trivial $Z$-action on  $\mft^F_\mr{reg} /W$ gives rise to 
   a closed substack 
   \begin{align} \label{EE900}
  \mfR \mfa \mfd :=  \overline{\mcI}_\mu ([(\mft_\mr{reg}^F/W) /Z])
  \end{align} 
   of $\overline{\mcI}_\mu ([\mfc /Z])$.
 Let us  define a finite set $\Delta$ to be
  \begin{align} \label{ee142}
  \Delta := (\mft_\mr{reg}^F/W) \times \mr{Inj} (\mu, Z).
  \end{align}
  Decomposition (\ref{ee034}) determines  a decomposition
  \begin{align} \label{ee141}
  \mfR \mfa \mfd & = \overline{\mcI}_\mu ((\mft_\mr{reg}^F/W) \times_k \mcB Z) \\
  & = (\mft_\mr{reg}^F/W) \times_k \overline{\mcI}_\mu (\mcB Z) \notag \\
 &  = \coprod_{\rho \in \mft_\mr{reg}^F/W} \mr{Spec}(k)_\rho  \times \coprod_{\kappa \in \mr{Inj} (\mu, Z)} \overline{\mcI}_\mu (\mcB Z)_\kappa  \notag \\
  & = \coprod_{\pmb{\rho} \in \Delta} \mfR \mfa \mfd_{\pmb{\rho}},
  \notag
  \end{align}
 where 
 $\mfR \mfa \mfd_{(\rho, \kappa)} := \mr{Spec}(k)_\rho \times_k \overline{\mcI}_\mu (\mcB Z)_\kappa$.

\SSP
\bpr \label{P0001}
Let $r$ be a positive integer,  $S$ a $k$-scheme, $\msX := (X, \{ \sigma_{X, i} \}_{i=1}^r)$
an $r$-pointed nodal curve over $S$, and 
 $\mbE^\spadesuit$  a $G$-do'per on $\msX$.
Then, for each $i \in \{ 1, \cdots, r \}$,
the radius $\pmb{\rho}_{\mbE^\spadesuit, i} \in \mfc (S) \times \mr{Inj} (\mu, Z)$ of $\mbE^\spadesuit$ at $\sigma_{X, i}$ lies in $\Delta$.
\epr
\begin{proof}
To complete the proof, it suffices to consider the case where
$G = G_{\mr{ad}}$ and $\mbE^\spadesuit$ represents 
 a faithful twisted $G$($= G_\mr{ad}$)-oper $\msE^\spadesuit := (\mcE_B, \nabla)$ which is normal with respect to a fixed choice of generators $\{ x_\alpha \}_\alpha$ in $\mfg^\alpha$'s (cf. ~\cite[Definition 2.14]{Wak5}).
By ~\cite[Proposition 3.13, (i)]{Wak5}, $\rho_{\mbE^\spadesuit, i}$ belongs to $\mfc(k)$.
Hence, 
one may assume, without loss of generality,  that
 $S = \mr{Spec} (k)$ and $k$ is  algebraic closed.
Denote by  $\mu_{i}^\nabla  \left(\in \mfg (k)\right)$ the monodromy operator of $\nabla$ at $\sigma_{X, i}$.
Let us consider a Jordan decomposition  
 $\mu_{i}^\nabla = \mu_s  +\mu_n$  with $\mu_s$ semisimple and $\mu_n$ nilpotent.
 Denote by $\mr{ad} : \mfg \migi \mr{End} (\mfg)$ the adjoint representation of $\mfg$, which is injective and compatible with the respective restricted structures, i.e., $p$-power maps.
 One may find  an isomorphism   of restricted Lie algebras  $\alpha : \mr{End} (\mfg) \isom \mfg \mfl_{\mr{dim} (\mfg)}$
 which sends
 $\alpha (\mr{ad}(\mu_{i}^\nabla))  \left(= \alpha(\mr{ad} (\mu_s)) + \alpha (\mr{ad} (\mu_n))\right)$ to  a Jordan normal form.
 Hence,  $\alpha (\mr{ad} (\mu_s))$ is diagonal and  every  entry of $\alpha (\mr{ad} (\mu_n))$ except  the superdiagonal  is $0$.
Let us observe  the following sequence of equalities:
 \begin{align} \label{ee298}
 \alpha (\mr{ad}(\mu_s)) + \alpha (\mr{ad}(\mu_n)) & = \alpha (\mr{ad} (\mu_s +\mu_n)) = \alpha(\mr{ad} ((\mu_s + \mu_n)^{[p]}))\\ & = \alpha(\mr{ad}(\mu_s + \mu_n))^p = (\alpha (\mr{ad}(\mu_s))+ \alpha (\mr{ad} (\mu_n)))^p,  \notag
 \end{align}
 where the second equality follows from  the assumption 
 that $\msE^\spadesuit$ has
  vanishing $p$-curvature (and   ~\cite[Eq.\,(345)]{Wak5}).
An explicit  computation of   $(\alpha (\mr{ad}(\mu_s))+ \alpha (\mr{ad} (\mu_n)))^p$ shows  (by    (\ref{ee298})) that 
   $\alpha (\mr{ad}(\mu_n)) =0$ ($\Longleftrightarrow  \mu_n=0$), namely, $\mu_{i}^\nabla$ is conjugate to some $v \in \mft$.
Since the map $(-)^{[p]}$ on $\mft$ coincides the $p$-power Frobenius endomorphism,
  the equality  $v = v^{[p]}$ obtained from  (\ref{ee298}) implies that
  $v$ is a Frobenius-invariant element of $\mft$. 
Moreover,  note that (since $\msE^\spadesuit$ is $\{ x_\alpha\}_\alpha$-normal) $\mu_{i}^\nabla$ equals $p_{-1} + u$ for some $u \in  \mfg^{\mr{ad}(p_1)}(k)$ (cf. \S\,\ref{y033} for the definitions of $p_{-1}$ and $p_1$).
 Hence, $\mu_{i}^\nabla$ is regular (cf.  the comment  preceding ~\cite[ Lemma 1.2.3]{Ngo}).
It follows that $v \in \mft^F_\mr{reg} (k)$.
Consequently,  we have $\rho_{\mbE^\spadesuit, i}  \left(= \chi (\mu_{i}^\nabla) = \chi (v)\right) \in (\mft^F_\mr{reg}/W) (k)$, which completes the proof of the assertion. 
 \end{proof}
\SSP

Because of  Proposition \ref{P0001} above,
 the morphism $\mr{ev}^{\mfO \mfp}_i$ ($i =1, \cdots, r$) and $\mr{ev}^{\mfO \mfp}$ restrict to morphisms
\begin{align} \label{ER48}
\mr{ev}_i
: \mfO \mfp^{^\mr{Zzz...}}_{G, g,r} \migi 
\mfR \mfa \mfd
   \left(\stackrel{(\ref{ee141})}{=} \coprod_{\pmb{\rho} \in \Delta}
   \mfR \mfa \mfd_{\pmb{\rho}}  \right), \ \  
\mr{ev}
:  \mfO \mfp^{^\mr{Zzz...}}_{G, g,r} \migi 
\mfR \mfa \mfd^{\times r}   \left( =\!\! \coprod_{(\pmb{\rho}_i)_{i}\in \Delta^{\times r}} \prod_{i=1}^r
\mfR \mfa \mfd_{\pmb{\rho}_i}  \right)
\end{align}
respectively.
Given $\vec{\pmb{\rho}}:= (\pmb{\rho}_i)_{i=1}^r \in \Delta^{\times r}$, we obtain an open and closed substack
\begin{align} \label{ee800}
\mfO \mfp^{^\mr{Zzz...}}_{G, g,r, \vec{\pmb{\rho}}} := \mr{ev}^{-1} (\prod_{i=1}^r
\mfR \mfa \mfd_{\pmb{\rho}_i}) = \mfO \mfp_{G, g,r, \vec{\pmb{\rho}}} \cap \mfO \mfp^{^\mr{Zzz...}}_{G,g,r}
\end{align}
 of
 $\mfO \mfp^{^\mr{Zzz...}}_{G, g,r}$.
 The stack $\mfO \mfp^{^\mr{Zzz...}}_{G, g,r}$ decomposes into the disjoint union
 \begin{align} \label{ee801}
\mfO \mfp^{^\mr{Zzz...}}_{G, g,r} = \coprod_{\vec{\pmb{\rho}} \in \Delta^{\times r}} \mfO \mfp^{^\mr{Zzz...}}_{G, g,r, \, \vec{\pmb{\rho}}}.
 \end{align}
Depending on the choice of $\vec{\pmb{\rho}}$, 
$\mfO \mfp^{^\mr{Zzz...}}_{G, g,r, \, \vec{\pmb{\rho}}}$ may be empty.
The following two lemmas will be applied to conclude 
Theorem \ref{P05}, (iv), described later.

\SSP
\ble \label{W0001} Let 
$S$,  $\msX$, and $\msX_\mr{tail}$ be  as in \S\,\ref{y0344},  and let $\vec{\pmb{\rho}}$  be an element of $\Delta^{\times r}$.
Denote by   $\mfO \mfp_{\msX, (\vec{\pmb{\rho}}, \,\pmb{\varepsilon})}^{^\mr{Zzz...}}$
(resp., $\mfO \mfp_{\msX_\mr{tail}, \vec{\pmb{\rho}}}^{^\mr{Zzz...}}$) the full subgroupoid  of 
$\mfO \mfp_{\msX, (\vec{\pmb{\rho}}, \,\pmb{\varepsilon})}$ (resp., $\mfO \mfp_{\msX_\mr{tail}, \vec{\pmb{\rho}}}$) consisting of $G$-do'pers.
Then, $\mfO \mfp_{\msX, (\vec{\pmb{\rho}}, \,\pmb{\varepsilon})}^{^\mr{Zzz...}}$ is contained in $\mfO \mfp_{\msX, (\vec{\pmb{\rho}}, \,\pmb{\varepsilon})}^{\mr{triv}}$.
In particular, the fully faithful functor 
  $\mfO \mfp_{\msX_\mr{tail}, \vec{\pmb{\rho}}} \migi  \mfO \mfp_{\msX, (\vec{\pmb{\rho}}, \,\pmb{\varepsilon})}$ constructed in \S\,\ref{y0344}
  restricts to an equivalence of categories 
  $\mfO \mfp_{\msX_\mr{tail}, \vec{\pmb{\rho}}}^{^\mr{Zzz...}} \isom \mfO \mfp_{\msX, (\vec{\pmb{\rho}}, \,\pmb{\varepsilon})}^{^\mr{Zzz...}}$.
\ele
\begin{proof}
Let us consider the former assertion.
To complete the proof, it suffices  to consider
the case where $G$ is of adjoint type and $S = \mr{Spec}(R)$ for some $k$-algebra $R$.
Moreover, by considering 
the formal neighborhood $D \left(:=\mr{Spec}(R[[t]])\right)$
 of $\sigma_{X, r+1}$ in $X$,
we can reduce the problem   to proving   the following assertion: 
let $\nabla$ be  a dormant $G$-oper on $D$ having  regular singularity along the divisor $t = 0$ 
 with residue $\varepsilon$ (in the sense of ~\cite[\S\,9.1]{Fr}),  
then $\nabla$ becomes, after the  gauge transformation by some element of $B (R((t)))$,
a $G$-oper on $D$. 
In what follows, let us prove this assertion.

We shall write  $N := [B, B]$, i.e., the unipotent radical of $B$,  and write $\mfn$ for  its Lie algebra.
Also, let us fix a collection of generators $\{ x_\alpha \}_\alpha$ of $\mfg^\alpha$'s, which induces an element $p_{-1} \in \mfg$ as mentioned in \S\,\ref{y033}.
Since we have assumed that  $G$ satisfies the condition $(*)_{G}$, 
 the exponential map $\mr{Exp} : \mfn \isom  N $ resulting from    ~\cite[Proposition 1.31]{Wak5} can be defined.
By means of this, we can  apply 
 an argument similar to the argument  in the proof of ~\cite[Proposition 9.2.1]{Fr} (for $\check{\lambda} =0$) to our discussion  even though we are working in positive characteristic (see also ~\cite[Proposition 2.19]{Wak5}).
It follows that, 
by  the gauge transformation by some element of $B (R((t)))$,
$\nabla$ may be brought to  a connection of the form 
$\frac{d}{dt} + p_{-1}+ v(t) + \frac{u}{t}$
 for some $ v(t) \in \mfb(R[[t]])$ and $u \in \mfn$.
The mod $t$ reduction of the  $p$-curvature of $\nabla$  with respect to this expression is given by $u^{[p]} - u$, which
 must be equal to 
$0$  because $\nabla$ has vanishing $p$-curvature.
But, 
 since 
  $u \in \mfn$ (which implies  $u^{[p]}=0$),  the equality  $u = 0$ holds.
Therefore,  
the connection $\nabla$ is gauge equivalent to $\frac{d}{dt} + p_{-1}+ v(t)$ and forms
 a $G$-oper on $D$.
This completes the proof of the former assertion.
The latter assertion follows from the former assertion and Proposition \ref{P301},  (i).
 \end{proof}

\SSP
\ble \label{W0002} 
The unique (up to isomorphism)
faithful twisted $G$-oper $\mbE^\spadesuit$ on $\msP'_{k} \left(:=\msP'_{\mr{Spec}(k)}  \right)$ with $\pmb{\rho}_{\mbE^\spadesuit, 1} = \pmb{\rho}_{\mbE^\spadesuit, 2}^\veebar \in \Delta$ (cf. Proposition \ref{Pff01},  (ii))
 is dormant.
\ele
\begin{proof}
To complete the proof, it suffices to consider the case where
$G= G_\mr{ad}$, $k$ is algebraically closed, 
and $\mbE^\spadesuit$ is normal with respect to a fixed collection of generators $\{x_\alpha \}_\alpha$ of $\mfg^\alpha$ (cf. ~\cite[Definition 2.14]{Wak5}).
Denote by $\rho \left(\in \mft^F_\mr{reg}/W\right)$ the radius  of $\mbE^\spadesuit$ at $\sigma_{\mbP, 1}$.
Let $p_1$ and $p_{-1}$ be as in  \S\,\ref{y033}.
Then, the  connection $\nabla$ defining  $\mbE^\spadesuit$ may be expressed as (\ref{ee290}), where  the element $v$ is taken to be   $v = p_{-1} + v_0$ for a unique  $v_0 \in \mfg^{\mr{ad}(p_1)}$ with $\chi (v) = \rho$ (cf. ~\cite[Lemma 1.2.1]{Ngo}).
It follows from 
~\cite[Lemma 1.2.3]{Ngo}
that 
$\nabla$ becomes, after the gauge transformation by some element of $G$,
a connection of the form (\ref{ee290}) with $v = \widetilde{\rho}$ for some lifting  $\widetilde{\rho} \in \mft_\mr{reg}^F$ of $\rho$.
The $p$-curvature
of  $\nabla$ with respect to this form  is given by  
  $\displaystyle
 \left(s_1\frac{d}{ds_1} \right)^{\otimes p} 
 \mapsto 1 \otimes (\widetilde{\rho}^{[p]}-\widetilde{\rho})$.
But, the equality $\widetilde{\rho}^{[p]} = \widetilde{\rho}$ holds since $\widetilde{\rho}$ belongs to $\mft_\mr{reg}^F$, so
$\nabla$ has vanishing $p$-curvature.
 That is to say,  $\mbE^\spadesuit$ turns out to be  dormant.
  \end{proof}

\LSP
\subsection{The moduli space  of  $G$-do'pers} \label{W003}
In what follows, let us describe  an  assertion concerning  
the structure of the moduli stack  of $G$-do'pers.
This  assertion contains a part of Theorem \ref{Tqqq}.

\SSP
\bt \label{P05}
\begin{itemize}
\item[(i)]
 Isomorphism  (\ref{e052}) restricts to an isomorphism
 \begin{align} \label{eq113}
\mfO \mfp^{^\mr{Zzz...}}_{G, g,r}  \isom  \mfO \mfp^{^\mr{Zzz...}}_{G_{\mr{ad}}, g,r}  \times_{\overline{\mfM}_{g,r}} \mfS \mfp_{Z, \delta^\sharp, g,r}.
 \end{align}
 In particular, the morphism $\mr{op}_\mr{ad}$ (cf. (\ref{e051})) restricts to a finite, flat, and generically \'{e}tale morphism 
 \begin{align} \label{ER44}
 \mr{op}_\mr{ad}^{^\mr{Zz.}} : \mfO \mfp_{G,g,r}^{^\mr{Zzz...}} \migi \mfO \mfp_{G_\mr{ad},g,r}^{^\mr{Zzz...}}, 
 \end{align}
 and  the following commutative diagram is  cartesian: 
 \begin{align} \label{ee214}
 \begin{CD}
 \mfO \mfp_{G,g,r}^{^\mr{Zzz...}} @> \mr{op}_\mr{ad}^{^\mr{Zz.}} >> \mfO \mfp_{G_\mr{ad},g,r}^{^\mr{Zzz...}}
 \\
 @V \iota_{g,r} VV @VV \iota_{g,r}V
 \\
 \mfO \mfp_{G,g,r}@>> \mr{op}_\mr{ad} > \mfO \mfp_{G_\mr{ad},g,r}.
 \end{CD}
 \end{align}
 Moreover, for each $\vec{\pmb{\rho}} := ((\rho_i, \kappa_i))_{i=1}^r \in \Delta^{\times r}$,
    isomorphism  (\ref{eq113})
   restricts to an isomorphism
  \begin{align} \label{eq114}
\mfO \mfp^{^\mr{Zzz...}}_{G, g,r, \vec{\pmb{\rho}}}  \isom  \mfO \mfp^{^\mr{Zzz...}}_{G_{\mr{ad}}, g,r, (\rho_i)_{i=1}^r}  \times_{\overline{\mfM}_{g,r}} \mfS \mfp_{Z, \delta^\sharp,  g,r, (\kappa_i)_{i=1}^r}.
 \end{align}
 \item[(ii)]
Both  $\mfO \mfp^{^\mr{Zzz...}}_{G, g,r}$ and  $\mfO \mfp^{^\mr{Zzz...}}_{G, g,r, \vec{\pmb{\rho}}}$ (for every  $\vec{\pmb{\rho}} \in \Delta^{\times r}$) may be represented by (possibly empty)  proper   Deligne-Mumford stacks over $k$
 which are  finite over over $\overline{\mfM}_{g,r}$.
Moreover,  $\mfO \mfp^{^\mr{Zzz...}}_{G, g,r}$ is nonempty and 
 has an irreducible component that dominates $\overline{\mfM}_{g,r}$.
 \item[(iii)]
 Let us assume further  that $G$ satisfies the condition $(**)_{G}$ described in Introduction.
 Then, 
 $\mfO \mfp^{^\mr{Zzz...}}_{G_{}, g,r, \vec{\pmb{\rho}}}$ is (finite and) flat  over the points of $\overline{\mfM}_{g,r}$ classifying pointed  totally degenerate curves (cf. ~\cite[Definition 7.15]{Wak5}  for the definition of a pointed totally degenerate curve).
 Moreover, $\mfO \mfp^{^\mr{Zzz...}}_{G_{}, g,r, \vec{\pmb{\rho}}}$ is generically \'{e}tale  over $\overline{\mfM}_{g,r}$ (i.e., any irreducible component that dominates $\overline{\mfM}_{g,r}$ admits a dense open subscheme which is \'{e}tale  over $\overline{\mfM}_{g,r}$)  and has  generic stabilizer isomorphic to the center $Z$ of $G$.
 \item[(iv)]
 The following assertion  holds (without the assumption imposed in $\mr{(iii)}$):
 \begin{align} 
 \mfO  \mfp^{^\mr{Zzz...}}_{G, 0, 3, (\pmb{\rho}_1, \pmb{\rho}_2, \pmb{\varepsilon})} \cong 
 \begin{cases} \mcB Z & \text{if $\pmb{\rho}_1 \in \Delta$  and $\pmb{\rho}_1 = \pmb{\rho}_2^\veebar$};  \\ \emptyset & \text{if otherwise}. \end{cases}
 \end{align}
 \end{itemize}
  \et
\begin{proof}
Assertion (i) follows from the various definitions involved together with Theorem \ref{P08f}.
Assertion (ii) follows from assertion (i) and ~\cite[Theorem C]{Wak5}.
Assertion (iii) follows from assertion (i), 
 Theorem \ref{P08f} (in the present paper),  Theorem G in {\it loc.\,cit.}, and ~\cite[Chap.\,II,   Theorem 2.8]{Mzk2}.
Finally, assertion (iv) follows from Corollary \ref{cff01}, Proposition \ref{P0001}, and Lemmas \ref{W0001} and  \ref{W0002}.
\end{proof}
\SSP

By the above theorem, it makes sense to speak of the {\it generic degree}
\begin{align} \label{Ww600}
 \mr{deg}^\mr{gen}(\mfO \mfp^{^\mr{Zzz...}}_{G, g,r, \vec{\pmb{\rho}}}/\overline{\mfM}_{g,r}) 
   \left(\in \mbQ\right)
\end{align}
of 
$\mfO \mfp^{^\mr{Zzz...}}_{G, g,r, \vec{\pmb{\rho}}}$ (for each $\vec{\pmb{\rho}} := ((\rho_i, \kappa_i))_{i=1}^r \in \Delta^{\times r}$)
over $\overline{\mfM}_{g,r}$.
Moreover,  isomorphism (\ref{eq114}) yields  
 the equality 
\begin{align} \label{Ww550}
 & \ \mr{deg}^\mr{gen}(\mfO \mfp^{^\mr{Zzz...}}_{G, g,r, \vec{\pmb{\rho}}}/\overline{\mfM}_{g,r}) 
 =\mr{deg}^\mr{gen}(\mfO \mfp^{^\mr{Zzz...}}_{G_\mr{ad}, g,r, (\rho_i)_{i=1}^r}/\overline{\mfM}_{g,r}) \cdot \mr{deg} (\mfS \mfp_{Z, \delta^\sharp,g,r, (\kappa_i)_{i=1}^r}/\overline{\mfM}_{g,r}). 
\end{align}
Then, we obtain  the following three propositions. 

\SSP
\bpr \label{P0191} 
 Assume that $G$ satisfies the condition $(**)_{G}$.
 Let  $g_1$, $g_2$, $r_1$, and $r_2$ be nonnegative integers with $2g_j -1 +r_j > 0$ ($j =1, 2$),
and let $\vec{\pmb{\rho}}_j \in \Delta^{\times r_j}$.
 Write $g := g_1 + g_2$, $r = r_1 + r_2$.
\begin{itemize}
\item[(i)]
For each  $\pmb{\rho} \in \Delta$, there exists a  morphism 
\begin{align}
\Phi^{^\mr{Zz.}}_{\mr{tree}, \, \pmb{\rho}} :  \ & \mfO \mfp^{^\mr{Zzz...}}_{G, g_1,r_1 +1, (\vec{\pmb{\rho}}_{1}, \pmb{\rho})}  \times_{\mr{ev}_{r_1 +1}, 
\mfR \mfa \mfd_{\pmb{\rho}}, \mr{ev}_{r_2 +1}} \mfO \mfp^{^\mr{Zzz...}}_{G, g_2,r_2 +1, (\vec{\pmb{\rho}}_{2}, \pmb{\rho})}   \migi \mfO \mfp^{^\mr{Zzz...}}_{G, g, r, (\vec{\pmb{\rho}}_{1}, \vec{\pmb{\rho}}_{2})}  
\end{align}
obtained by gluing together   two $G$-do'pers    along  the fibers over  the  respective last marked points of  the underlying curves.
 Moreover, this morphism makes  the following square diagram commute:
\begin{align}
\begin{CD}
\displaystyle 
\mfO \mfp^{^\mr{Zzz...}}_{G, g_1,r_1 +1, (\vec{\pmb{\rho}}_{1}, \pmb{\rho})}  \times_{
\mfR \mfa \mfd_{\pmb{\rho}}} \mfO \mfp^{^\mr{Zzz...}}_{G, g_2,r_2 +1, (\vec{\pmb{\rho}}_{2}, \pmb{\rho})} @>
  \Phi^{^\mr{Zz.}}_{\mr{tree}, \pmb{\rho}} >> \mfO \mfp^{^\mr{Zzz...}}_{G, g, r, (\vec{\pmb{\rho}}_{1}, \vec{\pmb{\rho}}_{2})}
\\
@V 
\pi_{g_1,r_1} \times \pi_{g_2,r_2} VV @VV \pi_{g, r} V
\\
\overline{\mfM}_{g_1, r_1+1} \times_k \overline{\mfM}_{g_2, r_2+1} @>> \Phi_{\mr{tree}}> \overline{\mfM}_{g, r}.
\end{CD}
\end{align}
\item[(ii)]
The following equality holds:
\begin{align}
\mr{deg}^\mr{gen}(\mfO \mfp^{^\mr{Zzz...}}_{G, g, r, (\vec{\pmb{\rho}}_{1}, \vec{\pmb{\rho}}_{2})}/\overline{\mfM}_{g,r})
= |Z| \cdot \sum_{\pmb{\rho} \in \Delta} \prod_{j=1}^2 \mr{deg}^\mr{gen}(\mfO \mfp^{^\mr{Zzz...}}_{G, g_j,r_j +1, (\vec{\pmb{\rho}}_{j}, \pmb{\rho})}/\overline{\mfM}_{g_j, r_j+1}).
\end{align}
 \end{itemize}
  \epr
\begin{proof}
The assertions of the case where $G = G_\mr{ad}$ follow from ~\cite[Theorem 7.13]{Wak5}.
The general case can be reduced to  that case because of Proposition \ref{P08}, Theorem \ref{P05}, (i),  and equality (\ref{Ww550}).
\end{proof}

\SSP
\bpr \label{P0192}
 Assume that $G$ satisfies the condition $(**)_{G}$.
  Let $(g, r)$ be a pair of  nonnegative integers with $2g + r >0$, and let $\vec{\pmb{\rho}} \in \Delta^{\times r}$.
\begin{itemize}
\item[(i)]
For each $\pmb{\rho} \in \Delta$, 
 there exists  
 a morphism
 \begin{align}
 \Phi^{^\mr{Zz.}}_{\mr{loop}, \pmb{\rho}} : \mr{Ker}\left(\mfO \mfp^{^\mr{Zzz...}}_{G, g,r+2, (\vec{\pmb{\rho}}, \pmb{\rho}, \pmb{\rho}^\veebar)} \stackrel{\mr{ev}_{r+1}, \mr{ev}_{r+2}}{\rightrightarrows} 
 \mfR \mfa \mfd_{\pmb{\rho}}\right) \migi \mfO \mfp^{^\mr{Zzz...}}_{G, g+1,r, \vec{\pmb{\rho}}}
 \end{align}
 obtained by gluing each $G$-do'per along the fibers over the last two marked 
 points of the underlying curve.
 Moreover, this morphism makes the following square diagram commute:
 \begin{align}
 \begin{CD}
 \displaystyle 
\mr{Ker}\left(\mfO \mfp^{^\mr{Zzz...}}_{G, g,r+2, (\vec{\pmb{\rho}}, \pmb{\rho}, \pmb{\rho}^\veebar)} 
 \stackrel{\mr{ev}_{r+1}, \mr{ev}_{r+2}}{\rightrightarrows} 
  \mfR \mfa \mfd_{\pmb{\rho}}\right)   @> 
  \Phi^{^\mr{Zz.}}_{\mr{loop}, \pmb{\rho}} >> \mfO \mfp^{^\mr{Zzz...}}_{G, g+1,r, \vec{\pmb{\rho}}} 
 \\
 @V 
  \pi_{g,r+2} VV @VV \pi_{g+1, r} V
 \\
 \overline{\mfM}_{g, r+2} @>> \Phi_{\mr{loop}}> \overline{\mfM}_{g+1, r}.
\end{CD}
\end{align}
 \item[(ii)]
 The following equality holds:
 \begin{align}
 \mr{deg}^\mr{gen}(\mfO \mfp^{^\mr{Zzz...}}_{G, g+1,r, \vec{\pmb{\rho}}}/\overline{\mfM}_{g+1, r})
 =
 |Z| \cdot \sum_{\pmb{\rho} \in \Delta} \mr{deg}^\mr{gen} (\mfO \mfp^{^\mr{Zzz...}}_{G, g,r+2, (\vec{\pmb{\rho}}, \pmb{\rho}, \pmb{\rho}^\veebar)}/\overline{\mfM}_{g, r+2}).
 \end{align}
  \end{itemize}
  \epr
\begin{proof}
The assertions of the case where $G = G_\mr{ad}$ follow from ~\cite[Theorem 7.13]{Wak5}.
The general case can be reduced to that case because of Proposition \ref{Pw081}, Theorem \ref{P05}, (i), and equality (\ref{Ww550}).
\end{proof}

\SSP
\bpr \label{P0193} 
 Assume that $G$ satisfies the condition $(**)_{G}$.
  Let $(g, r)$ be a pair of nonnegative integers with $2g -1+r >0$, and let $\vec{\pmb{\rho}} \in \Delta^{\times r}$.
\begin{itemize}
\item[(i)]
  There exists  an isomorphism
 \begin{align}
 \mfM_{g, r+1} \times_{\overline{\mfM}_{g,r}}\mfO \mfp^{^\mr{Zzz...}}_{G, g,r, \, \vec{\pmb{\rho}}}  \isom \mfM_{g,r+1} \times_{\overline{\mfM}_{g, r+1}} \mfO \mfp^{^\mr{Zzz...}}_{G, g,r+1, (\vec{\pmb{\rho}}, \, \pmb{\varepsilon})} 
 \end{align}
over $\mfM_{g, r+1}$.
\item[(ii)]
The following equality holds:
\begin{align}
\mr{deg}^\mr{gen}( \mfO \mfp^{^\mr{Zzz...}}_{G, g,r+1, (\vec{\pmb{\rho}}, \, \pmb{\varepsilon})} /\overline{\mfM}_{g,r+1})
=
\mr{deg}^\mr{gen}(\mfO \mfp^{^\mr{Zzz...}}_{G, g,r, \, \vec{\pmb{\rho}}}/\overline{\mfM}_{g,r}).
\end{align}
 \end{itemize}
  \epr
\begin{proof}
The assertions  follow from Proposition \ref{P020}, Theorem \ref{P05}, (i), and equality (\ref{Ww550}).
\end{proof}
\SSP

\begin{rema} \label{R036}
At the time of writing the present paper, the author does not know to what extent  one  can weaken  the condition  $(**)_{G}$ imposed in Theorem \ref{P05} (iii),  and Propositions \ref{P0191}, \ref{P0192}, and  \ref{P0193}.
 \end{rema}

\vspace{10mm}
\section{The virtual fundamental class on the moduli of do'pers} \label{SSS555}\SSP

In this section, we construct a {\it perfect obstruction theory} for the moduli  stack $\mfO \mfp^{^\mr{Zzz...}}_{G,g,r}$ (cf. Theorem \ref{P01}).
As a result, we obtain
a {\it virtual fundamental class} 
on that moduli (cf. (\ref{ee226})).
This virtual class will be used to construct a CohFT (forming, in fact, a $2$d TQFT) of $G$-do'pers.

Let us keep the notation at the beginning of \S\,\ref{SSS444}.

\LSP
\subsection{General definition of a perfect obstruction theory} \label{SS3}
First, we shall recall  from ~\cite{BF1} the notions of a perfect obstruction theory and the  virtual fundamental class associated with it.

Let $\mfX$ be  a separated Deligne-Mumford stack,  locally of finite type  over  $k$.
Denote by $D (\mcO_\mfX)$ the derived category of the category $\mr{Mod} (\mcO_\mfX)$ of $\mcO_\mfX$-modules and by $L_\mfX^\bullet \in \mr{Ob} (D (\mcO_\mfX))$ the cotangent complex  of $\mfX$ relative to $k$.

For  a morphism $E^0_\mr{fl} \migi E^1_\mr{fl}$  of abelian sheaves on the big  fppf site $\mfX_{\mr{fl}}$ of  $\mfX$,  one obtains  the quotient stack $[E^1_\mr{fl} / E^0_\mr{fl}]$.
That is to say, for an object $T \in \mr{Ob} (\mfX_\mr{fl})$ the groupoid  $[E^1_\mr{fl} / E^0_\mr{fl}] (T)$ of sections over $T$ is the category of pairs $(P, f)$, where $P$ is an $E^0_\mr{fl}$-torsor over $T$ and $f$ is an $E^0_\mr{fl}$-equivariant morphism $P \migi E^1_\mr{fl} |_T$ of sheaves on $\mfX_\mr{fl}$.
If  $E^\bullet_\mr{fl}$ is a complex of arbitrary length of abelian sheaves on  $\mfX_\mr{fl}$, we shall write  
\begin{align}
h^1/h^0 (E^\bullet_\mr{fl}) := [Z^1/ C^0],
\end{align} 
where $Z^1 := \mr{Ker} (E^1_\mr{fl} \migi E^2_\mr{fl})$, $C^0 := \mr{Coker} (E^{-1}_\mr{fl}\migi E^0_\mr{fl})$.
Denote by $\mfN_\mfX$ the stack defined to be
$\mfN_\mfX := h^1/ h^0 (((L_\mfX^\bullet)_\mr{fl})^\vee)$, where for each complex $E^\bullet$ on (the small \'{e}tale site of) $\mfX$, we shall write $E^\bullet_\mr{fl}$ for the complex on $\mfX_\mr{fl}$ associated with $E^\bullet$.

Let us take a
diagram
\begin{align} \label{ee211}
\xymatrix{
U\ar[r]^{\iota} \ar[d]_{\pi} &  M
\\
\mfX,  &
}
\end{align}
where $U$ denotes
  an affine scheme of finite type over $k$, $M$ denotes  a smooth affine scheme of finite type over $k$,  $\iota$ is  a closed immersion, and
   $\pi$ is an \'{e}tale morphism.
 Let  $I$ denote the ideal sheaf on $M$  defining the closed subscheme  $U$, and we consider the natural morphism $I/I^2  \migi \iota^*(\Omega_{M/k})$ as a complex $[I/I^2  \migi \iota^*(\Omega_{M/k})]$ concentrated in degrees $-1$ and $0$.
Then, there exists a natural   quasi-isomorphism 
\begin{align} \label{ee210}
\phi : L_\mfX^\bullet |_U \isom [I/I^2  \migi \iota^*(\Omega_{M/k})].
\end{align}
Write $T_M := \mcS pec (\mr{Sym}_{\mcO_M} (\Omega_{M/k}))$ (i.e., the total space of the tangent bundle $\mcT_{M/k}$)  and $N_{U/M} := \mcS pec (\mr{Sym}_{\mcO_U} (I/I^2))$.
Note that  $N_{U/M}$ admits an $\iota^* (T_M)$-action induced from the morphism $I/I^2  \migi \iota^*(\Omega_{M/k})$.
 Quasi-isomorphism  (\ref{ee210}) gives rise to an isomorphism
\begin{align} \label{ee212}
[N_{U/M}/ \iota^*(T_M)] \isom \mfN_\mfX |_U.
\end{align}

The {\bf intrinsic normal cone}  of $\mfX$ (cf. ~\cite[Definition 3.10]{BF1}) is defined as a unique closed substack $\mfC_\mfX$  of $\mfN_\mfX$ determined by the condition that
if 
we are given a diagram as in (\ref{ee211}), 
then $\mfC_\mfX |_U$ may be identified, via (\ref{ee212}), with  
the closed substack $[C_{U/M}/\iota^*(T_M)]$, where $C_{U/M}$ denotes the normal cone $\mcS pec (\bigoplus_{n \geq 0} I^n/I^{n+1})$ of $U$ in $M$.  
A {\bf perfect obstruction theory} for $\mfX$ (cf. ~\cite[Definitions 
  4.4  and 5.1]{BF1}) is a morphism $\phi : E^\bullet \migi L_\mfX^\bullet$ in $D (\mcO_\mfX)$ satisfying the following conditions:
\begin{itemize}
\item[$\bullet$]
$h^0 (\phi)$ is an isomorphism and $h^{-1} (\phi)$  is surjective.
\item[$\bullet$]
$E^\bullet$ is of perfect amplitude contained in $[-1, 0]$, i.e., is  locally isomorphic (in $D (\mcO_\mfX)$) to a complex $[E^{-1} \migi E^0]$ of locally free sheaves of finite rank.
\end{itemize} 

Let $\phi : E^\bullet \migi L_\mfX^\bullet$ be a perfect obstruction theory for $\mfX$.
 The {\bf virtual dimension} of $\mfX$ with respect 
to $\phi : E^\bullet \migi L_\mfX^\bullet$ is a well-defined locally constant function on $\mfX$, denoted by $\mr{rk} (E^\bullet)$,  defined in such a way that if 
$E^\bullet$ is locally  written as a complex of vector bundles $[E^{-1} \migi E^0]$,
then $\mr{rk} (E^\bullet) := \mr{dim} (E^0) -\mr{dim} (E^{-1})$.
Let us suppose further that  $\mr{rk} (E^\bullet)$ is constant
and that 
$E^\bullet$ has a global resolution, i.e., has a morphism of vector bundles $F^\bullet := [F^{-1} \migi F^0]$ considered as a complex concentrated in degrees $-1$ and $0$  together with  an isomorphism $F^\bullet \isom E^\bullet$ in $D (\mcO_\mfX)$.
Since $h^1/h^0 ((E_\mr{fl}^\bullet)^\vee)$ is isomorphic to $[(F_\mr{fl}^{-1})^\vee/(F_\mr{fl}^0)^\vee]$, the relative affine space $\mfV$ associated with  $(F^{-1}_\mr{fl})^\vee$ specifies a global presentation $\mfV \migi h^1/h^0 ((E_\mr{fl}^\bullet)^\vee)$.
Let $\mfW$ be the fiber product
\begin{align}
\begin{CD}
\mfW @>>> \mfV 
\\
@VVV @VVV
\\
\mfC_\mfX @>>> h^1/h^0 ((E_\mr{fl}^\bullet)^\vee),
\end{CD}
\end{align}
where the lower horizontal arrow denotes the composite of the closed immersion $\mfC_\mfX \migi \mfN_\mfX$ and  the morphism  $\mfN_\mfX \migi h^1/h^0 ((E_\mr{fl}^\bullet)^\vee)$ induced by $\phi$,  being  a closed immersion (cf. ~\cite[Theorem 4.5]{BF1}).
We define the {\bf virtual fundamental class $[\mfX, E^\bullet]^\mr{virt}$} to be the intersection of $\mfW$ with the zero section $0_\mfV : \mfX \migi \mfV$, i.e., 
\begin{align} 
[\mfX, E^\bullet]^\mr{virt} := 0^!_\mfV [\mfW], 
\end{align}
which is an element of $A_{\mr{rk} (E^\bullet)} (\mfX)_\mbQ$.
This class is independent of the global resolution $F^\bullet$ used to construct it.
If $\mfX$ is smooth, then the virtual fundamental class $[\mfX, L_\mfX^\bullet]^\mr{virt}$ is equal to the usual fundamental class $[\mfX]$.

\LSP
\subsection{The perfect obstruction theory for the moduli space of $G$-do'pers.} \label{SS3}
Let $(g, r)$ be a pair of nonnegative integers with $2g-2+r>0$.
The goal of this subsection is to  construct a perfect obstruction theory for $\mfO \mfp^{^\mr{Zzz...}}_{G,g,r}$.

Denote by $\mfC \mfo \mfn \mfn_{G, g,r}$  the category  classifying
triples
$(\msX, \mcE, \nabla)$
consisting of 
 an $r$-pointed stable curve $\msX := (X/S, \{ \sigma_{X, i} \}_{i=1}^r)$ of genus $g$    over a  $k$-scheme $S$,
  a $G$-bundle $\mcE$ on $X$, 
  and an  $S^\mr{log}$-connection $\nabla$ on $\mcE$.
This category may be represented by a smooth algebraic stack over $k$.
 Also, denote by $\mfC \mfo \mfn \mfn_{G, g,r}^{^\mr{Zzz...}}$ the closed substack of $\mfC \mfo \mfn \mfn_{G, g,r}$ classifying 
 triples $(\msX, \mcE, \nabla)$
 with  $\psi_{(\mcE, \nabla)} =0$.
The assignment from each pair $(\msX, \mbE^\spadesuit)$ (where  $\msX$ denotes a pointed stable curve and $\mbE^\spadesuit$ denotes a  faithful twisted $G$-oper  on  $\msX$) to the triple 
$(\msX, {^\dagger}\mcE_{G, \msX}, {^\dagger}\nabla)$ (where
${^\dagger}\nabla$ denotes the $S^\mr{log}$-connection defining the $\{x_\alpha \}_\alpha$-normalization of $\mbE^\spadesuit$ asserted in ~\cite[Proposition 2.19]{Wak5})
  determines  morphisms
  \begin{align}
 \xi : \mfO \mfp_{G,g,r} \migi \mfC \mfo \mfn \mfn_{G,g,r}, \hspace{5mm}  \xi^{^\mr{Zz.}} : \mfO \mfp^{^\mr{Zzz...}}_{G, g,r} \migi \mfC \mfo \mfn \mfn_{G,g,r}^{^\mr{Zzz...}} 
  \end{align}
over $\overline{\mfM}_{g,r}$.
The following square diagram is commutative and cartesian:
\begin{align} \label{ee213}
\begin{CD}
\mfC \mfo \mfn \mfn_{G, g,r}^{^\mr{Zzz...}} @> \iota^{\mfC \mfo}_{g,r} >> \mfC \mfo \mfn \mfn_{G, g,r}
\\
@A  \xi^{^\mr{Zz.}} AA @AA \xi A
\\
\mfO \mfp^{^\mr{Zzz...}}_{G, g,r} @>> \iota_{g,r}  > \mfO \mfp_{G, g,r},
\end{CD}
\end{align}
 where the upper horizontal arrow $\iota^{\mfC \mfo}_{g,r}$ denotes the natural closed immersion.
 Denote by $\mcI_G$ (resp.,  $\widehat{\mcI}_{G}$)
the ideal sheaf on $\mfO \mfp_{G, g,r}$ (resp., $\mfC \mfo \mfn \mfn_{G,g,r}$) defining the closed substack $\mfO \mfp_{G, g,r}^{^\mr{Zzz...}}$ (resp., $\mfC \mfo \mfn \mfn_{G,g,}^{^\mr{Zzz...}}$).
Then, diagram  (\ref{ee213}) induces
 a commutative  diagram of coherent sheaves on $\mfO \mfp_{G, g,r}^{^\mr{Zzz...}}$: 
\begin{align} \label{ee220}
\begin{CD}
\xi^{^\mr{Zz.}*} (\widehat{\mcI}_{G}/\widehat{\mcI}_{G}^2) @> \underline{\iota}^{\mfC \mfo}_{g,r}>> (\iota^{\mfC \mfo}_{g,r}\circ \xi^{^\mr{Zz.}})^*(\Omega_{\mfC \mfo \mfn \mfn_{G, g,r}/k})
\\
@V \underline{\xi}^{^\mr{Zz.}} VV @VV \underline{\xi}V
\\
\mcI_G/\mcI_G^2 @>> \underline{\iota}_{g,r} > \iota_{g,r}^* (\Omega_{\mfO \mfp_{G, g,r}/k}).
\end{CD}
\end{align}
Since  (\ref{ee213}) is cartesian, $\underline{\xi}^{^\mr{Zz.}}$ is verified to be surjective. 
Let us consider the composite $\underline{\iota}_{g,r} \circ \underline{\xi}^{^\mr{Zz.}}$ ($=\underline{\xi} \circ  \underline{\iota}^{\mfC \mfo}_{g,r}$) as a complex
\begin{align}
E^\bullet := [\xi^{^\mr{Zz.}*}(\widehat{\mcI}_{G}/\widehat{\mcI}_{G}^2) \migi \iota_{g,r}^{*} (\Omega_{\mfO \mfp_{G, g,r}/k})]
\end{align}
concentrated in degrees $-1$ and $0$.
The Deligne-Mumford stack $\mfO \mfp_{G,g,r}$ is smooth over $k$ (cf. Theorem \ref{P0012},  (ii)),  so the cotangent complex $L_{\mfO \mfp^{^\mr{Zzz...}}_{G,g,r}}^\bullet$ is naturally isomorphic (in $D (\mcO_{\mfO  \mfp_{G, g,r}^{^\mr{Zzz...}}})$) to the complex $[\mcI_G/\mcI_G^2 \migi \iota_{g,r}^{*} (\Omega_{\mfO \mfp_{G, g,r}/k})]$.
Thus, 
the pair of $\underline{\xi}^{^\mr{Zz.}}$ and the identity morphism of   $\iota_{g,r}^{*} (\Omega_{\mfO \mfp_{G, g,r}/k})$
 specifies 
 a morphism
\begin{align}
\phi : E^\bullet \migi L_{\mfO \mfp^{^\mr{Zzz...}}_{G,g,r}}^\bullet 
\end{align} 
 in   $D (\mcO_{\mfO  \mfp_{G, g,r}^{^\mr{Zzz...}}})$.
 
\SSP
\bt \label{P01} 
Suppose that $G$ is of adjoint type.
Then, the morphism  $\phi$ just obtained   forms a perfect obstruction theory for $\mfO \mfp^{^\mr{Zzz...}}_{G,g,r}$ of constant virtual dimension $3g-3+r$ with $E^\bullet$ perfect.
  \et
\begin{proof}
It follows from Theorem \ref{P0012}, (ii),  that $\Omega_{\mfO \mfp_{G,g,r}/k}$ is locally free of rank $3g-3+r +N$, where $N := (g-1) \cdot \mr{dim} (G) + \frac{r}{2} \cdot (\mr{dim}(G) +\mr{rk} (G))$.
Hence,  the problem is reduced to proving 
 that  $\xi^{^\mr{Zz.}*}(\widehat{\mcI}_{G}/\widehat{\mcI}_{G}^2)$ is locally free of rank  $N$.

 For simplicity, let us denote by $\msX := (f : X \migi S, \{ \sigma_{X,i} \}_{i=1}^r)$ the tautological family of pointed stable curves over $S := \mfO \mfp_{G,g,r}^{^\mr{Zzz...}}$.
 Also, denote by $\msE^\spadesuit := (\mcE_B, \nabla)$ the tautological  dormant $G$-do'per on $\msX$.
 Let  $\nabla^\mr{ad} : \mfg_{\mcE_G} \migi \Omega_{X^\mr{log}/S^\mr{log}} \otimes \mfg_{\mcE_G}$ be  the $S^\mr{log}$-connection on the adjoint vector bundle $\mfg_{\mcE_G}$ induced by $\nabla$; we regard it as a complex  $\mcK^\bullet [\nabla^\mr{ad}]$ concentrated in degrees $0$ and $1$.
 Note that  there exists a  short  exact sequence of $\mcO_{S}$-modules
  \begin{align}
0 \migi  \mbR^1 f_*(\mr{Ker} (\nabla^\mr{ad})) \migi \mbR^1 f_* (\mcK^\bullet [\nabla^\mr{ad}]) \migi f_* (\mr{Coker} (\nabla^\mr{ad})) \migi 0 
\end{align}
arising from the conjugate spectral sequence of $\mcK^\bullet [\nabla^\mr{ad}]$ (cf. ~\cite[Eq.\,(758)]{Wak5}).
It follows from well-known generalities  of deformation theory 
 that
one may  construct a canonical isomorphism
\begin{align} \label{ee222}
(\iota^{\mfC \mfo}_{g,r} \circ \xi^{^\mr{Zz.}})^*(\mcT_{\mfC \mfo \mfn \mfn_{G, g,r}/\overline{\mfM}_{g,r}}) \isom \mbR^1 f_* (\mcK^\bullet [\nabla^\mr{ad}]).
\end{align}
This isomorphism  restricts to an isomorphism
 \begin{align} \label{ee223}
 \xi^{^\mr{Zz.}*} (\mcT_{\mfC \mfo \mfn \mfn^{^\mr{Zzz...}}_{G, g,r}/\overline{\mfM}_{g,r}}) \isom \mbR^1 f_* (\mr{Ker} (\nabla^\mr{ad}))
 \end{align}
(cf. ~\cite[Proposition 6.11]{Wak5}).
According to ~\cite[Proposition 6.18]{Wak5}, 
 $f_* (\mr{Coker} (\nabla^\mr{ad}))$ is locally free  of rank $N$.
Hence,  the duals of (\ref{ee222}) and (\ref{ee223}) induce  a commutative diagram
\begin{align} \label{ee225}
\begin{CD}
 @. 0
\\
@. @VVV
\\
\xi^{^\mr{Zz.}*}(\widehat{\mcI}_{G}/\widehat{\mcI}_{G}^2)@> >> f_* (\mr{Coker} (\nabla^\mr{ad}))^\vee
\\
@VVV @VVV
\\
(\iota^{\mfC \mfo}_{g,r} \circ \xi^{^\mr{Zz.}})^*(\Omega_{\mfC \mfo \mfn \mfn_{G, g,r}/\overline{\mfM}_{g,r}}) @> ((\ref{ee222})^\vee)^{-1} >> \mbR^1 f_* (\mcK^\bullet [\nabla^\mr{ad}])^\vee
\\
@VVV @VVV
\\
 \xi^{^\mr{Zz.}*} (\Omega_{\mfC \mfo \mfn \mfn^{^\mr{Zzz...}}_{G, g,r}/\overline{\mfM}_{g,r}}) @> ((\ref{ee223})^\vee)^{-1} >> \mbR^1 f_* (\mr{Ker} (\nabla^\mr{ad}))^\vee
\\
@VVV @VVV
\\
0 @. 0,
\end{CD}
\end{align}
where the both sides of vertical sequences are exact.
Since we have  the equalities $\mbR^2 f_* (\mcK^\bullet [\nabla^\mr{ad}])$ $=0$ (cf. ~\cite[Proposition 6.5, (iii)]{Wak5}) and $\mbR^2 f_* (\mr{Ker} (\nabla^\mr{ad})) =0$ (which follows from $\mr{dim} (X/S) =1$),
  the constructions of (\ref{ee222}) and (\ref{ee223}) implies that
both $\mfC \mfo \mfn \mfn_{G, g,r}$ and $\mfC \mfo \mfn \mfn^{^\mr{Zzz...}}_{G, g,r}$ are  smooth over $\overline{\mfM}_{g,r}$ at the points lying over $S$.
Hence,  the left-hand  top vertical arrow 
 in (\ref{ee225}) is injective.
It follows that the top horizontal arrow 
$\xi^{^\mr{Zz.}*}(\widehat{\mcI}_{G}/\widehat{\mcI}_{G}^2) \migi  f_* (\mr{Coker} (\nabla^\mr{ad}))^\vee$
is an isomorphism.
Consequently,  $\xi^{^\mr{Zz.}*}(\widehat{\mcI}_{G}/\widehat{\mcI}_{G}^2)$
is verified to be locally free of rank  $N$.
This completes the proof of the theorem.
\end{proof}
\SSP

By the above result, we obtain the virtual fundamental class
\begin{align}
[\mfO \mfp^{^\mr{Zzz...}}_{G_\mr{ad}, g,r}]^\mr{vir} :=  [\mfO \mfp^{^\mr{Zzz...}}_{G_\mr{ad}, g,r}, E^\bullet]^\mr{vir} \in A_{3g-3+r} (\mfO \mfp^{^\mr{Zzz...}}_{G_\mr{ad}, g,r})_\mbQ.
\end{align}
Moreover, since the morphism $\mr{op}_\mr{ad}^{^\mr{Zz.}}$ (cf. (\ref{ER44})) is finite and flat (cf. Theorem \ref{P05}, (i)),
the pull-back of this cycle class  by $\mr{op}_\mr{ad}^{^\mr{Zz.}}$ can be defined; we shall denote  it  by
 \begin{align} \label{ee226}
 [\mfO \mfp^{^\mr{Zzz...}}_{G, g,r}]^\mr{vir} :=  (\mr{op}_\mr{ad}^{^\mr{Zz.}})^*([\mfO \mfp^{^\mr{Zzz...}}_{G_\mr{ad}, g,r}]^\mr{vir}) \in A_{3g-3+r} (\mfO \mfp^{^\mr{Zzz...}}_{G, g,r})_\mbQ.
 \end{align}

\SSP
\bde \label{DD34}
We shall refer to $[\mfO \mfp^{^\mr{Zzz...}}_{G, g,r}]^\mr{vir}$ as the {\bf virtual fundamental class} on $\mfO \mfp^{^\mr{Zzz...}}_{G, g,r}$.
\ede

\SSP
\begin{rema} \label{R039}
In the case of   $G = \mr{PGL}_2$,
 the moduli stack $\mfO \mfp^{^\mr{Zzz...}}_{\mr{PGL}_2, g,r}$ is smooth over $k$ (cf. ~\cite[Chap.\,II, Theorem 2.8]{Mzk2}).
In particular,  the virtual fundamental class $[\mfO \mfp^{^\mr{Zzz...}}_{\mr{PGL}_2, g,r}]^\mr{vir}$ coincides  with the usual  fundamental class $[\mfO \mfp^{^\mr{Zzz...}}_{\mr{PGL}_2, g,r}]$.
 \end{rema}

\vspace{10mm}
\section{2d TQFTs for do'pers} \label{SSS666}
\SSP

In this section, we define  the notion of a 
{\it CohFT}
 (= a {\it cohomological field theory}) mapped into the $l$-adic \'{e}tale cohomology groups  of $\overline{\mfM}_{g, r}$'s.
After that,  a CohFT for $G$-do'pers is constructed 
by using the virtual fundamental class obtained in the previous section and the factorization properties  resulting from Propositions \ref{P0191}, \ref{P0192}, and  \ref{P0193}.
It also forms a {\it $2$d TQFT} (= a {\it $2$-dimensional topological quantum field theory}) such that the corresponding Frobenius algebra is semisimple.
 
 Let us keep the notation at the beginning of \S\,\ref{SSS444}.
 Moreover, we suppose that $k$ is algebraically closed.
 Also, let us fix a prime $l$ different from $p \left(= \mr{char}(k) \right)$.
 
\LSP
\subsection{$l$-adic  cohomology and Borel-Moore homology} \label{SS3}
Let   $\mfM$ be   a separated Deligne-Mumford stack of finite type  over $k$.
Denote by $D^b_c (\mfM, \overline{\mbQ}_l)$ the derived category of constructible $\overline{\mbQ}_l$-modules on $\mfM$. Hence, for each complex $\mcL$ in  $D^b_c (\mfM, \overline{\mbQ}_l)$, we obtain its  \'{e}tale cohomology $H_{\text{\'{e}t}}^i (\mfM, \mcL)$ ($i \in \mbZ$).
Also, we write
\begin{align}
 \widetilde{H}^i_{\text{\'{e}t}} (\mfM, \mcL) := H^i_{\text{\'{e}t}} (\mfM, \mcL(\lfloor {\textstyle \frac{i}{2}} \rfloor)),  \hspace{5mm} \widetilde{H}^*_{\text{\'{e}t}} (\mfM, \mcL) := \bigoplus_{i\in \mbZ} \widetilde{H}^{i}_{\text{\'{e}t}} (\mfM, \mcL).
\end{align}
For each $i \in \mbZ$, write  
$H^{\mr{BM}}_{i} (\mfM, \overline{\mbQ}_l)$
 for  the {\bf $i$-th Borel-Moore homology} of $\mfM$
   defined in ~\cite[Definition 2.2]{O2}.
That is to say, 
we set
\begin{align}
H^{\mr{BM}}_{i} (\mfM, \overline{\mbQ}_l) := H_{\text{\'{e}t}}^{-i} (\mfM, \omega_{\mfM}),
\end{align}
 where $\omega_\mfM \in D_c^b (\mfM, \overline{\mbQ}_l)$ denotes   the $l$-adic dualizing complex of $\mfM$.
Also, set  
\begin{align}
\widetilde{H}^\mr{BM}_i(\mfM, \overline{\mbQ}_l) := H^\mr{BM}_i (\mfM, \overline{\mbQ}_l) (- \lceil {\textstyle \frac{i}{2}}\rceil),  \hspace{5mm} \widetilde{H}^{\mr{BM}}_{*} (\mfM, \overline{\mbQ}_l) := \bigoplus_{i \in \mbZ} \widetilde{H}^{\mr{BM}}_{i} (\mfM, \overline{\mbQ}_l).
\end{align}

The {\it $i$-th cycle map} is the $\mbQ$-linear map
\begin{align} \label{ee405}
\mr{cl}^i  : A_i (\mfM)_\mbQ \migi \widetilde{H}^{\mr{BM}}_{2i} (\mfM, \overline{\mbQ}_l)
\end{align}
 mentioned in ~\cite[\S\,2.10]{O2}.
 If $\mfM$ is smooth of dimension $d$, then $\omega_\mfM \cong \overline{\mbQ}_l (d) [2d]$ and hence, we obtain a composite of natural isomorphisms
\begin{align} \label{ee370}
(-)^\blacklozenge : \widetilde{H}^{\mr{BM}}_{*} (\mfM, \overline{\mbQ}_l) & \isom 
\bigoplus_{i\geq 0}H_{\text{\'{e}t}}^{-i} (\mfM, \overline{\mbQ}_l (d)[2d]) ( - \lceil {\textstyle \frac{i}{2}} \rceil) \\
& \isom \bigoplus_{i\geq 0}H_{\text{\'{e}t}}^{2d-i} (\mfM, \overline{\mbQ}_l (\lfloor {\textstyle \frac{2d -i}{2}} \rfloor)) \notag \\
& \isom \widetilde{H}_{\text{\'{e}t}}^{2d-*} (\mfM, \overline{\mbQ}_l). \notag
\end{align}
Then, the cycle map $\mr{cl}^i$ coincides, via  (\ref{ee370}), with  the classical definition of the cycle map $A_i (\mfM)_\mbQ \migi \widetilde{H}_{\text{\'{e}t}}^{2d-2i} (\mfM, \overline{\mbQ}_l)$.

By the definition of $\widetilde{H}^\mr{BM}_* (-, \overline{\mbQ}_l)$,
the cup product in  $l$-adic  cohomology induces
a natural pairing 
\begin{align} \label{ee435}
(-) \cap (-) : \widetilde{H}^{2j}_{\text{\'{e}t}} (\mfM, \overline{\mbQ}_l) \otimes_{\overline{\mbQ}_l} \widetilde{H}^{\mr{BM}}_{2i} (\mfM, \overline{\mbQ}_l) \migi \widetilde{H}^{\mr{BM}}_{2(i-j)} (\mfM, \overline{\mbQ}_l).
\end{align}
If $\mfN$ is another  separated Deligne-Mumford stack of finite type over $k$ and $f : \mfM \migi \mfN$ is a proper morphism over $k$, 
then there exists  the  {\it pushforward map} 
\begin{align}
f^\mr{hom}_* : \widetilde{H}^{\mr{BM}}_{*} (\mfM, \overline{\mbQ}_l) \migi \widetilde{H}^{\mr{BM}}_{*} (\mfN, \overline{\mbQ}_l) 
\end{align}
along $f$ described in ~\cite[\S\,2.10]{O2}.
The following {\it projection formula} may be immediately verified:
\begin{align} \label{ee436}
f_*^\mr{hom}(f^*(\alpha) \cap \beta) = \alpha \cap f_*^\mr{hom} (\beta) \in \widetilde{H}^{\mr{BM}}_{2(i-j)} (\mfN, \overline{\mbQ}_l),
\end{align}
where $\alpha \in \widetilde{H}^{2j}_{\text{\'{e}t}}(\mfN, \overline{\mbQ}_l)$ and $\beta \in \widetilde{H}^\mr{BM}_{2i} (\mfM, \overline{\mbQ}_l)$ (cf. ~\cite[Proposition 5.2]{Li}).
If, moreover, $\mfN = \mr{Spec}(k)$, then by composing  $f^\mr{hom}_*$ and   (\ref{ee435}),
  we obtain,
for each class $\alpha \in A_i (\mfM)_\mbQ $, a $\overline{\mbQ}_l$-linear  morphism
\begin{align} \label{ee402}
\int_\alpha :  \widetilde{H}^{2i}_{\text{\'{e}t}} (\mfM, \overline{\mbQ}_l) \xrightarrow{(-)\cap \mr{cl}^i (\alpha)} \widetilde{H}^{\mr{BM}}_{0} (\mfM, \overline{\mbQ}_l) \xrightarrow{f^\mr{hom}_*} \widetilde{H}^{\mr{BM}}_{0} (\mr{Spec}(k), \overline{\mbQ}_l) \isom \overline{\mbQ}_l.
\end{align}
By assigning $v \mapsto 0$ for any $v \in \widetilde{H}^{j}_{\text{\'{e}t}} (\mfM, \overline{\mbQ}_l)$ with $j \neq 2i$,  we shall regard  
$\displaystyle \int_\alpha$ as a morphism  $\widetilde{H}^{*}_{\text{\'{e}t}} (\mfM, \overline{\mbQ}_l) \migi \overline{\mbQ}_l$.

\LSP
\subsection{$l$-adic CohFTs and  $2$d TQFTs} \label{SS3}
We shall describe  the definition of a  cohomological field theory
 by means of the $l$-adic \'{e}tale cohomologies $\widetilde{H}_{\text{\'{e}t}}^* (\overline{\mfM}_{g,r}, \overline{\mbQ}_l)$ of  the stacks $\overline{\mfM}_{g,r}$.
 Also, we  recall the notion of a $2$-dimensional topological quantum field theory as  a special kind of cohomological field theory.

\SSP
\bde \label{D07}
\begin{itemize}
\item[(i)]
 A(n) {\bf ($l$-adic) cohomological field theory (with flat identity)},  abbreviated  {\bf CohFT},  is a collection of data
\begin{align} \label{ee357}
\Lambda:= (\mcH, \eta, {\bf 1},  \{ \Lambda_{g,r} \}_{g, r \geq 0, 2g-2+r >0})
\end{align}
consisting of 
\begin{itemize}
\item
  a finite dimensional  $\overline{\mbQ}_l$-vector space $\mcH$  (called the {\bf state space}) with basis $\mfe := \{ e_1, \cdots, e_{\mr{dim} (\mcH)} \}$;
\item
 a symmetric nondegenerate  pairing 
$\eta : \mcH \times \mcH \migi \overline{\mbQ}_l$  (called the {\bf metric}), where we shall write $(\eta^{e_a e_b})_{a, b}$ for the inverse of the matrix corresponding to  the metric $\eta$ with respect to the basis $\mfe$;
\item
 an element ${\bf 1}$ of $\mcH$;
\item
 $\overline{\mbQ}_l$-linear morphisms $\Lambda_{g,r} : \mcH^{\otimes r} \migi \widetilde{H}^*_{\text{\'{e}t}} (\overline{\mfM}_{g,r}, \overline{\mbQ}_l)$ (called the {\bf correlators}), where $\mcH^{\otimes 0} := \overline{\mbQ}_l$,
\end{itemize}
and satisfying the following conditions:
\begin{itemize}
\item[$\bullet$]
Each $\Lambda_{g,r}$ is 
compatible with 
 the respective actions of the  symmetric group $\mfS_r$ on $\mcH^{\otimes r}$ and $\widetilde{H}^*_{\text{\'{e}t}} (\overline{\mfM}_{g,r}, \overline{\mbQ}_l)$ arising from permutations of the $r$ factors in $\mcH^{\otimes r}$ and the $r$ punctures  in the tautological family of curves over $\overline{\mfM}_{g,r}$ respectively.
\item[$\bullet$]
For any $v_1$, $v_2 \in \mcH$, the following equality holds: 
\begin{align} \label{ee356}
\eta (v_1,   v_2) = \int_{\overline{\mfM}_{0,3}} \Lambda_{0,3} (v_1 \otimes v_2 \otimes {\bf 1}).
\end{align}
\item[$\bullet$]
For any $v_1, \cdots, v_{r_1 + r_2} \in \mcH$, the the following equality holds: 
\begin{align} \label{ee355}
&  \Phi_{\mr{tree}}^* (\Lambda_{g_1 + g_2, r_1 + r_2} (v_1 \otimes \cdots \otimes v_{r_1 + r_2})) \\
 = &  \sum_{e_1, e_2 \in \mfe}  \Lambda_{g_1, r_1 +1} (v_1 \otimes \cdots \otimes v_{r_1} \otimes e_1) \eta^{e_1 e_2}\otimes \Lambda_{g_2, r_2 +1} (e_2 \otimes v_{r_1 +1} \otimes \cdots \otimes v_{r_1 +r_2}), \notag 
\end{align}
where $\Phi_{\mr{tree}}^*$ denotes the morphism
\begin{align} \label{ee354}
\widetilde{H}^*_{\text{\'{e}t}} (\overline{\mfM}_{g_1 +g_2,r_1+r_2}, \overline{\mbQ}_l) \migi \widetilde{H}^*_{\text{\'{e}t}} (\overline{\mfM}_{g_1,r_1+1}, \overline{\mbQ}_l) \otimes_{\overline{\mbQ}_l} \widetilde{H}^*_{\text{\'{e}t}} (\overline{\mfM}_{g_2,r_2+1}, \overline{\mbQ}_l) 
\end{align}
induced by $\Phi_{\mr{tree}}$ (cf. (\ref{e072})).
\item[$\bullet$]
For any $v_1, \cdots, v_r \in \mcH$, the following equality holds:
\begin{align} \label{ee353}
\Phi_\mr{loop}^* (\Lambda_{g+1,r} (v_1\otimes  \cdots \otimes  v_r)) = \sum_{e_1, e_2 \in \mfe}\Lambda_{g, r+1} (v_1 \otimes  \cdots  \otimes v_r \otimes e_1 \otimes e_2) \eta^{e_1 e_2},
\end{align}
where $\Phi_\mr{loop}^*$ denotes the morphism $\widetilde{H}^*_{\text{\'{e}t}} (\overline{\mfM}_{g +1,r}, \overline{\mbQ}_l) \migi \widetilde{H}^*_{\text{\'{e}t}} (\overline{\mfM}_{g,r+2}, \overline{\mbQ}_l)$
induced by $\Phi_\mr{loop}$ (cf. (\ref{e071})).
\item[$\bullet$]
For any $v_1, \cdots, v_r \in \mcH$, the following equality holds:
\begin{align} \label{ee351}
\Phi_\mr{tail}^* (\Lambda_{g,r} (v_1 \otimes  \cdots \otimes v_r)) = \Lambda_{g, r+1} (v_1 \otimes \cdots \otimes  v_r \otimes {\bf 1}),
\end{align}
where $\Phi_\mr{tail}^*$ denotes the morphism $\widetilde{H}^*_{\text{\'{e}t}} (\overline{\mfM}_{g,r}, \overline{\mbQ}_l) \migi \widetilde{H}^*_{\text{\'{e}t}} (\overline{\mfM}_{g,r+1}, \overline{\mbQ}_l)$
induced by $\Phi_\mr{tail}$ (cf. (\ref{e070})).
\end{itemize}
\item[(ii)]
A {\bf $2$-dimensional topological quantum field theory (over $\overline{\mbQ}_l$)},  abbreviated  {\bf $2$d TQFT}, 
is an $l$-adic cohomological field theory whose correlators $\Lambda_{g,r}$ are all valued in $\widetilde{H}^0_{\text{\'{e}t}} (\overline{\mfM}_{g,r}, \overline{\mbQ}_l) = \overline{\mbQ}_l$ (cf. ~\cite[\S\,1.3.32]{Kock} for the naive definition of an $n$-dimensional topological quantum field theory). 
\end{itemize}
  \ede
\SSP

It  is well-known that   $2$d TQFTs correspond to Frobenius algebras.
Here, 
by a {\bf Frobenius algebra}  over $\overline{\mbQ}_l$, we mean (cf. ~\cite[\S\,2.2.5]{Kock})  a pair $(\mcH, \eta)$ consisting of a unital, associative, and commutative $\overline{\mbQ}_l$-algebra $\mcH$ of finite dimension and  a nondegenerate $\overline{\mbQ}_l$-bilinear pairing $\eta : \mcH \times \mcH \migi \overline{\mbQ}_l$ such that 
\begin{align}
\eta (v_1, (v_2 \times v_3)) = \eta ((v_1 \times v_2),  v_3)
\end{align}
 for any $v_1$, $v_2$, $v_3 \in \mcH$, where $\times$ denotes the multiplication in $\mcH$.
We shall say that a Frobenius algebra $(\mcH, \eta)$ is {\bf semisimple} if there exists a basis $\mfe^\dagger := \{ e^\dagger_a\}_{a \in I}$
  of $\mcH$ such that 
\begin{align} \label{ee430}
e^\dagger_{a} \times e^\dagger_{b} = \delta_{ab}  e^\dagger_{a} \hspace{5mm} \text{and}\hspace{5mm} \eta (e^\dagger_{a},  e^\dagger_{b}) = \delta_{ab}  \nu_a
\end{align}
 for any $a$, $b \in I$, where  each $\nu_{a}$ is  some {\it nonzero} element of $\overline{\mbQ}_l$.
 We shall refer to $\mfe^\dagger$ as  a {\bf canonical   basis} of $(\mcH, \eta)$.
 
 Now, let  
  $\Lambda := (\mcH, \eta, {\bf 1}, \{ \Lambda_{g,r} \}_{g,r})$ be a $2$d TQFT over $\overline{\mbQ}_l$ (where let $\mfe := \{ e_a  \}_{a \in I}$
  denote  the distinguished  basis of $\mcH$).
 There exists 
 a multiplication $\times : \mcH \times  \mcH \migi \mcH$ on $\mcH$ given by
 \begin{align}
 u \times v := \sum_{a,b \in I} \Lambda_{0, 3} (u \otimes  v \otimes  e_a) \eta^{ab} e_b. 
 \end{align}
 The $\overline{\mbQ}_l$-vector space  $\mcH$ together with  this multiplication forms  a unital, associative, and commutative $\overline{\mbQ}_l$-algebra, in which 
  ${\bf 1}$ is    the  unit.
  Moreover,  the pair $(\mcH, \eta)$ forms a Frobenius algebra over $\overline{\mbQ}_l$, which we shall refer to as the {\bf Frobenius algebra associated with $\Lambda$}.

\LSP
\subsection{The $2$d TQFT  associated with the moduli space of  do'pers} \label{SS3}
As described in Introduction, we set 
\begin{align} \label{ee140}
\mcV := \widetilde{H}^*_{\text{\'{e}t}} (\mfR \mfa \mfd, \overline{\mbQ}_l). 
\end{align}
Decomposition (\ref{ee141}) gives rise to 
a composite isomorphism
\begin{align} \label{EE888}
\mcV & \isom \widetilde{H}^*_{\text{\'{e}t}} (\coprod_{\pmb{\rho} \in \Delta} 
\mfR \mfa \mfd_{\pmb{\rho}}, \overline{\mbQ}_l)  
 \isom \bigoplus_{\pmb{\rho} \in \Delta} \widetilde{H}^*_{\text{\'{e}t}} (
 \mfR \mfa \mfd_{\pmb{\rho}}, \overline{\mbQ}_l).  
\end{align}
For each $\pmb{\rho} \in \Delta$,  denote by 
$e_{\pmb{\rho}}$
the element of 
$\mcV$ corresponding, via  (\ref{EE888}),
to $1 \in \widetilde{H}^0_{\text{\'{e}t}} (\mfR \mfa \mfd_{\pmb{\rho}}, \overline{\mbQ}_l)$
$\left(\subseteq \widetilde{H}^*_{\text{\'{e}t}} (\mfR \mfa \mfd_{\pmb{\rho}}, \overline{\mbQ}_l)\right)$.
In particular, we obtain 
\begin{align} \label{fff444}
e_{\pmb{\varepsilon}} \in \mcV.
\end{align}

Next, 
let $(g, r)$ be a pair of nonnegative integers with $2g-2+r >0$.
We  define a $\overline{\mbQ}_l$-linear morphism
\begin{align} \label{ee501}
\Lambda_{G, g,r} : \mcV^{\otimes r} \migi  \widetilde{H}^*_{\text{\'{e}t}} (\overline{\mfM}_{g,r}, \overline{\mbQ}_l)
\end{align}
to be the morphism determined uniquely by 
\begin{align} \label{ee502}
\Lambda_{G, g,r}(\bigotimes_{i=1}^r v_i) := \left(\pi_{g,r *}^\mr{hom}\left(
\left(\prod_{i=1}^r \mr{ev}_i^{*}(v_i)\right) \cap  \mr{cl}^{3g-3+r}([\mfO \mfp^{^\mr{Zzz...}}_{G, g,r}]^{\mr{vir}})  \right)\right)^\blacklozenge
\end{align}
(cf. (\ref{ee216}) and (\ref{ER48}) for the definitions of $\pi_{g,r}$ and  $\mr{ev}_i$ respectively) for every  $v_1, \cdots, v_r \in \mcV$.

Finally, let us  write
\begin{align} \label{ee80255}
\eta : \mcV \times  \mcV \migi \overline{\mbQ}_l
\end{align}
for the $\overline{\mbQ}_l$-bilinear pairing determined by
\begin{align} \label{ee802}
\eta (v_1, v_2) 
= \left(\pi_{0,3 *}^\mr{hom}\left(
\left(  \mr{ev}_1^{*}(v_1) \mr{ev}_2^{*}(v_2) \mr{ev}_3^{*}(e_{\pmb{\varepsilon}}) \right) \cap  \mr{cl}^{0}([\mfO \mfp^{^\mr{Zzz...}}_{G, 0,3}]^{\mr{vir}})  \right)\right)^\blacklozenge.
\end{align}
Thus, we have obtained  a collection of data
\begin{align} \label{ee148}
\Lambda_{G} :=(\mcV, \eta, e_{\pmb{\varepsilon}}, \{\Lambda_{G, g,r} \}_{g,r \geq 0, 2g-2+r >0}).
\end{align}

\SSP
\bt \label{T10} Suppose that $G$ satisfies the condition  $(**)_{G}$.
\begin{itemize}
\item[(i)]
The set of elements $\mfe_\Delta := \{ e_{\pmb{\rho}}\}_{\pmb{\rho} \in \Delta}$ forms a basis of $\mcV$.
Also, for each $\vec{\pmb{\rho}} := (\pmb{\rho}_i)_{i=1}^r \in \Delta^{\times r}$, 
$\Lambda_{G, g,r} (\bigotimes_{i=1}^r e_{\pmb{\rho}_i})$ lies in 
$\widetilde{H}^0_{\text{\'{e}t}} (\overline{\mfM}_{g,r}, \overline{\mbQ}_l)$
  and 
the following equality holds:
\begin{align}
\Lambda_{G, g,r} (\bigotimes_{i=1}^r e_{\pmb{\rho}_i}) = \mr{deg}^\mr{gen} (\mfO \mfp^{^\mr{Zzz...}}_{G,g,r, \vec{\pmb{\rho}}}/ \overline{\mfM}_{g,r})
\end{align}
(cf. (\ref{Ww600})).
 Finally,  
 the pairing $\eta$ is symmetric, nondegenerate, and moreover, given by
 \begin{align}
 \eta_{\pmb{\rho}_1 \pmb{\rho}_2} := \eta (e_{\pmb{\rho}_1},  e_{\pmb{\rho}_2}) = \begin{cases} \frac{1}{| Z |}
  & \text{if $\pmb{\rho}_1 = \pmb{\rho}_2^\veebar$}; \\  \ 0 & \text{if otherwise.}\end{cases}
 \end{align}

\item[(ii)]
The collection of data
$\Lambda_{G}$
 forms a CohFT  valued in $\widetilde{H}^0_{\text{\'{e}t}} (\overline{\mfM}_{g,r}, \overline{\mbQ}_l)$,
  namely, forms a $2$d TQFT over  $\overline{\mbQ}_l$.
 Moreover, the corresponding  Frobenius algebra $(\mcV, \eta)$ is semisimple.
\end{itemize}
  \et
\begin{proof}
Let us consider  assertion (i).
The first assertion follows from 
the following sequence of isomorphisms:
\begin{align} \label{ee402}
\widetilde{H}^*_{\text{\'{e}t}} (
\mfR \mfa \mfd_{\pmb{\rho}}, \overline{\mbQ}_l) \isom \widetilde{H}^*_{\text{\'{e}t}} (\mcB \mr{Coker}(\kappa), \overline{\mbQ}_l) \isom \widetilde{H}^0_{\text{\'{e}t}} (\mcB \mr{Coker}(\kappa), \overline{\mbQ}_l) =  \overline{\mbQ}_l,
\end{align}
where $\pmb{\rho} := (\rho, \kappa) \in \Delta$ and the first isomorphism follows from Proposition \ref{Lw001}.

Next, let us consider the second assertion.
The definition of $\Lambda_{G, g,r}$ implies 
the equality
\begin{align} \label{ee750}
\Lambda_{G, g,r} (\bigotimes_{i=1}^r e_{\pmb{\rho}_i}) = \left( \pi^\mr{hom}_{g,r*} \left( \mr{cl}^{3g-3+r}  \left([\mfO \mfp^{^\mr{Zzz...}}_{G,g,r}]^\mr{vir} \Big|_{\mfO \mfp^{^\mr{Zzz...}}_{G,g,r, \vec{\pmb{\rho}}}}\right)\right) \right)^\blacklozenge,
\end{align}
where $(-) |_{\mfO \mfp^{^\mr{Zzz...}}_{G,g,r, \vec{\pmb{\rho}}}}$ denotes the restriction of the class $(-)$ to the component $\mfO \mfp^{^\mr{Zzz...}}_{G,g,r, \vec{\pmb{\rho}}}$.
Also, since the square diagram
\begin{align}
\begin{CD}
A_{3g-3+r} (\mfO \mfp^{^\mr{Zzz...}}_{G, g,r})_\mbQ @> \mr{cl}^{3g-3+r} >> \widetilde{H}^\mr{BM}_{6g-6+2r} (\mfO \mfp^{^\mr{Zzz...}}_{G, g,r}, \overline{\mbQ}_l)
\\
@V \pi_{g,r *} VV @VV \pi_{g,r *}^{\mr{hom}} V
\\
A_{3g-3+r} (\overline{\mfM}_{g,r})_\mbQ @>> \mr{cl}^{3g-3+r} > \widetilde{H}^\mr{BM}_{6g-6+2r} (\overline{\mfM}_{g,r}, \overline{\mbQ}_l)
\end{CD}
\end{align}
is commutative (cf. ~\cite[\S\,2.10]{O2}), we have
\begin{align} \label{ee411}
& \ \pi^\mr{hom}_{g,r*} \left( \mr{cl}^{3g-3+r} \left([\mfO \mfp^{^\mr{Zzz...}}_{G,g,r}]^\mr{vir} \Big|_{\mfO \mfp^{^\mr{Zzz...}}_{G,g,r, \vec{\pmb{\rho}}}} \right)\right) = 
\mr{cl}^{3g-3+r} \left( \pi_{g,r*}\left([\mfO \mfp^{^\mr{Zzz...}}_{G,g,r}]^\mr{vir} \Big|_{\mfO \mfp^{^\mr{Zzz...}}_{G,g,r, \vec{\pmb{\rho}}}}\right)\right). 
\end{align}
Since
    $\mfO \mfp^{^\mr{Zzz...}}_{G,g,r, \vec{\pmb{\rho}}}$ is finite over $\overline{\mfM}_{g,r}$ (cf. Theorem \ref{P05}, (ii))  and $\overline{\mfM}_{g,r}$ is an irreducible stack of  dimension $3g-3+r$,
 any prime cycle in 
$A_{3g-3+r} (\mfO \mfp^{^\mr{Zzz...}}_{G,g,r, \vec{\pmb{\rho}}})_\mbQ$
 dominates $\overline{\mfM}_{g,r}$.
On the other hand, 
 the generic \'{e}taleness of $\mfO \mfp^{^\mr{Zzz...}}_{G_\mr{ad}, g,r, (\rho_i)_{i=1}^r}$ (where we set  $((\rho_i, \kappa_i))_{i=1}^r := \vec{\pmb{\rho}}$)
 over $\overline{\mfM}_{g,r}$ (cf. ~\cite[Theorem G]{Wak5}) implies the generic \'{e}taleness of 
 $\mfO \mfp^{^\mr{Zzz...}}_{G, g,r, \vec{\pmb{\rho}}}$ ($\cong \mfO \mfp^{^\mr{Zzz...}}_{G_\mr{ad}, g,r, (\rho_i)_{i}} \times_{\overline{\mfM}_{g,r}} \mfS \mfp_{Z, \delta^\sharp, g,r, (\kappa_i)_{i}}$ by Theorem \ref{P05}, (i)) over $\mfS \mfp_{Z, \delta^\sharp, g,r, (\kappa_i)_{i}}$.
 Hence, the smoothness of $\mfS \mfp_{Z, \delta^\sharp, g,r, (\kappa_i)_{i}}$  (cf. Theorem \ref{P08f}) implies that 
 $\mfO \mfp^{^\mr{Zzz...}}_{G, g,r, \vec{\pmb{\rho}}}$ is    generically smooth over $k$.
To be precise, 
any irreducible component $\mfN$ of 
$\mfO \mfp^{^\mr{Zzz...}}_{G, g,r, \vec{\pmb{\rho}}}$
   dominating $\overline{\mfM}_{g,r}$ has a 
     dense open substack $\mfN^\circ$ which (does not intersects any other irreducible components and) is smooth over $k$.
The restriction of $[\mfO \mfp^{^\mr{Zzz...}}_{G,g,r}]^\mr{vir} \Big|_{\mfO \mfp^{^\mr{Zzz...}}_{G,g,r, \vec{\pmb{\rho}}}}$ to  $\mfN^\circ$ coincides with the usual fundamental class  $[\mfN^\circ]$ in the usual sense.
By
the observations made so far and  
  the definition of the pushforward map $\pi_{g,r*}$ between rational Chow groups,
  the following equality turns out to be satisfied: 
 \begin{align} \label{ee410}
 \pi_{g,r*}\left([\mfO \mfp^{^\mr{Zzz...}}_{G,g,r}]^\mr{vir} \Big|_{\mfO \mfp^{^\mr{Zzz...}}_{G,g,r, \vec{\pmb{\rho}}}}\right) 
 =
  \mr{deg}^{\mr{gen}} (\mfO \mfp^{^\mr{Zzz...}}_{G,g,r, \vec{\pmb{\rho}}}/\overline{\mfM}_{g,r}) \cdot [\overline{\mfM}_{g,r}].
 \end{align}
 Thus,  (\ref{ee750}),  (\ref{ee411}),  and (\ref{ee410}) give  the following sequence of equalities:
 \begin{align}
 \Lambda_{G, g,r} (\bigotimes_{i=1}^r e_{\pmb{\rho}_i}) 
 &  \stackrel{(\ref{ee750})}{=}  \left(\pi_{g,r *}^\mr{hom}\left( \mr{cl}^{3g-3+r} \left([\mfO \mfp^{^\mr{Zzz...}}_{G,g,r}]^\mr{vir} \Big|_{\mfO \mfp^{^\mr{Zzz...}}_{G,g,r, \vec{\pmb{\rho}}}} \right) \right)\right)^\blacklozenge \\
 & \stackrel{(\ref{ee411})}{=} \left(\mr{cl}^{3g-3+r} \left( \pi_{g,r*}\left([\mfO \mfp^{^\mr{Zzz...}}_{G,g,r}]^\mr{vir} \Big|_{\mfO \mfp^{^\mr{Zzz...}}_{G,g,r, \vec{\pmb{\rho}}}}\right)\right) \right)^\blacklozenge  \notag  \\
 & \stackrel{(\ref{ee410})}{=}  \left(\mr{cl}^{3g-3+r} \left(\mr{deg} (\mfO \mfp^{^\mr{Zzz...}}_{G,g,r, \vec{\pmb{\rho}}}/\overline{\mfM}_{g,r}) \cdot [\overline{\mfM}_{g,r}]\right) \right)^\blacklozenge \notag \\
 &  \  = \mr{deg}^\mr{gen} (\mfO \mfp^{^\mr{Zzz...}}_{G,g,r, \vec{\pmb{\rho}}}/\overline{\mfM}_{g,r})   \cdot \mr{cl}^{3g-3+r} \left([\overline{\mfM}_{g,r}]\right)^\blacklozenge \notag \\
 & \  = \mr{deg}^\mr{gen} (\mfO \mfp^{^\mr{Zzz...}}_{G,g,r, \vec{\pmb{\rho}}}/\overline{\mfM}_{g,r}). \notag
 \end{align}
This completes the proof of the second assertion.
The third assertion follows from Theorem \ref{P05}, (iv), and hence,  the proof of assertion (i) is completed.

Finally, let us consider assertion (ii).
The former assertion  follows from   assertion (i) and Propositions \ref{P0191}, (ii), \ref{P0192}, (ii), and \ref{P0193}, (ii).
In what follows, we shall consider the latter assertion, i.e.,
 the semisimplicity  of the Frobenius algebra $(\mcV, \eta)$.
 To this end, it suffices to prove that the $\overline{\mbQ}_l$-algebra $\mcV$ is reduced.
 Let us fix an isomorphism $\overline{\mbQ}_l \isom \mbC$.
Denote by $\gamma$ the involution on $\mcV$ (viewed as an $\mbR$-algebra) given by $\sum_{\pmb{\rho} \in \Delta} v_{\pmb{\rho}} e_{\pmb{\rho}} \mapsto \sum_{\pmb{\rho} \in \Delta} \overline{v}_{\pmb{\rho}} e_{\pmb{\rho}^\veebar}$, 
where 
each $v_{\pmb{\rho}}$
 is an element of $\mbC \ (\cong \overline{\mbQ}_l)$
and  $\overline{(-)}$ denotes the complex conjugation.
Note that 
 the equality  $\eta (x,  \gamma (x)) = \frac{1}{|Z|}\sum_{\pmb{\rho} \in \Delta} | v_{\pmb{\rho}}|^2$ holds for each 
 $x := \sum_{\pmb{\rho} \in \Delta} v_{\pmb{\rho}} e_{\pmb{\rho}} \in \mcV$.
 Hence, $\eta (x,  \gamma (x)) = 0$ implies $x = 0$.
 Now, let us take 
  an element $x :=  \sum_{\pmb{\rho} \in \Delta} v_{\pmb{\rho}} e_{\pmb{\rho}} \in \mcV$ with $x \times x =0$.
Then,  
\begin{align}
\eta (x \times \gamma (x),  \gamma (x \times \gamma (x)))  & = \eta (x \times \gamma (x), \gamma (x) \times x) \\
 & = \eta (x \times x,  \gamma (x)\times \gamma (x)) \notag \\
 & = 0, \notag 
\end{align}
 which implies $x \times \gamma(x) =0$.
It follows that $\eta (x,   \gamma (x)) =0$, and hence, that $x =0$.
Consequently, $\mcV$ turns out to  be reduced.
 This completes the proof of the latter assertion of (ii).
\end{proof}
\SSP

\begin{rema} \label{R03jj6} 
\begin{itemize}
\item[(i)]
Since  the Frobenius algebra $(\mcV, \eta)$ associated with
$\Lambda_G$ is semisimple, there exists a canonical 
basis $\{ e^\dagger_{\pmb{\rho}}\}_{\pmb{\rho} \in \Delta}$
  of $\mcV$.
Let us write $e^\ddagger_{\pmb{\rho}} := (\nu_{\pmb{\rho}}^{1/2})^{-1} e^\dagger_{\pmb{\rho}}$, where each $\nu_{\pmb{\rho}}^{1/2}$ is an element of $\overline{\mbQ}_l$ with  $(\nu_{\pmb{\rho}}^{1/2})^2 =  \nu_{\pmb{\rho}}$.
The {\bf $S$-matrix} $S := (S_{\pmb{\rho} \pmb{\lambda}})_{\pmb{\rho} \pmb{\lambda}} \in \mr{GL} (\mcV)$
 is defined in such a way that 
$e_{\pmb{\rho}} = \sum_{\pmb{\lambda} \in \Delta} S_{\pmb{\rho} \pmb{\lambda}} e^\ddagger_{\pmb{\lambda}}$ (i.e., $S_{\pmb{\rho} \pmb{\lambda}} := \underline{S}^{-1}_{\pmb{\rho} \pmb{\lambda}}$ in the sense of ~\cite[\S\,3.3.5, Eq.\,(20)]{ABO}) for any $\pmb{\rho} \in \Delta$.
Hence, $\nu_{\pmb{\rho}} = S_{\pmb{\varepsilon} \pmb{\rho}}^2$ ($\pmb{\rho} \in \Delta$) (cf. ~\cite[\S\,3.3.8]{ABO}).
By applying the discussion in ~\cite[\S\,4.3]{ABO}, we obtain the {\it Verlinde formula for $G$-do'pers}.
That is to say, for each $(\pmb{\rho}_i)_{i=1}^r \in \Delta^{\times r}$,  the following equality holds: 
\begin{align} \label{ee700}
\Lambda_{G, g,r} (\bigotimes_{i=1}^r e_{\pmb{\rho}_i})  \left(=\mr{deg}^\mr{gen} (\mfO \mfp^{^\mr{Zzz...}}_{G, g,r, (\pmb{\rho}_i)_{i=1}^r}/\overline{\mfM}_{g,r}) \right)  = \sum_{\pmb{\lambda} \in \Delta} \frac{\prod_{i=1}^r S_{\pmb{\rho}_i \pmb{\lambda}}}{S^{2g-2+r}_{\pmb{\varepsilon \lambda}}}.
\end{align}
This formula contains the case of $g =0$, i.e.,
\begin{align}
\Lambda_{G, g,r} (1)  \left(=\mr{deg}^\mr{gen} (\mfO \mfp^{^\mr{Zzz...}}_{G, g,0}/\overline{\mfM}_{g,0}) \right)  = \sum_{\pmb{\rho} \in \Delta}  \nu^{1-g}_{\pmb{\rho}} = \sum_{\pmb{\rho} \in \Delta} S_{\pmb{\varepsilon} \pmb{\rho}}^{2-2g}.
\end{align}
\item[(ii)]
Assume further that $G$ is of adjoint type, i.e., $|Z| = 1$.
Then, 
$\Lambda_G$ defines  (under the fixed isomorphism $\overline{\mbQ}_l \isom \mbC$) a {\it fusion ring} in the sense of ~\cite[Definition 11.10]{Sch}.
Formula (\ref{ee700}) may be expressed  as 
\begin{align} \label{ee771}
\Lambda_{G, g,r} (\bigotimes_{i=1}^r e_{\pmb{\rho}_i}) 
 = \sum_{\chi \in \mfS} \chi (\mcC as_G)^{g-1} \prod_{i=1}^r \chi (e_{\pmb{\rho}_i})
 \end{align}
 (cf. ~\cite[Theorem F]{Wak5}), where $\mcC as_G := \sum_{\pmb{\rho} \in \Delta} \pmb{\rho} \times \pmb{\rho}^\veebar$ and $\mfS$ denotes the set of  characters (i.e.,  morphisms of $\overline{\mbQ}_l$-algebras) $\mcV \migi \overline{\mbQ}_l$.
Recall that the forgetting morphism  $\mfO \mfp^{^\mr{Zzz...}}_{G, g,r, (\pmb{\rho}_i)_{i=1}^r} \migi \overline{\mfM}_{g,r}$  is
 representable,  finite,  and generically \'{e}tale.
Hence, for a sufficiently general curve $\msX$ in $\overline{\mfM}_{g,r}$,
the number of $G$-do'pers of radii $(\pmb{\rho}_i )_{i=1}^r$ on $\msX$
 is exactly equal to the value  
 $\Lambda_{G, g,r} (\bigotimes_{i=1}^r e_{\pmb{\rho}_i})$.

If, moreover,  $G = \mr{PGL}_n$  for a small $n$ (relative to $g$ and $p$), then the result of    ~\cite[Theorem H]{Wak5}  allows  us to 
   compute the value $\Lambda_{\mr{PGL}_n, g, 0}(1)$ without  explicit knowledge  of the characters $\mcV \migisurj \overline{\mbQ}_l$.
Indeed, 
under the assumption that
 $p>n \cdot \mr{max} \{g-1, 2 \}$,
 the following formula holds:
 \begin{align} \Lambda_{\mr{PGL}_n, g, 0}(1) & \left(=\mr{deg}^\mr{gen}(\mfO \mfp_{\mr{PGL}_n,  g,0}^{^\mr{Zzz...}}/\overline{\mfM}_{g, 0})  \right) \\
 & =    \frac{ p^{(n-1)(g-1)-1}}{n!} \cdot  
 \sum_{\genfrac{.}{.}{0pt}{}{(\zeta_1, \cdots, \zeta_n) \in \mbC^{\times n} }{ \zeta_i^p=1, \ \zeta_i \neq \zeta_j (i\neq j)}}
\frac{(\prod_{i=1}^n\zeta_i)^{(n-1)(g-1)}}{\prod_{i\neq j}(\zeta_i -\zeta_j)^{g-1}}.  \notag
 \end{align}
\end{itemize}
 \end{rema}

\vspace{10mm}
\section{The  Witten-Kontsevich theorem for do'pers} \label{SSS6} \SSP
In this final section, we introduce the {\it partition function for
do'pers}
  and 
apply the discussion in ~\cite[\S\,4]{JK}  in order to obtain a result analogous   to the Witten-Kontsevich theorem; this result   gives nontrivial relationships among the  intersection numbers of the psi classes on $\mfO \mfp^{^\mr{Zzz...}}_{G, g,r}$ (cf. Theorem \ref{TT01}).

Let us keep the notation at the beginning of \S\,\ref{SSS444}. 
Also, suppose that $G$ satisfies 
the condition $(**)_{G}$.

\LSP
\subsection{Correlator functions} \label{SS5e}
Denote by $\psi_i \in \widetilde{H}_{\text{\'{e}t}}^2 (\overline{\mfM}_{g,r}, \overline{\mbQ}_l)$ ($i= 1, \cdots, r$)  the $i$-th psi class on $\overline{\mfM}_{g,r}$.
Given a pair of nonnegative integers $(g, r)$ and  an $r$-tuple of nonnegative integers  $d_1, \cdots, d_r$,
we shall recall (cf., e.g., ~\cite[Introduction]{Mirz}) the invariants
\begin{align} \label{ee150}
\langle \tau_{d_1} \cdots \tau_{d_r} \rangle_{g} \left(= \langle \prod_{i=1}^r \tau_{d_i} \rangle_g \right) := \int_{\overline{\mfM}_{g,r}} \prod_{i=1}^r \psi_i^{d_i}
\in \overline{\mbQ}_l,
\end{align}
where $\langle - \rangle_g := 0$  if either $r =0$ or $2g-2+r \leq 0$ is satisfied.

\SSP
\begin{rema} \label{REm2}
In most cases, the intersection theory of  psi classes is discussed in terms of 
the orbifold 
$\overline{\mfM}_{g,r, \mbC}^\mr{top}$
(or the corresponding topological stack in the sense of ~\cite{Noo}) associated with $\overline{\mfM}_{g,r, \mbC} := \overline{\mfM}_{g, r, \mbZ} \times_\mbZ \mbC$ and  
 the usual complex cohomology $H_{}^* (\overline{\mfM}_{g,r, \mbC}^\mr{top}, \mbC)$  of $\overline{\mfM}_{g,r, \mbC}^\mr{top}$.
However, after fixing  an isomorphism $\overline{\mbQ}_l \isom \mbC$,
the argument in ~\cite[Chap.\,III, Theorem 3.12]{Mil} together with  Riemann's existence theorem for  stacks (cf. ~\cite[Theorem 20.4]{Noo}) yields  an isomorphism
$H^*_{\text{\'{e}t}} (\overline{\mfM}_{g, r, \mbC}, \overline{\mbQ}_l) \isom H_{}^* (\overline{\mfM}_{g,r, \mbC}^\mr{top}, \mbC)$ preserving the cup product and  the Chern class maps.
On the other hand,
if we set  $\overline{\mfM}_{g, r, \mbZ_p} := \overline{\mfM}_{g, r} \times_\mbZ \mbZ_p$, then  the natural morphism $\mbZ_p \migi k$ and a fixed inclusion  $\mbZ_p \migiincl \mbC$ induce
the base change maps $H_{\text{\'{e}t}}^* (\overline{\mfM}_{g,r, \mbZ_p}, \overline{\mbQ}_l) \migi H_{\text{\'{e}t}}^* (\overline{\mfM}_{g,r}, \overline{\mbQ}_l)$ and
 $H_{\text{\'{e}t}}^* (\overline{\mfM}_{g,r, \mbZ_p}, \overline{\mbQ}_l) \migi H_{\text{\'{e}t}}^* (\overline{\mfM}_{g,r, \mbC}, \overline{\mbQ}_l)$ respectively.
These maps  preserve the cup product  and moreover preserve the $1$-st Chern class map $c_1 (-)$ because of the construction using the Kummer sequence.
Thus, we conclude that
the invariants $\langle \tau_{d_1} \cdots \tau_{d_r} \rangle_{g}$ (which in fact belong to $\mbQ$) defined in the $l$-adic \'{e}tale cohomology of  $\overline{\mfM}_{g, r}$ as above  coincide with the usual intersection numbers of  the corresponding psi classes on $\overline{\mfM}_{g,r, \mbC}^\mr{top}$.
\end{rema}
\SSP

Next, let us  define 
 the class $\widehat{\psi}_i  \in \widetilde{H}_{\text{\'{e}t}}^2 (\mfO \mfp_{G,g,r}^{^\mr{Zzz...}}, \overline{\mbQ}_l)$ on $\mfO \mfp_{G,g,r}^{^\mr{Zzz...}}$ to be   the pull-back of  $\psi_i$.
Given  a pair of  nonnegative integers $(g, r)$, an $r$-tuple of nonnegative integers 
   $(d_1, \cdots, d_r)$, and an $r$-tuple of elements $(v_1, \cdots, v_r)$ of $\mcV$,  
we write 
\begin{align}
 \langle \tau_{d_1} (v_1) \cdots \tau_{d_r} (v_r)\rangle_{G, g} \ & \left(= \langle \prod_{i=1}^r \tau_{d_i} (v_i)\rangle_{G, g}\right) 
    := \int_{ [\mfO \mfp^{^\mr{Zzz...}}_{G, g,r}]^\mr{vir} } 
 \prod_{i=1}^r \mr{ev}_i^{*}(v_i) \widehat{\psi}^{d_i}_i  
 \in \overline{\mbQ}_l,
\end{align}
where  $\langle - \rangle_{G, g} := 0$  if either $r =0$ or $2g-2+r \leq 0$ is satisfied.
The invariants  $\langle \tau_{d_1} (v_1) \cdots \tau_{d_r} (v_r)\rangle_{G, g}$ are called 
the {\bf $r$-point correlators}.

\SSP
\bpr \label{T02} 
Let  $(g, r)$, $(d_1, \cdots, d_r)$, and $(v_1, \cdots, v_r)$ be as above.
 Then, the following equality holds:
 \begin{align}
 \langle \tau_{d_1} (v_1), \cdots, \tau_{d_r} (v_r)\rangle_{G, g} = 
    \Lambda_{G, g,r} (\bigotimes_{i=1}^r v_i)
    \langle \tau_{d_1} \cdots \tau_{d_r} \rangle_{g}.
 \end{align}
\epr
\begin{proof}
The assertion follows from the following sequence of equalities:
\begin{align}
& \ \ \  \  \langle \tau_{d_1} (v_1) \cdots \tau_{d_r} (v_r)\rangle_{G, g} \\
 & =  
 (f \circ \pi_{g,r})_*^\mr{hom} \left(\left(\prod_{i=1}^r \mr{ev}_i^{*}(v_i) \widehat{\psi}^{d_i}_i \right) \cap  \mr{cl}^{3g-3+r} ([\mfO \mfp^{^\mr{Zzz...}}_{G,g,r}]^\mr{vir}) \right)
  \notag \\
 & = f_*^\mr{hom} \circ \pi_{g,r}^\mr{hom} \left( \left(\prod_{i=1}^r \widehat{\psi}^{d_i}_i \right)  \cap \left( \left(\prod_{i=1}^r \mr{ev}_i^*(v_i) \right)\cap \mr{cl}^{3g-3+r} ([\mfO \mfp^{^\mr{Zzz...}}_{G,g,r}]^\mr{vir})\right) \right)\notag \\
  & = f_*^\mr{hom} \left( \left( \prod_{i=1}^r \psi^{d_i}_i \right) \cap \pi_{g,r*}^\mr{hom} \left(\left( \prod_{i=1}^r \mr{ev}_i^* (v_i)\right) \cap \mr{cl}^{3g-3+r} ([\mfO \mfp^{^\mr{Zzz...}}_{G,g,r}]^\mr{vir})  \right) \right)   \notag \\
& = f_*^\mr{hom} \left(\left(\prod_{i=1}^r \psi^{d_i}_i \right) \cap \left(\Lambda_{G, g,r}(\bigotimes_{i=1}^r v_i) [ \overline{\mfM}_{g,r}]\right)\right) \notag \\
& = \Lambda_{G, g,r} (\bigotimes_{i=1}^r v_i) f_*^\mr{hom}  \left(\left( \prod_{i=1}^r \psi^{d_i}_i \right)\cap \mr{cl}^{3g-3+r} ([ \overline{\mfM}_{g,r}]) \right) \notag \\
& =  \Lambda_{G, g,r} (\bigotimes_{i=1}^r v_i) \langle \tau_{d_1} \cdots \tau_{d_r} \rangle_g,  \notag
\end{align}
where $f$ denotes the structure morphism $\overline{\mfM}_{g,r} \migi \mr{Spec}(k)$ of $\overline{\mfM}_{g,r}$ and the third equality follows from the projection formula (cf. (\ref{ee436})).
\end{proof}

\LSP
\subsection{The partition function of $G$-do'pers} 
   \label{SS5d}
Let $\hbar$ and $t_{d, \, \pmb{\rho}}$  ($d \in  \mbZ_{\geq 0}, \pmb{\rho}\in \Delta$) be  formal parameters.
Given a basis $\mfe' := \{ e'_{\pmb{\rho}}\}_{\pmb{\rho} \in \Delta}$ of $\mcV$ and  an integer  $g \geq 0$, we set 
\begin{align} \label{ee810}
\Phi_{G, g, \mfe'}    &:= \langle \mr{exp} (\sum_{d \in \mbZ_{\geq 0}, \, \pmb{\rho} \, \in \Delta}\tau_d (e'_{\pmb{\rho}}) t_{d, \, \pmb{\rho}})\rangle_{G, g}\\
& \  = \sum_{r \geq 0} \frac{1}{r!} \sum_{\genfrac{.}{.}{0pt}{}{d_1, \cdots, d_r \geq 0,}{\pmb{\rho}_1, \cdots, \pmb{\rho}_r \in \Delta}} \langle \prod_{i=1}^r \tau_{d_i} (e'_{\pmb{\rho}_i})\rangle_{G, g} \prod_{i=1}^r t_{d_i,  \, \pmb{\rho}_i} \notag
 \\
&  \ =  \sum_{(s_{d,  \pmb{\rho}})_{d, \pmb{\rho}}}  \langle \prod_{\genfrac{.}{.}{0pt}{}{d \geq 0}{\pmb{\rho} \in \Delta}} \tau_{d}(e'_{\pmb{\rho}})^{s_{d, \pmb{\rho}}}\rangle_{G, g} \prod_{\genfrac{.}{.}{0pt}{}{d \geq 0}{\pmb{\rho} \in \Delta}} \frac{t^{s_{d, \pmb{ \rho}}}_{d, \, \pmb{\rho}}}{s_{d,  \pmb{\rho}}!} 
  \left(\in \overline{\mbQ}_l [[\{ t_{d, \, \pmb{\rho}} \}_{d \in \mbZ_{\geq 0}, \,  \pmb{\rho} \in \Delta}]] \right), \notag 
\end{align}
where the sum in the rightmost of this sequence  runs over the set of sequences of nonnegative integers $(s_{d, \, \pmb{\rho}})$ indexed by the elements of $\mbZ_{\geq 0} \times \Delta$ with finitely many nonzero integers.
Also, write
\begin{align} \label{ee811}
\Phi_{G, \mfe'} &:= \sum_{g \geq 0} \Phi_{G, g, \mfe'}  \hbar^{2g-2}  \left(\in \hbar^{-2} \overline{\mbQ}_l [[\hbar ]][[\{ t_{d, \pmb{\rho}} \}_{d \in \mbZ_{\geq 0}, \pmb{\rho} \, \in \Delta}]] \right), \\
Z_{G, \mfe'}  &:= \mr{exp} \left(\Phi_{G, \mfe'} \right)  \left(\in  \overline{\mbQ}_l ((\hbar ))[[\{ t_{d, \pmb{\rho}} \}_{d \in \mbZ_{\geq 0}, \pmb{\rho} \, \in \Delta}]] \right). \notag 
\end{align}
If $\mfe''$ is another basis of $\mcV$, then the change of basis from $\mfe'$ to $\mfe''$ induces naturally an automorphism of  the $\overline{\mbQ}_l ((\hbar ))$-algebra  $\overline{\mbQ}_l ((\hbar ))[[\{ t_{d, \pmb{\rho}} \}_{d \in \mbZ_{\geq 0}, \pmb{\rho} \, \in \Delta}]]$,  by which $Z_{G, \mfe'}$ is  mapped  to $Z_{G, \mfe''}$.
Notice that  $Z_{G, \mfe_\Delta}$ coincides with   $Z_G$   described  in Introduction.

\SSP
\bde \label{DFF061} 
We shall refer to $Z_G \left(= Z_{G, \mfe_\Delta}\right)$ as the {\bf partition function
of $G$-do'pers}.
  \ede
\SSP

Next,  let us fix a canonical base  $\mfe^\dagger := \{ e^\dagger_{\pmb{\rho}} \}_{\pmb{\rho} \in \Delta}$ (with $\eta (e^\dagger_{\pmb{\rho}}, e^\dagger_{\pmb{\rho}}) = \nu_{\pmb{\rho}}$) of the Frobenius algebra corresponding to the $2$d TQFT 
$\Lambda_G := (\mcV, \eta, e_{\pmb{\varepsilon}}, \{ \Lambda_{G, g,r} \}_{g,r})$.
Also,   fix  elements  $\nu_{\pmb{\rho}}^{1/3}$ of $\overline{\mbQ}_l$ with $(\nu_{\pmb{\rho}}^{1/3})^3 = \nu_{\pmb{\rho}}$.
For each $\pmb{\rho} \in \Delta$ and $n \geq -1$,
we shall write
\begin{align}
L_{n}^{(\pmb{\rho})} := &  - \frac{(2n+3)!!}{2^{n+1}} 
(\nu_{\pmb{\rho}}^{1/3})^n
 \frac{\partial}{\partial t_{\pmb{\rho}, n+1}} + \sum_{i=0}
^\infty \frac{(2i+2n+1)!!}{(2i-1)!!2^{n+1}} 
(\nu_{\pmb{\rho}}^{1/3})^n 
t_{\pmb{\rho}, i} \frac{\partial}{\partial t_{\pmb{\rho}, i+n}} \\
& + \frac{\hbar^2}{2} \sum_{i=0}^{n-1}\frac{(2i+1)!! (2n-2i-1)!!}{2^{n+1}} (\nu_{\pmb{\rho}}^{1/3})^{n-3}  \frac{\partial^2}{\partial t_{\pmb{\rho}, i} \partial t_{\pmb{\rho}, n-1-i}} \notag \\
& + \delta_{n, -1} \frac{\hbar^{-2}}{2} 
(\nu_{\pmb{\rho}}^{1/3})^2
 t^2_{\pmb{\rho}, 0} + \delta_{n, 0} \frac{1}{16}.\notag
\end{align}
There operators satisfy  $[L_{n}^{(\pmb{\rho}_1)}, L_{m}^{(\pmb{\rho}_2)}] = (n-m)  \delta_{\pmb{\rho}_1, \pmb{\rho}_2} L_{n+m}^{(\pmb{\rho}_1)}$ for any $n, m \geq -1$ and any $\pmb{\rho}_1, \pmb{\rho}_2 \in \Delta$ (cf. ~\cite[\S\,4.3]{JK}).
Similarly, 
if  $L_n$  ($n \geq -1$) are  the  differential operators defined in Introduction,
then 
 the equality $[L_n, L_m] = (n-m) L_{n+m}$ holds for any $n$, $m \geq -1$.

\SSP
\bt \label{TT01}
\begin{itemize}
\item[(i)]
 Given   $\pmb{\rho} \in \Delta$, $n \geq -1$,
   we obtain the following equality:
 \begin{align}
 L_{n}^{(\pmb{\rho})} Z_{G, \mfe^\dagger} = 0.
 \end{align} 
 Moreover, these equalities for various $(\pmb{\rho}, n)$'s completely determine $Z_G$ (cf. the discussion preceding Definition \ref{DFF061}).
\item[(ii)]
 For any $n \geq -1$, the following equality holds:
 \begin{align}
 L_n Z_{G} = 0.
 \end{align}
\end{itemize}
 \et
\begin{proof}
The assertions follow from  Proposition \ref{T02}  and  ~\cite[Proposition 4.4]{JK} applied to
$\Lambda_G$.
Indeed, 
by a fixed  isomorphism $\overline{\mbQ}_l \isom \mbC$,  Proposition 4.4  in  {\it loc.\,cit.}\,is available because 
the invariants $\langle \tau_1, \cdots, \tau_{d_r}\rangle$ defined in (\ref{ee150})  coincide with the usual intersection numbers of the corresponding psi classes on the module space of complex curves (cf. Remark \ref{REm2}).
    In particular, the Witten-Kontsevich theorem, which is used in the proof of that proposition, also holds   for the psi classes in the $l$-adic \'{e}tale cohomology groups of $\overline{\mfM}_{g, r}$'s.
\end{proof}
\SSP

Let $(g,r)$ be a pair of nonnegative integers, $(d_1, \cdots, d_r)$ an $r$-tuple of nonnegative integers, and $(v_1, \cdots, v_r)$ an $r$-tuple of elements  of $\mcV$. 
Then, we shall write
\begin{align}
\langle \langle \tau_{d_1} (v_1) \cdots \tau_{d_r} (v_r) \rangle \rangle_{G, g} &:= \langle \tau_{d_1} (v_1) \cdots \tau_{d_r} (v_r) \mr{exp} (\sum_{d \in \mbZ_{\geq 0}, \, \pmb{\rho} \, \in   \Delta} \tau_d (e_{\pmb{\rho}}) t_{d,  \pmb{\rho}}) \rangle_{G, g}, 
\end{align}
where  $\langle \langle - \rangle \rangle_{G, g} := 0$ if either $r =0$ or $2g-2+r \leq 0$ is satisfied.
Also,  write 
\begin{align}
\langle \langle \tau_{d_1} (v_1) \cdots \tau_{d_r} (v_r) \rangle \rangle_G &:=
\sum_{g \geq 0} \hbar^{2g-2} \langle \langle \tau_{d_1} (v_1) \cdots \tau_{d_r} (v_r) \rangle \rangle_{G, g}.
\end{align}

Then, by ~\cite[Proposition 4.6]{JK}, the following proposition holds.

\SSP
\bt \label{TT01rrr} 
 For any  $d \in \mbZ_{\geq 0}$, $v \in \mcV$, and $\pmb{\rho}_1, \pmb{\rho}_2, \pmb{\rho}_3, \pmb{\rho}_4 \in \Delta$, the following equality holds:
\begin{align} \label{ee780}
& \  \frac{2d+1}{\hbar^2} \langle \langle \tau_d (v) \tau_0 (e_{\pmb{\rho}_1}) \tau_0 (e_{\pmb{\rho}_2}) \rangle \rangle_G \eta^{e_{\pmb{\rho}_1} e_{\pmb{\rho}_2}}  \\
= & \   \langle \langle \tau_{d-1} (v) \tau_0 (e_{\pmb{\rho}_1}) \rangle  \rangle_G \eta^{e_{\pmb{\rho}_1} e_{\pmb{\rho}_2}} \langle \langle \tau_0 (e_{\pmb{\rho}_2}) \tau_0 (e_{\pmb{\rho}_3}) \tau_0 (e_{\pmb{\rho}_4}) \rangle \rangle_G \eta^{e_{\pmb{\rho}_3} e_{\pmb{\rho}_4}}  \notag \\
+ &  \ 2 \langle \langle \tau_{d-1} (v) \tau_0 (e_{\pmb{\rho}_1}) \tau_0 (e_{\pmb{\rho}_3}) \rangle \rangle_G \eta^{e_{\pmb{\rho}_1} e_{\pmb{\rho}_2}} \eta^{e_{\pmb{\rho}_3} e_{\pmb{\rho}_4}} \langle \langle \tau_0 (e_{\pmb{\rho}_2}) \tau_0 (e_{\pmb{\rho}_4}) \rangle \rangle_G \eta^{e_{\pmb{\rho}_3} e_{\pmb{\rho}_4}} \notag   \\
+ &  \  \frac{1}{4} \langle \langle \tau_{d-1} (v) \tau_0 (e_{\pmb{\rho}_1}) \tau_0 (e_{\pmb{\rho}_2}) \tau_0 (e_{\pmb{\rho}_3}) \tau_0 (e_{\pmb{\rho}_4})    \rangle \rangle_G\eta^{e_{\pmb{\rho}_1} e_{\pmb{\rho}_2}} \eta^{e_{\pmb{\rho}_3} e_{\pmb{\rho}_4}}.    \notag
\end{align}
Moreover, equation (\ref{ee780}) and the equalities  $L^{(\pmb{\rho})}_{-1} Z_{G, \mfe^\dagger}= 0$ ($\pmb{\rho} \in \Delta$) resulting from Theorem \ref{TT01} completely determine $\Phi_{G, \mfe^\dagger}$.
 \et

\LSP
\subsection{The $2$d TQFT for $\mr{PGL}_2$-do'pers} \label{SSjj5d}
In this  subsection,  we   compare 
   the  $2$d TQFT  for $\mr{PGL}_2$-do'pers
  with 
the $\mfs \mfl_2 (\mbC)$-WZW (= Wess-Zumino-Witten) conformal field theory.

 Write $(\mcU, \eta^\mcU)$ for   the fusion ring associated with the $\mfs \mfl_2 (\mbC)$-WZW conformal field theory at level $p-2$ (i.e., ``$\mcR_{p-2}(\mfs \mfl_2)$" in the notation of ~\cite{Beau}), and 
 $Z^{\mr{WZW}, p-2}_{\mfs \mfl_2}$ for  the partition function associated with the corresponding $2$d TQFT.
 
 First, let us recall the structure of the $\mbC$-algebra $\mcU$.
Denote by $\alpha$ the
unique root in  a Borel subgroup of $\mr{PGL}_2$ 
  and
  by $(-, -)$ the Killing form on $\mfs \mfl_2 (\mbC)$ normalized as $(H_\alpha, H_\alpha) = 2$, where $H_{\alpha}$ denotes the coroot  (considered as an element of $\mft$ via differentiation)  associated to $\alpha$.
Also, denote by $P_{p-2}$ the set of dominant weights $\lambda$ of $\mfs \mfl_2 (\mbC)$ with $0 \leq \lambda (H_\alpha) \leq p-2$.
 That is to say,  we have
 $P_{p-2} = 
\{ 0, \frac{1}{2} \alpha, \alpha, \cdots, \frac{p-2}{2} \alpha \}$.
We shall identify $P_{p-2}$ with the set $\left\{ 0, \frac{1}{2}, \cdots, \frac{p-2}{2} \right\} \left(\subseteq \frac{1}{2} \mbZ\right)$ via the correspondence $j \alpha \leftrightarrow j$.
 Then,   the underlying $\mbC$-vector space of $\mcU $ is isomorphic to   the direct sum   $\bigoplus_{j  \in P_{p-2}} \mbC e^\mcU_{j}$ with basis $\{e^\mcU_j \}_{j \in P_{p-2}}$ indexed by $P_{p-2}$.
(In particular, $Z^{\mr{WZW}, p-2}_{\mfs \mfl_2}$ may be considered as an element of $\mbC ((\hbar ))[[ \{ t_{d, m}\}_{d \geq0, m  \in P_{p-2}}]]$.)
The multiplication ``$\times$" in $\mcU$ is determined by the following conditions (cf. ~\cite[Lemma 4.2]{Beau}):
\begin{itemize}
\item
 If
   $\{ N_{a, b, c} \in \mbC \, | \, a, b, c \in P_{p-2}\}$ denotes  a collection defined in such a way that 
$e_a^\mcU \times e_b^\mcU = \sum_{c \in P_{p-2}} N_{a, b, c} e_c^\mcU$, then
 $N_{a, b, c}  \in \{0, 1\}$ for any triple $(a ,b,c)$;
 \item
 The equality    $N_{a, b,c} = 1$ holds if and only if 
$(a, b,c)$ satisfies\footnote{In ~\cite[Eq.\,(1032)]{Wak}, this condition is incorrectly described.} 
 \begin{align} \label{ee577}
 a + b + c \in \mbZ, \hspace{5mm} a + b + c \leq p-2, \hspace{5mm} \text{and} \hspace{5mm} |b-c| \leq a \leq b +c. 
 \end{align}
 \end{itemize}
 
 On the other hand, 
 let 
 $(\mcV, \eta^\mcV)$  be the Frobenius  algebra associated with $\Lambda_{\mr{PGL}_2}$.
 We regard $(\mcV, \eta^\mcV)$ as a $\mbC$-algebra by  fixing an isomorphism $\overline{\mbQ}_l \isom \mbC$.
 We shall  write 
  \begin{align}
  P^\mbZ_{p-2} := \left\{ m \in \mbZ \ \bigg| \ 0 \leq j \leq \frac{p-2}{2} \right\}   \left(
  = P_{p-2} \cap \mbZ \right).
  \end{align}
Then, the finite set $\Delta$ (cf. (\ref{ee142})) in the case of $G = \mr{PGL}_2$ (where the maximal torus $T$ is taken to be the subgroup consisting of the images of invertible diagonal matrices via the quotient $\mr{GL}_2 \migisurj \mr{PGL}_2$) may be identified with $P^\mbZ_{p-2}$ via the  bijection
\begin{align} \label{ER52}
P^\mbZ_{p-2} \isom \Delta
\end{align}
given by assigning, to each $a \in P^\mbZ_{p-2}$, the element of $\mft_{\mr{reg}}^F/W$ represented by the diagonal matrix 
\begin{align}
\begin{pmatrix}  \overline{(2a+1)/2} & 0 \\ 0 & - \overline{(2a+1)/2}\end{pmatrix}.
\end{align}
 Then, it follows from ~\cite[Introduction, Theorem 1.3]{Mzk2}) that, under the identification  $P^\mbZ_{p-2} = \Delta$ given by (\ref{ER52}),
 $\mcV$ is isomorphic to 
  the $\mbC$-vector space 
   $\bigoplus_{m \in P^\mbZ_{p-2}} \mbC e^\mcV_{m}$
   (where 
 $\{ e^\mcV_m \}_{m \in P^\mbZ_{p-2}}$ is a basis 
  indexed by  $P^\mbZ_{p-2}$) with 
  multiplication characterized uniquely by
  the condition that the assignment   $e^\mcV_{m} \mapsto e^\mcU_m$
  defines a $\mbC$-algebra homomorphism
  \begin{align}
  \mr{incl} : \mcV \migi \mcU.
  \end{align}

Just as in the case of the pseudo-fusion ring for dormant $\mr{PGL}_2$-opers discussed in  ~\cite[\S\,7.8.2]{Wak5}, we can describe   the characters $\mcV \migi \mbC$ of $\mcV$ from those of $\mcU$ via  this homomorphism $\mr{incl}$.
In particular, by ~\cite[Proposition 6.3]{Beau}, the explicit knowledge  of the characters allows us to perform some  computations that we need in the ring $\mcU$, e.g., the formula displayed in (\ref{ee80ee2}).
Also, we obtain the following assertion regarding the partition function $Z_{\mr{PGL}_2}$:
 
\SSP
\bt \label{TT11}  
 Let us consider the surjective morphism of $\mbC$-algebras
 \begin{align}
 \theta : \mbC (( \hbar )) [[\{ t_{d, m}\}_{d \geq0, m  \in P_{p-2}}]] \migisurj \mbC (( \hbar )) [[ \{ t_{d, m}\}_{d \geq0, m  \in P^\mbZ_{p-2}}]]
 \end{align}
 given by $\hbar \mapsto \frac{\hbar}{2}$, $t_{d, m} \mapsto t_{d, m}$ for any $m \in P_{p-2}^\mbZ$, and $t_{d, m} \mapsto 0$ for any $m \in P_{p-2} \setminus P^\mbZ_{p-2}$.
 Then, the following equality holds:
 \begin{align}
 Z_{\mr{PGL}_2} = \theta (Z^{\mr{WZW}, p-2}_{\mfs \mfl_2}).
 \end{align}
  \et
\begin{proof}
The assertion follows from the fact that
 if  $\Lambda_{g,r}^{\hspace{-0.7mm} \mr{WZW}, p-2}$ ($g$, $r \geq 0$, $2g-2+r >0$) denote the correlators of the $2$d  TQFT corresponding to  the Frobenius algebra $(\mcU, \eta^\mcU)$, then 
 the equality 
 \begin{align}
 \Lambda_{g,r}^{\hspace{-0.7mm} \mr{WZW}, p-2} (\bigotimes_{i=1}^r e^\mcU_{n_i}) =  2^{g-1} \cdot \Lambda_{\mr{PGL}_2, g, r} (\bigotimes_{i=1}^r e^\mcV_{n_i})
 \end{align}
  holds for any $n_1, \cdots, n_r \in P^\mbZ_{p-2}$ (cf. (\ref{ee80ee2}) and ~\cite[Corollary 9.8]{Beau}).
\end{proof}
\LSP

\end{document}